\documentstyle[twoside,12pt]{article}

\newwrite\citeno
\immediate\openout\citeno=cite
\def\site#1{{\cite{#1}}
    \write\citeno{#1  \thepage}}
\newcommand{\indexentry}[2]{\par\noindent #1 \hrulefill\  #2}
\def\dex#1{{\it #1}\index{#1}}
\def\sdex#1{\index{#1}}
\newtheorem{theorem}{Theorem}[section]
\newtheorem{lemma}[theorem]{Lemma}
\newtheorem{corollary}[theorem]{Corollary}
\newtheorem{question}[theorem]{{\bf Question}}
\newcommand{\sector}[1]{\newpage\section{#1}\markright{Descriptive Set Theory
and Forcing}}
\newcommand{\claim}{\par\bigskip\noindent {\bf Claim: }}
\newcommand{\proof}{\par\noindent proof:\par}

\newcommand{\val}[1]{\left[\kern-2pt|  \; #1\;  |\kern-2pt\right]}
\mathchardef\wtilde="0365
\newcommand{\bpi}[2]{\vtop{\hbox{${\bf\Pi}^{#1}_{#2}$} \kern-4pt
       \hbox{\kern1pt$\wtilde$}\kern-6pt}}
\newcommand{\bsig}[2]{\vtop{\hbox{${\bf\Sigma}^{#1}_{#2}$} \kern-4pt
       \hbox{\kern1pt$\wtilde$}\kern-6pt}}
\newcommand{\bdel}[2]{\vtop{\hbox{${\bf\Delta}^{#1}_{#2}$} \kern-4pt
       \hbox{\kern1pt$\wtilde$}\kern-6pt}}
\newcommand{\name}[1]{ \stackrel{\circ}{#1} }
\newcommand{\finseq}{\omega^{<\omega} }
\newcommand{\irrationals}{\Bbb P}
\newcommand{\integers}{\Bbb Z}
\newcommand{\cantorset}{\cantorspace}
\newcommand{\comp}[1]{\sim{#1}\-}
\newcommand{\borel}{{\rm Borel}}
\newcommand{\ord}{{\rm ord}}
\newcommand{\rank}{{\rm rank}}
\newcommand{\st}{{\rm stone}}
\newcommand{\dom}{{\rm domain}}
\newcommand{\base}{{\cal B}}
\newcommand{\ordinals}{{\rm OR}}
\newcommand{\tzero}{T^0}
\newcommand{\tone}{T^{>0}}
\newcommand{\bool}{{\Bbb B}}
\newcommand{\emptyseq}{\langle\rangle}
\newcommand{\disj}{\mathop{\vee}}
\newcommand{\MA}{{\rm MA}}
\newcommand{\dense}{{\cal D}}
\newcommand{\ideal}{{\cal I}}
\newcommand{\famly}{{\cal F}}
\newcommand{\fbool}{{\Bbb F}}
\newcommand{\posetq}{{\Bbb Q}}

\newcommand{\interior}{{\rm int}}
\newcommand{\cl}{{\rm cl}}
\newcommand{\meager}{{\rm meager}}
\newcommand{\fin}{{\rm FIN}}
\newcommand{\cof}{{\rm cof}}
\newcommand{\cov}{{\rm cov}}
\newcommand{\sleq}{\preceq}
\newcommand{\embed}{\preceq}
\newcommand{\strictembed}{\prec}
\newcommand{\fleq}{\unlhd}
\newcommand{\fless}{\lhd}
\newcommand{\rmiff}{\mbox{ iff }}
\newcommand{\rmor}{\mbox{ or }}
\newcommand{\rmand}{\mbox{ and }}
\newcommand{\diam}{{\rm diam}}
\newcommand{\lin}{{\rm line}}
\newcommand{\notE}{\not\!\!E}
\newcommand{\hyp}{{\rm hyp}}
\newcommand{\semihyp}{{\rm semihyp}}
\newcommand{\arrow}{\rightarrow}
\newcommand{\bairespace}{\omega^{\omega}}
\newcommand{\cantorspace}{2^{\omega}}
\newcommand{\concat}{{\;\hat{ }\;}}
\newcommand{\continuum}{{\goth c}}
\newcommand{\cross}{\times}
\newcommand{\forces}{\mid\vdash}
\newcommand{\infsets}{[\omega]^{\omega}}
\newcommand{\Intersect}{\bigcap}
\newcommand{\intersect}{\cap}
\newcommand{\poset}{{\Bbb P}}
\newcommand{\rationals}{{\Bbb Q}}
\newcommand{\reals}{{\Bbb R}}

\newcommand{\Union}{\bigcup}
\newcommand{\union}{\cup}
\newcommand{\Conj}{\mathop{\bigwedge\kern -6pt\bigwedge}}
\newcommand{\conj}{\mathop{\wedge}}
\newcommand{\eq}{\approx}
\newcommand{\implies}{\rightarrow}
\newcommand{\isom}{\simeq}

\newcommand{\prm}{}%%%
\newcommand{\res}{\upharpoonright}

\newcommand{\qed}{\nopagebreak\par\noindent\nopagebreak$\blacksquare$\par}
\input amssym.def
\input amssym
\makeindex

\begin{document}
\pagestyle{empty}

\vskip 3in

\begin{center}
 \Huge
 Descriptive Set Theory
\end{center}

\begin{center}
 \Huge
  and
\end{center}

\begin{center}
 \Huge
 Forcing:
\end{center}

\begin{center}
 {\Large How to prove theorems about Borel sets\\ }
 {\Large the hard way. }
\end{center}

\vskip 3in

\begin{flushright}
  Arnold W. Miller\\
  Department of Mathematics\\
  480 Lincoln Dr.\\
  Van Vleck Hall\\
  University of Wisconsin\\
  Madison, WI. 53706\\
  miller@math.wisc.edu\\
  January 1994
\end{flushright}

\newpage

\begin{center}
Note to the readers
\end{center}

\bigskip

Departing from the usual
author's statement-I would like to say that I am not responsible
for any of the mistakes in this document.  Any mistakes here are
the responsibility of the reader.  If anybody wants to point out
a mistake to me, I promise to respond by saying ``but you know what
I meant to say, don't you?"

\bigskip

These are lecture notes from a course I gave at the University of
Wisconsin during the Spring semester of 1993.  Some knowledge of
forcing is assumed as well as a modicum of elementary
Mathematical Logic, for example, the Lowenheim-Skolem Theorem.

Many of the classical theorems of descriptive set theory are
presented ``just-in-time'' for when they are needed.  Questions
like ``Who proved what?'' always interest me, so I have included
my best guess here.  Hopefully, I have managed to offend a large
number of mathematicians.

\bigskip
The  results in section \ref{cohen} and
\ref{random} are new and answer a questions from my thesis.
I have also included (without permission)
an unpublished result of Fremlin (Theorem \ref{frem}).  I think
the proof given of Louveau's Theorem \ref{louv} is also a little different.
\footnote{In a randomly infinite Universe, any event occurring here and
   now with finite probability must be occurring simultaneously at
   an infinite number of other sites in the Universe. It is hard to
   evaluate this idea any further, but one thing is certain: if it
   is true then it is certainly not original!--
    The Anthropic Cosmological
    Principle, by John Barrow and Frank Tipler.}

\newpage
% \input grading
% \newpage
\tableofcontents
\newpage

\pagestyle{headings}

\sector{What are the reals, anyway?}

Definitions.
Let $\omega=\{0,1,\ldots\}$ and let $\bairespace$ \sdex{$\bairespace$}
(\dex{Baire space}) be the
set of functions from $\omega$ to $\omega$. Let $\finseq$ \sdex{$\finseq$}
be the set
of all finite sequences of elements of $\omega$.
$|s|$\sdex{$|s|$}  is the length of $s$, $\langle\rangle$ is the empty sequence,
and for $s\in\finseq$  and $n\in\omega$ let $s\concat n$\sdex{$s?concat n$}
 denote the
sequence which starts out with $s$ and has one more element $n$ concatenated
onto the end.
The basic open sets
of $\bairespace$ are the sets of the form: \sdex{$[s]$}
$$[s] =\{x\in\bairespace : s\subseteq x\}$$
for $s\in\finseq$.  A subset of $\bairespace$ is open iff it is the union
of basic open subsets. It is \dex{separable} (has a countable dense subset)
since it is \dex{second countable} (has a countable basis).
The following defines a complete metric on
$\bairespace$:
$$d(x,y)=\left\{\begin{array}{ll}
               0   & \mbox{if $x=y$}\\
      {1\over n+1} & \mbox{ if $x\res n = y\res n$ and $x(n)\not= y(n)$}
\end{array}\right. $$
\dex{Cantor space} $\cantorspace$\sdex{$\cantorspace$}
 is the subspace  of $\bairespace$ consisting
of all functions from $\omega$ to $2=\{0,1\}$.  It is compact.

\begin{theorem} \label{eul} (Baire \site{baire})
  $\bairespace$ is homeomorphic to the irrationals $\irrationals$.
\end{theorem}
\proof
First replace $\omega$ by the integers $\integers$.  We will construct
a mapping from $\integers^\omega$ to $\irrationals$.
Enumerate the rationals $\rationals=\{q_n: n\in\omega\}$.
Inductively construct a sequence of open intervals
$\langle I_s : s\in \integers^{<\omega}\rangle$ satisfying the
following:

\begin{enumerate}

  \item $I_{\langle\rangle}=\reals$,
  and for $s\not= \langle\rangle$ each $I_s$ is a nontrivial open interval
  in $\reals$ with rational endpoints,

  \item for every $s\in \integers^{<\omega}$ and
  $n\in \integers\;\;\; I_{s\concat n}\subseteq I_s$,

  \item the  right end point of $I_{s\concat n}$ is the left end point of
  $I_{s\concat n+1}$,

  \item  $\{I_{s\concat n}: n\in \integers\}$ covers
  all of $I_s$ except for their endpoints,

  \item the length of $I_s$ is less than $1\over |s|$ for
  $s\not= \langle\rangle$, and

  \item the $n^{th}$ rational $q_n$ is an endpoint of
   $I_t$ for some $|t|\leq n+1$.

\end{enumerate}
Define the function $f:\integers^{\omega}\arrow\irrationals$ as
follows.  Given $x\in\integers^{\omega}$ the set
$$\Intersect_{n\in\omega}I_{x\res n}$$
must consist of a singleton
irrational.  It is nonempty because
$$\mbox{closure}(I_{x\res n+1})\subseteq I_{x\res n}.$$
It is a singleton
because their diameters shrink to zero.

So we can define $f$ by
$$\{f(x)\}=\Intersect_{n\in\omega}I_{x\res n}.$$
The function $f$ is one-to-one because if $s$ and $t$ are incomparable
then $I_s$ and $I_t$ are disjoint.  It is onto since for every
$u\in\irrationals$ and $n\in\omega$ there is a unique $s$ of
length $n$ with $u\in I_s$. It is a homeomorphism because
$$f([s])=I_s\intersect\irrationals$$
and the sets of the form $I_s\intersect\irrationals$ form a basis
for $\irrationals$.
\qed

Note that the map given is also an order isomorphism from
$\integers^\omega$ with the lexicographical order to $\irrationals$
with it's usual order.

We can identify $\cantorspace$ with $P(\omega)$, the set of all
subsets of $\omega$, by identifying a subset with its characteristic
function.  Let $F=\{x\in\cantorspace: \forall^{\infty} n\; x(n)=0\}$
(the quantifier $\forall^{\infty}$ stands for ``for all but finitely
many $n$'').  $F$ corresponds to the finite sets and so
$\cantorspace\setminus F$ corresponds to the infinite subsets of
$\omega$ which we write as $\infsets$.\sdex{$\infsets$}

\begin{theorem}
  $\bairespace$ is homeomorphic to $\infsets$.
\end{theorem}

\proof
Let $f\in\bairespace$ and define $F(f)\in\cantorspace$ to be
the sequence of 0's and 1's determined by:
$$F(f)=0^{f(0)}\concat 1 \concat 0^{f(1)}\concat 1\concat 0^{f(2)}
\concat 1 \concat\cdots$$
where $0^{f(n)}$ refers to a string of length $f(n)$ of zeros.
The function $F$ is a one-to-one onto map from $\bairespace$ to
$\cantorspace\setminus F$. It is a homeomorphism because
$F([s])=[t]$ where
$t= 0^{s(0)}\concat 1 \concat 0^{s(1)}\concat 1\concat 0^{s(2)}
\concat 1 \concat\cdots\concat 0^{s(n)}\concat 1$ where $|s|=n+1$.
Note that  sets of the form $[t]$ where $t$ is a finite sequence ending
in a one form a basis for $\cantorspace\setminus F$.
\qed

I wonder why $\bairespace$ is called Baire space?
The earliest mention
of  this I have seen is in Sierpi\'{n}ski \site{sierp1928}
where he refers to $\bairespace$ as the 0-dimensional
space of Baire.  Sierpi\'{n}ski also says that Frechet was the first
to describe the  metric $d$ given above.  Unfortunately,
Sierpi\'{n}ski \site{sierp1928} gives very few references.

The classical
proof of Theorem \ref{eul}
is to use ``continued fractions'' to get the correspondence.
Euler \site{euler} proved that every
rational number gives rise to a finite continued fraction and
every irrational number gives rise to an infinite continued fraction.
Brezinski \site{brez} has more on the history of continued fractions.

My proof of Theorem \ref{eul} allows me to remain
blissfully ignorant\footnote{ It is impossible for a man to learn
what he thinks he already knows.-Epictetus }
of even the elementary theory of continued fractions.

Cantor space, $\cantorspace$, is clearly named so because it is
homeomorphic to Cantor's middle two thirds set.

\newpage\part{On the length of Borel hierarchies}
\section{Borel Hierarchy}

\bigskip

Definitions.
For $X$ a topological space define $\bsig01$  to
be the open
subsets of $X$. For $\alpha>1$ define
$A\in\bsig0\alpha$ iff there exists a sequence $\langle B_n:n\in\omega\rangle$

with each $B_n\in\bsig0{\beta_n}$ for some $\beta_n<\alpha$ such that
$$A=\Union_{n\in\omega}\comp{B_n}$$
where $\comp{B}$ is the complement of $B$ in $X$, i.e.,
$\comp{B}=X\setminus B$.  Define $\bpi0\alpha=\{\comp{B}: B\in\bsig0\alpha\}$
and $\bdel0\alpha=\bsig0\alpha\intersect\bpi0\alpha$.  The Borel
subsets of $X$ are defined
by $\borel(X)=\Union_{\alpha<\omega_1}\bsig0\alpha(X)$. It is clearly
the smallest family of sets containing the open subsets of $X$ and
closed under countable unions and complementation.

\sdex{$\comp{B}$}
\sdex{$?bpi0\alpha$}
\sdex{$?bdel0\alpha$}
\sdex{$\borel(X)$}
\sdex{$?bsig0\alpha$}

\begin{theorem}
$\bsig0\alpha$ is closed under countable unions and finite intersections,
$\bpi0\alpha$ is closed under countable intersections and finite unions,
and $\bdel0\alpha$ is closed under finite intersections, finite
 unions, and complements.
\end{theorem}
\proof
That $\bsig0\alpha$ is closed under countable unions is clear from its
definition. It follows from
DeMorgan's laws  by taking complements that
$\bpi0\alpha$ is closed under countable intersections.
Since
$$(\Union_{n\in\omega} P_n)\intersect (\Union_{n\in\omega} Q_n)
=\Union_{n,m\in\omega} (P_n \intersect Q_m) $$
$\bsig0\alpha$ is closed under finite intersections.  It follows by
DeMorgan's laws that
$\bpi0\alpha$ is closed under finite unions.
$\bdel0\alpha$ is closed under finite intersections, finite
unions, and complements since it is the intersection of the
two classes.
\qed

\begin{theorem} \label{map}
If $f:X\arrow Y$ is continuous and $A\in\bsig0\alpha(Y)$, then
$f^{-1}(A)\in \bsig0\alpha(X)$.
\end{theorem}
This is an easy induction since it is true for open sets
($\bsig01$) and $f^{-1}$ passes over complements and unions.
\qed
Theorem \ref{map} is also, of course,
true for $\bpi0\alpha$ or $\bdel0\alpha$
in place of $\bsig0\alpha$.

\begin{theorem} \label{subsp}
  Suppose $X$ is a subspace of $Y$, then
  $$\bsig0\alpha(X)=\{A\intersect X: A \in \bsig0\alpha(Y)\}.$$
\end{theorem}
\proof
For $\bsig01$ it follows from the definition of subspace.
For $\alpha>1$ it is an easy induction.
\qed

Theorem \ref{subsp} is
true for $\bpi0\alpha$
in place of $\bsig0\alpha$, but not in general for $\bdel0\alpha$.

The class of sets
$\bsig02$ is also referred to as $F_\sigma$ and the
class $\bpi02$ as $G_\delta$.

\sdex{$F_\sigma$}
\sdex{$G_\delta$}

\begin{theorem} For $X$ a topological space
\begin{enumerate}
  \item $\bpi0\alpha(X)\subseteq \bsig0{\alpha+1}(X)$,
  \item $\bsig0\alpha(X)\subseteq \bpi0{\alpha+1}(X)$, and
  \item   if $\bpi01(X)\subseteq \bpi02(X)$ (i.e., closed sets are
          $G_\delta$), then
  \begin{enumerate}
   \item $\bpi0\alpha(X)\subseteq \bpi0{\alpha+1}(X)$,
   \item $\bsig0\alpha(X)\subseteq \bsig0{\alpha+1}(X)$, and hence
   \item $\bpi0\alpha(X)\union \bsig0\alpha(X)\subseteq  \bdel0{\alpha+1}(X)$.
  \end{enumerate}
\end{enumerate}
\end{theorem}
\proof
Induction on $\alpha$.
\qed
In metric spaces closed sets are $G_\delta$, since
$$C=\Intersect_{n\in\omega}\{x: \exists y \in C\;\;d(x,y)<{1\over n+1}\}$$
for $C$ a closed set.

The assumption that closed sets are $G_\delta$ is necessary since
if
$$X=\omega_1+1$$
with the order topology, then the closed set
consisting of the singleton point $\{\omega_1\}$ is not $G_\delta$;
in fact, it is not in the \dex{$\sigma$-ring} generated by the open sets
(the smallest family containing the open sets and closed under
countable intersections and countable unions).

Williard \site{will} gives an example which
is a second countable Hausdorff space.
Let $X\subseteq\cantorspace$ be any nonBorel set. Let $\cantorspace_*$
be the space $\cantorspace$ with the smallest topology containing
the usual topology and $X$ as an open set.  The family of
all sets of the form $(B\intersect X)\union C$ where $B,C$ are
(ordinary) Borel subsets of $\cantorspace$ is
the $\sigma$-ring generated by the open subsets of $\cantorspace_*$,
because:
$$\Intersect_n(B_n\intersect X)\union C_n=
(\;(\Intersect_nB_n\union C_n)\intersect X)\union \Intersect_n C_n$$
$$\Union_n (B_n\intersect X)\union C_n=
((\Union_n B_n)\intersect X)\union \Union_n C_n.$$
Note that $\comp{X}$ is not in this $\sigma$-ring.

\begin{theorem}\label{leb}
  (Lebesgue \site{leb}) For every $\alpha$ with $1\leq\alpha<\omega_1$
  $\bsig0\alpha(\cantorspace)\not=\bpi0\alpha(\cantorspace)$.
\end{theorem}

The proof of this is a diagonalization argument applied to a universal set.
\sdex{universal set}
We will need the following two lemmas.
\begin{lemma}\label{univ}
  Suppose $X$ is second countable (i.e. has a countable base), then
  for every $\alpha$ with $1\leq\alpha<\omega_1$ there exists
  a \dex{universal}' $\bsig0\alpha$ set $U\subseteq \cantorspace\cross X$,
  i.e., a set $U$ which is $\bsig0\alpha(\cantorspace\cross X)$ such
  that for every $A\in\bsig0\alpha(X)$ there exists $x\in\cantorspace$
  such that $A=U_x$ where $U_x=\{y\in X: (x,y)\in U\}$.
\end{lemma}
\proof
The proof is by induction on $\alpha$.  Let
$\{B_n:n\in\omega\}$ be a countable base for $X$.
For $\alpha=1$ let
$$U=\{(x,y): \exists n \;(x(n)=1 \conj y\in B_n)\}=
\Union_n(\{x:x(n)=1\}\cross B_n).$$
For $\alpha>1$ let $\beta_n$ be a sequence which sups up to $\alpha$
if $\alpha$ a limit, or equals $\alpha-1$ if $\alpha$ is a successor.
 Let $U_n$ be a universal $\bsig0{\beta_n}$ set.
Let $$\langle n,m \rangle=2^n(2m+1)-1$$ be the usual pairing function
which gives a recursive bijection between $\omega^2$ and $\omega$.
For any $n$ the map $g_n:\cantorspace\cross X\arrow \cantorspace\cross X$
is defined by $(x,y)\mapsto (x_n,y)$ where $x_n(m)=x(\langle n,m\rangle)$.
This map is continuous so if we define $U_n^*=g_n^{-1}(U_n)$, then
$U_n^*$ is $\bsig0{\beta_n}$, and because the map $x\mapsto x_n$ is onto
it is also a universal $\bsig0{\beta_n}$ set.  Now define $U$ by:
$$U=\Union_n \comp{U_n^*}.$$
$U$ is universal for $\bsig0\alpha$ because given any sequence
$B_n\in\bsig0{\beta_n}$ for $n\in\omega$ there exists $x\in\cantorspace$ such
that for every $n\in\omega$ we have that $B_n=(U_n^*)_{x}=(U_n)_{x_n}$ (this
is because the map $x\mapsto \langle x_n:n<\omega\rangle$ takes
$\cantorspace$ onto $({\cantorspace})^\omega$.) But then
$$U_x=(\Union_n \comp{U_n^*})_x= \Union_n \comp{(U_n^*)_x}
=\Union_n\comp(B_n).$$
\qed

\par\noindent Proof of Theorem \ref{leb}:

Let $U\subseteq\cantorspace\cross\cantorspace$ be a universal
$\bsig0\alpha$ set.  Let
$$D=\{x:\langle x,x\rangle\in U\}.$$
$D$ is the continuous preimage of $U$ under the map
$x\mapsto \langle x,x\rangle$, so it is $\bsig0\alpha$, but
it cannot be $\bpi0\alpha$ because if it were, then there would
be $x\in\cantorspace$ with $\comp{D}=U_x$ and then
$x\in{D}$ iff $\langle x,x\rangle\in U$ iff $x\in U_x$ iff
$x\in \comp{D}$.
\qed
\sdex{$\ord(X)$}
Define $\ord(X)$ to be the least $\alpha$ such that
$\borel(X)=\bsig0\alpha(X)$. Lebesgue's theorem says
that $\ord(X)=\omega_1$.  Note that $\ord(X)=1$ if
$X$ is a discrete space and that $\ord(\rationals)=2$.

\begin{corollary} For any space $X$ which contains a homeomorphic
copy of $\cantorspace$ (i.e., a perfect set) we have that
$\ord(X)=\omega_1$, consequently
$\bairespace$, $\reals$, and any uncountable complete separable
metric space have $\ord=\omega_1$.
\end{corollary}
\proof
If the Borel hierarchy on $X$ collapses, then by Theorem \ref{subsp}
it also collapses on all subspaces of $X$.
Every uncountable complete separable
metric space contains a \dex{perfect set} (homeomorphic copy of $\cantorspace$).
To see this suppose $X$ is an uncountable complete separable
metric space.
Construct a family of open sets
$\langle U_s: s\in 2^{<\omega}\rangle$ such that
\begin{enumerate}
  \item $U_s$ is uncountable,
  \item $\cl(U_{s\concat 0})\intersect \cl(U_{s\concat 1})=\emptyset$,
  \item $\cl(U_{s\concat i})\subseteq U_s$ for i=0,1, and
  \item diameter of $U_s$ less than $1/|s|$
\end{enumerate}
Then the map $f:\cantorset\arrow X$ defined so that
$$\{f(x)\}=\Intersect_{n\in\omega}U_{x\res n}$$
gives an embedding of $\cantorset$ into $X$.
\qed

Lebesgue \site{leb} used universal functions instead of sets, but
the proof is much the same.  Corollary \ref{univdel}
of Louveau's Theorem shows that
there can be no  Borel set which is universal for all
$\bdel0\alpha$ sets.  Miller \site{ctblmod} contains
examples from model theory of Borel sets of arbitrary high
rank.

The notation ${\bf \Sigma^0_\alpha}, \Pi^0_\beta$ was first popularized
by Addison \site{addison}.  I don't know if the ``bold face''
and ``light face''
notation is such a good idea, some copy machines wipe it out.  Consequently,
I use
$$\bsig0\alpha$$
which is blackboard boldface.

\sector{Abstract Borel hierarchies}

Suppose $F\subseteq P(X)$ is a family of sets.
Most of the time we would like to think of $F$ as a countable
\dex{field of sets}
(i.e. closed under complements
and finite intersections) and so analogous to the family of
clopen subsets of some space.

We define the classes $\bpi0\alpha(F)$
analogously. Let $\bpi00(F)=F$ and for every $\alpha>0$ define
$A\in\bpi0\alpha(F)$ iff there exists $B_n\in\bpi0{\beta_n}$ for some
$\beta_n<\alpha$ such that
$$A=\Intersect_{n\in\omega}\comp{B_n}.$$
Define $\bsig0\alpha(F)=\{\comp{B}: B\in \bpi0\alpha(F)\}$,
$\bdel0\alpha(F)=\bpi0\alpha(F)\intersect \bsig0\alpha(F)$, and
$\borel(F)=\Union_{\alpha<\omega_1}\bsig0\alpha(F)$.
Also, as before, let $\ord(F)$ be the least $\alpha$ such
that $\borel(F)=\bsig0\alpha(F)$.

\sdex{$?bsig0\alpha(F)$}
\sdex{$?bpi0\alpha(F)$}
\sdex{$?bdel0\alpha(F)$}
\sdex{$\borel(F)$}

\begin{theorem} \label{bbm}
(Bing, Bledsoe, Mauldin \site{bbm}) Suppose $F\subseteq P(\cantorspace)$
  is a countable family such that
  $\borel(\cantorspace)\subseteq \borel(F)$.  Then $\ord(F)=\omega_1$.
\end{theorem}

\begin{corollary} \label{bbmcor}
  Suppose $X$ is any space containing a perfect set and
  $F\subseteq P(X)$   is a countable family such that
  $\borel(X)\subseteq Borel(F)$. Then $\ord(F)=\omega_1$.
\end{corollary}
\proof
Suppose $\cantorspace\subseteq X$ and let
$\hat{F}=\{A\intersect\cantorspace: A\in F\}$.
By Theorem \ref{subsp} we have that
$\borel(\cantorspace)\subseteq \borel(\hat{F})$ and
so by Theorem \ref{bbm} we know $\ord(\hat{F})=\omega_1$.
But this implies $\ord(F)=\omega_1$.
\qed

The proof of Theorem \ref{bbm} is a generalization of Lebesgue's
universal set argument.  We need to prove the following two
lemmas.

\begin{lemma}\label{univ2}
 (Universal sets) Suppose $H\subseteq P(X)$ is countable
and define $R=\{A\cross B : A\subseteq\cantorspace $ is clopen and $B\in H\}$.
Then for every $\alpha$ with $1\leq\alpha<\omega_1$ there exists
$U\subseteq\cantorspace\cross X$ with $U\in \bpi0\alpha(R)$ such that
for every $A\in \bpi0\alpha(H)$ there exists $x\in\cantorspace$ with
$A=U_x$.
\end{lemma}

\proof
This is proved exactly as Theorem \ref{univ}, replacing the basis
for $X$ with $H$.
Note that when we replace $U_n$ by $U_n^*$ it is necessary to
prove by induction on $\beta$ that for every set $A\in \bpi0\beta(R)$ and
$n\in\omega$ that
the set $$A^*=\{(x,y):(x_n,y)\in A\}$$ is also in $\bpi0\beta(R)$.
\qed

\begin{lemma} \label{diag}
Suppose $H\subseteq P(\cantorspace)$, $R$ is
defined as in Lemma \ref{univ2}, and
$$\borel(\cantorspace)\subseteq \borel(H).$$ Then for
every set $A\in\borel(R)$ the set
$D=\{x: (x,x)\in A\}$ is in $\borel(H)$.
\end{lemma}
\proof
If $A=B\cross C$ where $B$ is clopen and $C\in H$, then
$D=B\intersect C$ which is in $\borel(H)$ by assumption.
Note that
$$\{x: (x,x)\in \Intersect_n A_n\}=\Intersect_n\{x:(x,x)\in A_n\}$$
and
$$\{x: (x,x)\in \comp A\}=\comp\{x: (x,x)\in A\},$$  so the result follows
by induction.
\qed

\par\noindent Proof of Theorem \ref{bbm}:

Suppose $\borel(H)=\bpi0\alpha(H)$ and let
$U\subseteq\cantorspace\cross\cantorspace$ be universal for
$\bpi0\alpha(H)$ given by Lemma \ref{univ2}.   By Lemma \ref{diag}
the set $D=\{x:(x,x) \in U\}$ is in $\borel(H)$ and hence
its complement is in $\borel(H)=\bpi0\alpha(H)$. Hence we get
the same old contradiction: if $U_x=\comp D $, then
$x\in D$ iff $x\notin D$.
\qed

\begin{theorem}\label{rec}
  (Rec{\l}aw) If $X$ is a second countable space
  and $X$ can be mapped continuously onto the unit interval,
  $[0,1]$, then $\ord(X)=\omega_1$.
\end{theorem}
\proof
Let $f:X\arrow [0,1]$ be continuous and onto. Let
$\base$ be a countable base for $X$ and let
$H=\{f(B):B\in\base\}$.  Since the preimage of an open subset
of $[0,1]$ is open in $X$ it is clear that
$\borel([0,1])\subseteq\borel(H)$. So by Corollary~\ref{bbmcor}
it follows that $\ord(H)=\omega_1$.  But $f$ maps the Borel
hierarchy of $X$ directly over to the hierarchy generated by $H$,
so $\ord(X)=\omega_1$.
\qed

Note that if $X$ is a discrete space of cardinality the continuum
then there is a continuous map of $X$ onto $[0,1]$ but
$\ord(X)=1$.

The Cantor space $\cantorspace$ can be mapped continuously onto
$[0,1]$ via the map
$$x\mapsto \sum_{n=0}^{\infty} {x(n)\over 2^{n+1}}.$$
This map is even one-to-one except at countably many points
where it is two-to-one.  It is also easy to see that $\reals$ can
be mapped continuously onto $[0,1]$ and
$\bairespace$ can be mapped onto $\cantorspace$.  It follows that in
Theorem \ref{rec} we may replace $[0,1]$ by $\cantorspace$,
$\bairespace$, or $\reals$.

Myrna Dzamonja points out that any completely regular space $Y$
which contains a perfect set can be mapped onto $[0,1]$.  This is
true because if $P\subseteq Y$ is perfect, then there is
a continuous map $f$ from $P$ onto $[0,1]$.  But since $Y$ is completely
regular this map extends to $Y$.

Rec{\l}aw did not publish his result, but I did, see Miller \site{proj}
and \site{survey2}.

\sector{Characteristic function of a sequence}

The idea of a \dex{characteristic function} of a sequence of sets
is due to Kuratowski and generalized the notion of a characteristic
function of a set introduced by de le Vall\'{e}e-Poussin.
The general notion was introduced by
Spilrajn \site{marcz} who also exploited it
in \site{marcz2}.
I think that Spilrajn also wrote under the name Marczewski
during World War II to disguise the fact that he was Jewish
from the Nazis.

Suppose $F\subseteq P(X)$ is a countable field of sets
(i.e. $F$ is a family of sets which is closed under complements in $X$
and finite intersections).
Let $F=\{A_n:n\in\omega\}$.
Define $c:X\arrow\cantorspace$ by
$$c(x)(n)=\left\{\begin{array}{ll}
1 & \mbox{if $x\in A_n$}\\
0 & \mbox{if $x\notin A_n$}
\end{array}\right.$$
Let $Y=c(X)$, then there is a direct correspondence
between $F$ and
$$\{C\intersect Y:\; C\subseteq\cantorspace \mbox{ clopen } \}.$$

\begin{theorem} (Spilrajn \site{marcz2}) \label{charseq}
  If $F\subseteq P(X)$ is a
  countable field of sets,
  then there exists a subspace $Y\subseteq\cantorspace$ such
  that $\ord(F)=\ord(Y)$.
\end{theorem}
\proof
If we define $x\eq y$ iff $\forall n\; (x\in A_n$ iff $y\in A_n)$, then
we see that members of $\borel(F)$ respect $\eq$. The preimages
of points of $Y$ under $c$ are exactly the equivalence classes of
$\eq$.  The map  $c$ induces a bijection between $X/\eq$ and $Y$ which
takes the family $F$ exactly to a clopen basis for the topology on $Y$.
Hence $\ord(F)=\ord(Y)$.
\qed

The following theorem says that bounded Borel hierarchies must have a top.

\begin{theorem} (Miller \site{borelhier})
Suppose $F\subseteq P(X)$ is a field of sets and
$$ord(F)=\lambda$$
where $\lambda$ is a countable limit ordinal. Then there exists
$B\in\borel(F)$ which is not in $\bpi0\alpha(F)$ for any
$\alpha<\lambda$.
\end{theorem}

\proof
By the characteristic function of a sequence of sets argument
we may assume without loss of generality that
$$F=\{C\intersect Y: \mbox{ $C\subseteq 2^\kappa$ clopen }\}.$$
A set $C\subseteq 2^\kappa$ is clopen iff it is a finite union
of sets of the form $$[s]=\{x\in 2^\kappa: s\subseteq x\}$$ where
$s:D\arrow 2$ is a map with $D\in [\kappa]^{<\omega}$ (i.e. $D$ is a
finite subset of $\kappa$).  Note that by induction for every
$A\in\borel(F)$ there exists an $S\in [\kappa]^{\omega}$ (called
a support of $A$) with the property that
for every $x,y\in 2^\kappa$ if $x\res S =y\res S$ then
$(x\in A$ iff $y\in A)$.  That is to say, membership in $A$ is
determined by restrictions to $S$.
By a Lowenheim-Skolem kind of an argument, we can find a countable
$S\subseteq \kappa$ with the property that for every
$\alpha<\lambda$ and $s:D\arrow 2$ with $D\in [S]^{<\omega}$ if
$\ord(Y\intersect [s])>\alpha$ then there exists
$A\in\bsig0\alpha(F)\setminus \bdel0\alpha(F)$ such that
$A\subseteq [s]$ and $A$ is supported by $S$.
By permuting $\kappa$ around we may assume without loss of generality
that $S=\omega$.  Define
$$T=\{s\in\finseq: \ord(Y\intersect [s])=\lambda\}.$$
Note that $T$ is a tree, i.e., $s\subseteq t\in T$ implies $s\in T$.
Also for any $s\in T$ either $s\concat 0\in T$ or $s\concat 1 \in T$,
because
$$[s]=[s\concat 0]\union [s\concat 1].$$
Since $\langle\rangle\in T$ it must be that $T$ has an infinite branch.
Let $x:\omega\arrow 2$ be such that $x\res n\in T$ for all $n<\omega$.
For each $n$ define $$t_n=(x\res n)\concat (1-x(n))$$ and note that
$$2^\kappa=[x]\union \Union_{n\in\omega} [t_n]$$
is a partition of $2^\kappa$ into clopen sets and one closed set $[x]$.

\bigskip
\claim For every $\alpha<\lambda$ and $n\in\omega$
there exists $m>n$ with $$\ord(Y\intersect [t_m]) >\alpha.$$
\proof
Suppose not and let $\alpha$ and $n$ witness this.  Note that
$$[x\res n]=[x]\union\Union_{n\leq m<\omega} [t_m].$$
Since $\ord([x\res n]\intersect Y)=\lambda$ we know there exists
$A\in\bsig0{\alpha+3}(F)\setminus \bdel0{\alpha+3}(F)$ such
that $A\subseteq [x\res n]$ and $A$ is supported by $S=\omega$.
Since $A$ is supported by $\omega$ either $[x]\subseteq A$
or $A$ is disjoint from $[x]$.  But
if $\ord([t_m]\intersect Y)\leq\alpha$ for each $m>n$,
then
$$A_0=\Union_{n\leq m <\omega} (A\intersect [t_m])$$
is $\bsig0\alpha(F)$ and $A=A_0$ or $A=A_0\union [x]$ either of which
is $\bsig0\alpha(F)$ (as long as $\alpha>1$).
This proves the Claim. \qed

The claim allows us to construct a set which is not at a level below
$\lambda$ as follows. Let $\alpha_n<\lambda$ be a sequence unbounded
in $\lambda$ and let $k_n$ be a distinct sequence with
$\ord([t_{k_n}]\intersect Y)\geq \alpha_n$.  Let
$A_n\subseteq [t_{k_n}]$ be in $\borel(F)\setminus\bdel0{\alpha_n}(F)$.
Then $\union_n A_n$ is not at any level bounded below $\lambda$.
\qed

\begin{question}
 Suppose $R\subseteq P(X)$ is a ring of sets, i.e., closed
under finite unions and finite intersections.  Let $R_\infty$
be the $\sigma$-ring generated by $R$, i.e., the smallest family
containing $R$ and closed under countable unions and countable
intersections. For $n\in\omega$ define $R_n$ as follows.
$R_0=R$ and let $R_{n+1}$ be the family of countable
unions (if $n$ even) or family of countable intersections (if $n$ odd)
of sets from $R_n$.  If $R_\infty=\Union_{n<\omega}R_n$, then
must there be $n<\omega$ such that $R_\infty=R_n$?
\end{question}

\sector{Martin's Axiom}

The following result seems to be due to Rothberger \site{roth}
and Solovay \site{jensol}\site{marsol}.  The forcing we use is
due to Silver.
However, it is probably just another view of Solovay's `almost disjoint
sets forcing'. \sdex{Martin's Axiom}

\begin{theorem}
  \label{silverlem}
  Assuming Martin's Axiom if $X$ is any second countable
  Hausdorff space of cardinality less than the continuum,
  then $\ord(X)\leq 2$ and, in fact,
  every subset of $X$ is $G_\delta$.
\end{theorem}
\proof

Let $A\subseteq X$ be arbitrary and let $\base$ be a countable
base for the topology on $X$.
The partial order $\poset$ is defined as follows.
$p\in \poset$ iff $p$ is a finite consistent set of sentences of the
form
\begin{enumerate}
  \item  $``x \notin \name{U}_n\mbox{''}$ where $x\in X\setminus A$ or
  \item $``B\subseteq \name{U}_n\mbox{''}$ where $B\in\base$ and $n\in\omega$.
\end{enumerate}
Consistent means that there is not a pair of sentences
$``x \notin \name{U}_n\mbox{''},``B\subseteq \name{U}_n\mbox{''}$ in $p$
where $x\in B$.  The ordering on $\poset$ is reverse containment,
i.e. $p$ is stronger than $q$, $p\leq q$ iff $p\supseteq q$.

We call this forcing \dex{Silver forcing}.

\claim $\poset$ satisfies the ccc.
\proof
Note that since $\base$ is countable there are only countably many
sentences of the type $``B\subseteq \name{U}_n\mbox{''}$.  Also if
$p$ and $q$ have exactly the same sentences of this type then
$p\union q\in \poset$ and hence $p$ and $q$ are compatible.
It follows that $\poset$ is the countable union of filters and
hence we cannot find an uncountable set of pairwise incompatible
conditions.
\qed

For $x\in X\setminus A$ define
$$D_x=\{p\in\poset: \exists n \; ``x\notin \name{U}_n\mbox{''}\in p\}.$$
For $x\in A$ and $n\in\omega$ define
$$E_x^n=\{p\in\poset: \exists B\in \base \;
x\in B \rmand ``B\subseteq \name{U}_n\mbox{''}\in p\}.$$

\claim $D_x$ is dense for each $x \in X\setminus A$ and
$E_x^n$ is dense for each $x \in A$ and $n\in\omega$.
\proof
To see that $D_x$ is dense let $p\in\poset$ be
arbitrary.   Choose $n$ large enough so that
$\name{U_n}$ is not mentioned in $p$, then
$(p\union \{``x\notin \name{U}_n\mbox{''}\})\in \poset$.

To see that $E_x^n$ is dense let $p$ be arbitrary and let
$Y\subseteq X\setminus A$ be the set of elements of $X\setminus A$
mentioned by $p$.  Since $x\in A$ and $X$ is Hausdorff there
exists $B\in \base$ with $B\intersect Y=\emptyset$ and
$x\in B$. Then $q=(p\union\{``B\subseteq \name{U}_n\mbox{''}\})\in \poset$ and
$q\in E_x^n$.
\qed

Since the cardinality of $X$ is less than the continuum we can
find a $\poset$-filter $G$ with the property that
$G$ meets each $D_x$ for  $x \in X\setminus A$ and each
$E_x^n$ for $x \in A$ and $n\in\omega$.
Now define
$$U_n=\Union\{B: ``B\subseteq \name{U}_n\mbox{''}\in G\}.$$
Note that $A=\Intersect_{n\in\omega}U_n$ and so $A$ is
$G_\delta$ in $X$.
\qed

Spaces $X$ in which every subset is $G_\delta$ are called
\dex{Q-sets}.

The following question was raised during an email correspondence
with Zhou.

\begin{question}
  Suppose every set of reals of
  cardinality $\aleph_1$ is a $Q$-set.  Then is
  ${\goth p}>\omega_1$, i.e., is it true that for every
  family ${\cal F}\subseteq\infsets$ of size $\omega_1$
  with the finite intersection property there exists an
  $X\in\infsets$ with $X\subseteq^* Y$ for all $Y\in {\cal F}$?
\end{question}

It is a theorem of Bell \site{bell} that ${\goth p}$ is the
first cardinal for which MA for $\sigma$-centered forcing
fails.  Another result along this line due to Alan Taylor
is that ${\goth p}$ is the cardinality of the smallest non
$\gamma$-set, see Galvin and Miller \site{galmil}.
\sdex{${?goth p}$}

Fleissner and Miller \site{qsets} show it is consistent
to have a $Q$-set whose union with the rationals is not
a $Q$-set.

For more information on Martin's Axiom see Fremlin \site{frembook}.
For more on Q-sets, see Fleissner \site{fleiss1} \site{fleiss2},
Miller \site{survey1} \site{survey2}, Przymusinski \site{przym},
 Judah and Shelah \site{qsetzero} \site{qsetsier}, and
 Balogh \site{balogh}.

\sector{Generic $G_\delta$}

It is natural\footnote{  `Gentlemen, the great thing about this,
   like most of the
   demonstrations of the higher mathematics, is that it can
   be of no earthly use to anybody.'
                          -Baron Kelvin}
 to ask
\begin{center}
``What are the possibly lengths of
Borel hierarchies?''
\end{center}  In this section we present a way of
forcing a generic $G_\delta$.

Let $X$ be a Hausdorff space with a countable
base $\base$.  Consider the following forcing notion.
$p\in\poset$ iff it is a finite consistent set of sentences
of the form:
\begin{enumerate}
  \item ``$B\subseteq U_n$'' where $B\in\base$ and $n\in\omega$, or
  \item ``$x\notin U_n$'' where $x\in X$ and $n\in\omega$, or
  \item ``$x\in \Intersect_{n<\omega}U_n$'' where $x\in X$.
\end{enumerate}
Consistency means that we cannot say that both
``$B\subseteq U_n$''
and ``$x\notin U_n$'' if it happens that $x\in B$ and we cannot
say both ``$x\notin U_n$''
and ``$x\in \Intersect_{n<\omega}U_n$''. The ordering is reverse inclusion.
A $\poset$ filter $G$ determines a $G_\delta$ set $U$ as follows:
Let $$U_n=\Union\{B\in\base: ``B\subseteq U_n\mbox{''}\in G\}.$$
Let $U=\Intersect_n U_n$. If $G$ is $\poset$-generic over $V$, a
density argument shows that for every $x\in X$ we have that
$$x\in U \rmiff``x\in \Intersect_{n<\omega}U_n\mbox{''}\in G.$$

Note that $U$ is not in $V$ (as long as $X$ is infinite).
For suppose $p\in\poset$ and $A\subseteq X$ is in $V$ is such
that
$$p\forces \name{U}=\check{A}.$$
Since $X$ is infinite there exist $x\in X$ which is not mentioned
in $p$. Note that
$p_0=p\union\{``x\in \Intersect_{n<\omega}U_n\mbox{''}\}$ is consistent
and also
$p_1=p\union\{``x\notin U_n\mbox{''}\}$  is consistent
for all sufficiently large $n$ (i.e. certainly for $U_n$ not mentioned
in $p$.) But $p_0\forces x\in U$ and $p_1\forces x\notin U$,
and since $x$ is either in $A$ or not in $A$ we arrive at
a contradiction.

In fact, $U$ is not $F_\sigma$ in the extension (assuming $X$ is uncountable).
To see this we will first need to prove that $\poset$ has ccc.

\begin{lemma}
$\poset$ has ccc.
\end{lemma}
\proof
Note that $p$ and $q$ are compatible iff
$(p\union q)\in\poset$ iff $(p\union q)$ is a consistent set of
sentences.
Recall that there are three types of sentences:
\begin{enumerate}
  \item $B\subseteq U_n$
  \item $x\notin U_n$
  \item $x\in \Intersect_{n<\omega}U_n$
\end{enumerate}
where $B\in\base$, $n\in\omega$, and $x\in X$.
Now if for contradiction $A$ were an uncountable antichain, then
since there are only countably many sentences of type 1 above
we may assume that all $p\in A$ have the same set of type 1 sentences.
Consequently for each distinct pair $p,q\in A$ there must be
an $x\in X$ and $n$ such that either
$``x\notin U_n\mbox{''}\in p$ and
 $``x\in \Intersect_{n<\omega}U_n\mbox{''}\in q$ or vice-versa.
For each $p\in A$ let $D_p$ be the finitely many elements of $X$
mentioned by $p$ and let $s_p:D_p\arrow\omega$ be defined
by
$$s_p(x)=\left\{\begin{array}{ll}
0   & \mbox{if $``x \in \Intersect_{n<\omega}U_n\mbox{''}\in p$ }\\
n+1 & \mbox{if $``x\notin U_n\mbox{''}\in p$}
\end{array}\right.$$
But now $\{s_p:p\in A\}$ is an uncountable family of pairwise
incompatible finite partial functions from $X$ into $\omega$ which is
impossible. (FIN($X$,$\omega$) has the ccc.)
\qed

If $V[G]$ is a generic extension of a model $V$ which contains a
topological space $X$, then we let
$X$ also refer to the space in $V[G]$ whose
topology is generated by the open subsets of $X$ which are in $V$.

\begin{theorem} (Miller \site{borelhier}) \label{switch}
Suppose $X$ in $V$ is an uncountable
Hausdorff space with countable base $\base$ and
$G$ is $\poset$-generic over $V$. Then in $V[G]$
the $G_\delta$ set $U$ is not $F_\sigma$.
\end{theorem}
\proof
We call this argument the old switcheroo.
Suppose for contradiction
$$p\forces \Intersect_{n\in\omega} U_n=\Union_{n\in\omega} C_n
\mbox{ where $C_n$ are closed in $X$ }.$$
For $Y\subseteq X$ let $\poset(Y)$ be the elements of
$\poset$ which only mention $y\in Y$ in type 2 or 3 statements.
Let $Y\subseteq X$ be countable such that
\begin{enumerate}
\item $p\in \poset(Y)$ and
\item  for every $n$ and
$B\in\base$ there exists a maximal antichain $A\subseteq \poset(Y)$
which decides the statement ``$B\intersect C_n=\emptyset$''.
\end{enumerate}
Since $X$ is uncountable there exists $x\in X\setminus Y$.
Let
$$q=p\union \{``x\in \Intersect_{n\in\omega}U_n\mbox{''}\}.$$
Since $q$ extends $p$, clearly
$$q\forces x\in \Union_{n\in\omega} C_n$$
so there exists $r\leq q$ and $n\in\omega$ so that
$$r\forces x\in C_n.$$
Let $$r=r_0\union \{``x\in \Intersect_{n\in\omega}U_n\mbox{''}\}$$
where $r_0$ does not mention $x$.  Now we do the switch. Let
$$t=r_0\union \{``x\notin U_m\mbox{''}\}$$
where $m$ is chosen sufficiently large so that $t$ is a consistent
condition.   Since
$$t\forces x\notin \Intersect_{n\in\omega}U_n$$
we know that
$$t\forces x\notin C_n.$$
Consequently there exist $s\in \poset(Y)$ and $B\in \base$
such that
\begin{enumerate}
  \item $s$ and $t$ are compatible,
  \item $s\forces B\intersect C_n=\emptyset$, and
  \item $x\notin B$.
\end{enumerate}
But $s$ and $r$ are compatible, because $s$ does not mention $x$.
This is a contradiction since
$s\union r\forces x\in C_n$ and
$s\union r\forces x\notin C_n$.
\qed

\sector{$\alpha$-forcing}

In this section we generalize the forcing which produced a
generic $G_\delta$ to arbitrarily high levels of the Borel
hierarchy.  Before doing so we must prove some elementary
facts about well-founded trees.

Define $T\subseteq Q^{<\omega}$ to be a \dex{tree} iff
$s\subseteq t\in T$ implies $s\in T$.  Define the \dex{rank function}
$r:T\arrow \ordinals\union\{\infty\}$ of $T$ as follows:
\begin{enumerate}
  \item $r(s)\geq 0$ iff $s\in T$,
  \item $r(s)\geq \alpha+1$ iff $\exists q\in Q\;\; r(s\concat q)\geq \alpha$,
  \item $r(s)\geq \lambda$ (for $\lambda$ a limit ordinal)
  iff $r(s)\geq\alpha$ for every $\alpha<\lambda$.
\end{enumerate}
Now define $r(s)=\alpha$ iff $r(s)\geq\alpha$ but not
$r(s)\geq\alpha+1$ and $r(s)=\infty$ iff $r(s)\geq\alpha$
for every ordinal $\alpha$.

Define $[T]=\{x\in Q^\omega:\forall n \; \; x\res n\in T\}$. \sdex{$[T]$}
We say that $T$ is \dex{well-founded} iff $[T]=\emptyset$.

\begin{theorem} \label{tree}
  $T$ is well-founded iff $r(\langle\rangle)\in \ordinals$.
\end{theorem}
\proof
%Note that if $s\subseteq t\in T$ then $r(t)\leq r(s)$ (regarding
%$\infty \geq \alpha$ for every ordinal $\alpha$.

It follows
easily from the definition that if $r(s)$ is an ordinal,
 then
$$r(s)=\sup\{r(s\concat q)+1:q\in Q\}.$$
Hence, if $r(\langle\rangle)=\alpha\in \ordinals$ and $x\in [T]$, then
$$r(x\res (n+1))< r(x\res n)$$
 is a descending sequence of ordinals.

On the other hand, if $r(s)=\infty$ then for some $q\in Q$ we must
have $r(s\concat q)=\infty$.  So if $r(\langle\rangle)=\infty$
we can construct
(using the axiom of choice) a sequence $s_n\in T$ with
$r(s_n)=\infty$ and $s_{n+1}=s_{n}\concat x(n)$.  Hence $x\in\ [T]$.
\qed

Definition. $T$ is a \dex{nice $\alpha$-tree}   iff
\begin{enumerate}
  \item $T\subseteq \finseq$ is a tree,
  \item $r:T\arrow (\alpha+1)$ is its rank function
   ($r(\langle\rangle)=\alpha$),
  \item if $r(s)>0$, then for every $n\in\omega\;\; s\concat n\in T$,
  \item if $r(s)=\beta$ is a successor ordinal, then for every $n\in\omega
   \;\; r(s\concat n)=\beta-1$, and
  \item if $r(s)=\lambda$ is a limit ordinal, then $r(s\concat 0)\geq 2$
  and $r(s\concat n)$ increases to $\lambda$ as $n\arrow\infty$.
\end{enumerate}

It is easy to see that for every $\alpha<\omega_1$ nice $\alpha$-trees
exist.  For $X$ a Hausdorff space with countable base, $\base$, and $T$
a nice
$\alpha$-tree ($\alpha\geq 2$), define the partial order
$\poset=\poset(X,\base,T)$ which we call \dex{$\alpha$-forcing}
as follows:

\bigskip
\par\noindent $p\in\poset$ iff $p=(t,F)$ where
\begin{enumerate}
  \item $t:D\arrow\base$ where $D\subseteq \tzero=\{s\in T: r(s)=0\}$
  is finite,
  \item $F\subseteq \tone\cross X$ is finite where
  $$\tone=T\setminus \tzero=\{s\in T: r(s)>0\},$$
  \item if $(s,x),(s\concat n,y)\in F$, then $x\not=y$, and
  \item if $(s,x)\in F$ and $t(s\concat n)=B$, then $x\notin B$.
\end{enumerate}
The ordering on $\poset$ is given by $p\leq q$ iff $t_p\supseteq t_q$ and
$F_p\supseteq F_q$.
\sdex{$\tzero$}
\sdex{$\tone$}

\begin{lemma}
  $\poset$ has ccc.
\end{lemma}
\proof
Suppose $A$ is uncountable antichain.  Since there are only countably
many different $t_p$ without loss we may assume that there exists
$t$ such that $t_p=t$ for all $p\in A$.  Consequently for $p,q\in A$
the only thing that can keep $p\union q$ from being a condition is that
there must be an $x\in X$ and an $s,s\concat n\in \tone$ such that
$$(s,x),(s\concat n,x)\in (F_p\union F_q).$$
But now for each $p\in A$ let $H_p:X\arrow [\tone]^{<\omega}$ be
the finite partial function defined by
$$H_p(x)=\{s\in \tone: (s,x)\in F_p\}$$
where domain $H_p$ is $\{x:\exists s\in \tone \; (s,x)\in F_p\}$.
Then $\{H_p:p\in A\}$ is an uncountable antichain in the order
of finite partial functions from $X$ to $[\tone]^{<\omega}$, a countable set.
\qed

Define for $G$ a $\poset$-filter the set $U_s\subseteq X$  for $s\in T$
as follows:
\begin{enumerate}
  \item for $s\in \tzero$ let $U_s=B$ iff $\exists p\in G$ such that
  $t_p(s)=B$ and
  \item for $s\in \tone$ let
   $U_s=\Intersect_{n\in\omega} \comp U_{s\concat n}$
\end{enumerate}
Note that $U_s$ is a $\bpi0\beta(X)$-set where $r(s)=\beta$.

\begin{lemma} \label{easydense}
  If $G$ is $\poset$-generic over $V$ then in $V[G]$ we have
  that for every $x\in X$ and $s\in \tone$
  $$x\in U_s\iff\exists p\in G\; (s,x)\in F_p.$$
\end{lemma}
\proof
First suppose that $r(s)=1$ and note that the following set is dense:
$$D=\{p\in\poset: (s,x)\in F_p \rmor
\exists n\exists B\in\base \; x\in B \rmand t_p(s\concat n)=B\}.$$
To see this let $p\in\poset$ be arbitrary.
If $(s,x)\in F_p$ then $p\in D$ and we are already done.
If $(s,x)\notin F_p$ then let
$$Y=\{y:(s,y)\in F_p\}.$$
Choose $B\in \base$ with $x\in B$ and $Y$ disjoint from
$B$.  Choose $s\concat n$  not in the domain of $t_p$, and
let $q=(t_q,F_p)$ be defined by $t_q=t_p\union (s\concat n, B)$.
So $q\leq p$ and $q\in D$.  Hence $D$ is dense.

Now by definition $x\in U_s$ iff
$x\in\Intersect_{n\in\omega}\comp U_{s\concat n}$.  So let
$G$ be a generic filter and $p\in G\intersect D$.  If
$(s,x)\in F_p$ then we know that for every $q\in G$ and
for every $n$,  if $t_q(s\concat n)=B$ then $x\notin B$. Consequently,
$x\in U_s$.  On the other hand if $t_p(s\concat n)=B$ where
$x\in B$, then
$x\notin U_s$ and for every $q\in G$ it must be that $(s,x)\notin F_q$
(since otherwise $p$ and $q$ would be incompatible).

\bigskip
Now suppose $r(s)>1$.  In this case note that the following set is dense:
$$E=\{p\in \poset: (s,x)\in F_p \rmor \exists n\;
(s\concat n,x)\in F_p\}.$$
To see this let $p\in\poset$ be arbitrary.  Then either
$(s,x)\in F_p$ and already $p\in E$ or by choosing $n$ large enough
$q=(t_p,F_p\union\{(s\concat n,x)\})\in E$. (Note $r(s\concat n)>0$.)

Now assume the result is true for all $U_{s\concat n}$.
Let $p\in G\intersect E$. If $(s,x)\in F_p$ then for every
$q\in G$ and $n$ we have $(s\concat n, x)\notin F_q$ and so by
induction $x\notin U_{s\concat n}$ and so $x\in U_s$.
On the other hand if $(s\concat n, x)\in F_p$, then
by induction $x\in U_{s\concat n}$ and so $x\notin U_s$, and
so again for every $q\in G$ we have
 $(s,x)\notin F_q$.
\qed

The following lemma is the heart of the old \dex{switcheroo} argument
used in Theorem \ref{switch}.  Given any $Q\subset X$ define the
$\rank(p,Q)$ as follows:
$$\rank(p,Q)=\max\{r(s): (s,x)\in F_p \mbox{ for some } x\in X\setminus Q\}.$$

\begin{lemma}\label{rank} (Rank Lemma). For any $\beta\geq 1$ and
$p\in \poset$ there exists $\hat{p}$ compatible with $p$ such that
\begin{enumerate}
  \item $\rank(\hat{p},Q)<\beta+1$ and
  \item for any $q\in\poset$ if $\rank(q,Q)<\beta$, then
  $$\mbox{$\hat{p}$ and $q$ compatible implies $p$ and $q$ compatible.}$$
\end{enumerate}
\end{lemma}
\proof
Let $p_0\leq p$ be any extension which satisfies:
for any $(s,x)\in F_p$ and $n\in\omega$, if
 $r(s)=\lambda>\beta$ is a limit ordinal and $r(s\concat n)<\beta+1$, then
there exist $m\in\omega$ such that $(s\concat n\concat m,x)\in F_{p_0}$.
Note that since $r(s\concat n)$ is increasing to $\lambda$
there are only finitely many $(s,x)$ and $s\concat n$ to
worry about. Also $r(s\concat n\concat m)>0$ so
this is possible to do.

Now let $\hat{p}$ be defined as follows:
$$t_{\hat{p}}=t_p$$ and
$$F_{\hat{p}}=\{(s,x)\in F_{p_0}: x\in Q \rmor
r(s)<\beta+1\}.$$
Suppose for contradiction that there exists $q$ such that
$\rank(q,Q)<\beta$, $\hat{p}$ and $q$ compatible, but $p$
 and $q$ incompatible.
Since $p$ and $q$ are incompatible either
\begin{enumerate}
  \item there exists $(s,x)\in F_q$ and $t_p(s\concat n)=B$ with $x\in B$, or
  \item there exists $(s,x)\in F_p$ and $t_q(s\concat n)=B$ with $x\in B$, or
  \item there exists $(s,x)\in F_p$ and $(s\concat n, x)\in F_q$, or
  \item there exists $(s,x)\in F_q$ and $(s\concat n, x)\in F_p$.
\end{enumerate}
(1) cannot happen since $t_{\hat{p}}=t_p$ and so $\hat{p},q$ would be
incompatible. (2) cannot happen since $r(s)=1$ and $\beta\geq 1$ means that
$(s,x)\in F_{\hat{p}}$ and so again $\hat{p}$ and $q$ are incompatible.
 If (3) or (4) happens for $x\in Q$ then again (in case 3)
  $(s,x)\in F_{\hat{p}}$
or (in case 4) $(s\concat n, x)\in F_{\hat{p}}$ and so $\hat{p},q$
incompatible.

So assume $x\notin Q$.  In case (3) by the definition of
$\rank(q,Q)<\beta$ we know that $r(s\concat n)<\beta$.
Now since $T$ is a nice tree we know that either
$r(s)\leq\beta$ and so $(s,x)\in F_{\hat{p}}$ or $r(s)=\lambda$ a
limit ordinal.  Now if $\lambda\leq \beta$ then $(s,x)\in F_{\hat{p}}$.
If $\lambda>\beta$ then by our construction of $p_0$ there exist $m$ with
$(s\concat n\concat m,x)\in F_{\hat{p}}$ and so $\hat{p},q$ are incompatible.
Finally in case (4) since $x\notin Q$ and so $r(s)<\beta$ we have
that $r(s\concat n)<\beta$ and so
$(s\concat n, x)\in F_{\hat{p}}$ and
so $\hat{p},q$ are incompatible.
\qed

Intuitively, it should be that statements of small rank are forced by
conditions of small rank.  The next lemma will make this more precise.
Let $L_{\infty}(P_\alpha: \alpha<\kappa)$ be the infinitary propositional
logic with $\{P_\alpha: \alpha<\kappa\}$ as the atomic sentences.
\sdex{$L_{\infty}(P_\alpha: \alpha<\kappa)$}
Let  $\Pi_0$-sentences be the atomic ones, $\{P_\alpha: \alpha<\kappa\}$.
For any $\beta>0$ let $\theta$ be a
\dex{$\Pi_\beta$-sentence} iff there exists
$\Gamma\subseteq \Union_{\delta<\beta}\Pi_\delta$-sentences and
$$\theta=\Conj_{\psi\in \Gamma}\neg \psi.$$
Models for this propositional language can naturally be regarded
as subsets $Y\subseteq\kappa$ where we define
\begin{enumerate}
\item $Y\models P_\alpha$ iff $\alpha\in Y$,
\item $Y\models \neg\theta$ iff not $Y\models\theta$, and
\item $Y\models \Conj \Gamma$ iff
$Y\models\theta$ for every $\theta\in \Gamma$.
\end{enumerate}

\begin{lemma} (Rank and Forcing Lemma)\label{forcing}
Suppose $\rank:\poset\arrow\ordinals$ is any function on
a poset $\poset$ which satisfies the Rank Lemma \ref{rank}.
Suppose $\forces_{\poset}\name{Y}\subset\kappa$ and for
every $p\in\poset$ and $\alpha<\kappa$ if
$$p\forces \alpha\in \name{Y}$$
then there exist $\hat{p}$ compatible with $p$ such that
$\rank(\hat{p})=0$ and
$$\hat{p}\forces\alpha\in \name{Y}.$$

Then for every $\Pi_\beta$-sentence $\theta$ (in the ground model)
and every $p\in\poset$, if
$$p\forces``\name{Y}\models\theta\mbox{''}$$
then there exists $\hat{p}$ compatible with $p$ such
that $\rank(\hat{p})\leq\beta$ and
$$\hat{p}\forces``\name{Y}\models\theta\mbox{''}.$$
\end{lemma}
\proof
This is one of those lemmas whose statement is longer than its proof.
The proof is induction on $\beta$ and for $\beta=0$ the conclusion
is true by assumption.  So suppose $\beta>0$ and
$\theta=\Conj_{\psi\in \Gamma}\neg \psi$ where
$\Gamma\subseteq \Union_{\delta<\beta}\Pi_\delta$-sentences.
By the rank lemma there exists $\hat{p}$ compatible with
$p$ such that $\rank(\hat{p})\leq\beta$ and for every
$q\in\poset$ with $\rank(q)<\beta$ if
$\hat{p},q$ compatible  then  ${p},q$ compatible.
We claim that
$$\hat{p}\forces``\name{Y}\models\theta\mbox{''}.$$
Suppose not.   Then  there exists $r\leq \hat{p}$ and
$\psi\in\Gamma$ such that
$$r\forces``\name{Y}\models\psi\mbox{''}.$$
By inductive assumption there exists $\hat{r}$ compatible with
$r$ such that
$$\rank(\hat{r})<\beta$$
 such that
$$\hat{r}\forces``\name{Y}\models\psi\mbox{''}.$$
But $\hat{r},\hat{p}$ compatible implies
$\hat{r},p$ compatible, which is a contradiction because
$\theta\implies\neg\psi$ and so
$$p\forces``\name{Y}\models\neg\psi\mbox{''}.$$
\qed

Some earlier uses of rank in forcing arguments is Steel's
forcing \site{steel},\site{frieds},\site{harrsharp};
and Silver's analysis of the collapsing
algebra \site{silver}.

\sector{Boolean algebras}

In this section we consider the length of Borel hierarchies
generated by a subset of a {complete boolean algebra}.  We
find that the generators of the complete boolean algebra
associated with $\alpha$-forcing generate it in exactly
$\alpha+1$ steps.  We start by presenting some background
information.

Let $\bool$ be a \dex{cBa}, i.e, complete boolean algebra.  This means
that in addition to being a boolean algebra, infinite sums and
products, also exist; i.e., for any $C\subseteq \bool$ there
exists $b$ (denoted $\sum C$) such that
\begin{enumerate}
  \item $c\leq b$ for every $c\in C$ and
  \item for every $d\in \bool$ if
  $c\leq d$ for every $c\in C$, then $b\leq d$.
\end{enumerate}
Similarly we  define $\prod C=-\sum _{c\in C}-c$ where
$-c$ denotes the complement of $c$ in $\bool$.

A partial order $\poset$ is \dex{separative} iff for any
$p,q\in\poset$ we have
$$p\leq q \rmiff\forall r\in\poset (r\leq p
\mbox{ implies $q,r$ compatible}).$$

\begin{theorem} (Scott, Solovay see \site{jech}) \label{ss}
A partial order $\poset$ is separative iff
there exists a cBa $\bool$ such that $\poset\subseteq\bool$
is dense in $\bool$, i.e. for every $b\in\bool$
if $b>0$ then there exists $p\in\poset$ with $p\leq b$.
\end{theorem}

It is easy to check that the $\alpha$-forcing $\poset$ is
separative (as long as $\base$ is infinite):
If $p\not\leq q$ then either
\begin{enumerate}
  \item $t_p$ does not extend $t_q$, so there exists $s$ such
  that $t_q(s)=B$ and   either $s$ not in the
  domain of $t_p$ or $t_p(s)=C$ where $C\not=B$ and so in either case
  we can find $r\leq p$ with $r,q$ incompatible, or
  \item $F_p$ does not contain  $F_q$, so there exists
  $(s,x)\in (F_q\setminus F_p)$ and we can either add
  $(s\concat n,x)$ for sufficiently large $n$ or add
  $t_r(s\concat n)=B$ for some sufficiently large $n$ and
  some $B\in \base$ with $x\in B$ and get $r\leq p$ which
  is incompatible with $q$.
\end{enumerate}

The elegant (but as far as I am concerned mysterious) approach
to forcing using complete boolean algebras contains the following
facts:
\begin{enumerate}
  \item for any sentence $\theta$ in the forcing language
     $$\val{\theta}=\sum\{b\in\bool: b\forces\theta\}=
     \sum\{p\in\poset: p\forces\theta\}$$
     where $\poset$ is any dense subset of $\bool$,
   \item $p\forces\theta$ iff $p\leq\val{\theta}$,
   \item $\val{\neg\theta}=-\val{\theta}$,
   \item $\val{\theta\conj\psi}=\val{\theta}\conj\val{\psi}$,
   \item $\val{\theta\disj\psi}=\val{\theta}\disj\val{\psi}$,
   \item for any set $X$ in the ground model,
    $$ \val{\forall x\in \check{X}\;\theta(x)}=
    \prod_{x\in X}\val{\theta(\check{x}}.$$
\end{enumerate}

\par\noindent Definitions. For $\bool$ a cBa and
$C\subseteq \bool$ define:

$\bpi00(C)=C$ and

$\bpi0\alpha(C)=\{\prod\;\Gamma: \Gamma\subseteq
\{-c:c\in\Union_{\beta<\alpha}\bpi0\beta(C)\}\}$ for $\alpha>0$.

\medskip

$\ord(\bool)=\min\{\alpha:\exists C\subseteq \bool$ countable with
$\bpi0\alpha(C)=\bool\}$. \sdex{$?ord(?bool)$}

\begin{theorem}(Miller \site{borelhier}) \label{boolord}
  For every $\alpha\leq\omega_1$ there exists a countably generated
  ccc cBa $\bool$ with $\ord(\bool)=\alpha$.
\end{theorem}

\proof

Let $\poset$ be $\alpha$-forcing and $\bool$ be the cBa given by
the Scott-Solovay Theorem~\ref{ss}.  We will show that
$\ord(\bool)=\alpha+1$.

Let
$$C=\{p\in\poset: F_p=\emptyset\}.$$
$C$ is countable and we claim that $\poset\subseteq\bpi0\alpha(C)$.
Since $\bool=\bsig01(\poset)$ this will imply that
$\bool=\bsig0{\alpha+1}(C)$ and so  $\ord(\bool)\leq\alpha+1$.

First note that for any $s\in T$ with $r(s)=0$ and $x\in X$,
$$\val{x\in U_s}=\sum\{p\in C: \exists B\in\base\; t_p(s)=B\rmand
x\in B\}.$$

By Lemma \ref{easydense}
we know for generic filters $G$
  that for every $x\in X$ and $s\in \tone$
  $$x\in U_s\iff\exists p\in G\; (s,x)\in F_p.$$
Hence  $\val{x\in U_s}=\langle \emptyset,\{(s,x)\}\rangle$ since
if they are not equal, then
$$b= \val{x\in U_s}\;\Delta\;\langle \emptyset,\{(s,x)\}\rangle>0,$$
but letting $G$ be a generic ultrafilter with $b$ in it would lead
to a contradiction.
We get that for $r(s)>0$:
$$\langle \emptyset,\{(s,x)\}\rangle=\val{x\in U_s}=
\val{x\in \Intersect_{n\in\omega}\comp U_{s\concat n}}=
\prod_{n\in\omega} -\val{x\in U_{s\concat n}}.$$
Remembering that for $r(s\concat n)=0$ we have
$\val{x\in U_{s\concat n}}\in \bsig01(C)$, we see
by induction that for every $s\in \tone$ if $r(s)=\beta$ then
$$\langle \emptyset,\{(s,x)\}\rangle\in \bpi0\beta(C).$$
For any $p\in \poset$
$$p=\langle t_p,\emptyset\rangle\conj\prod_{(s,x)\in F_p}
\langle \emptyset,\{(s,x)\}\rangle .$$
So we have that $p\in\bpi0{\alpha}(C)$.

\bigskip

Now we will see that $\ord(\bool) > \alpha$. We use the following
Lemmas.

$\bool^+$ are the nonzero elements of $\bool$. \sdex{$?bool^+$}

\begin{lemma}\label{rankbool1}
If $r:\poset\arrow\ordinals$ is a rank function, i.e. it satisfies
the Rank Lemma~\ref{rank} and in addition $p\leq q$ implies
$r(p)\leq r(q)$, then if $\poset$ is dense in the cBa $\bool$
then $r$ extends to $r^*$ on $\bool^+$:
$$r^*(b)=\min\{\beta\in\ordinals: \exists C\subseteq \poset: b=\sum C
\rmand \forall p\in C\; r(p)\leq\beta\}$$
and still satisfies the Rank Lemma.
\end{lemma}
\proof
Easy induction.
\qed

\begin{lemma} \label{rankbool2}
  If $r:\bool^+\arrow\ord$ is a rank function and $E\subseteq \bool$
  is a countable collection of rank zero elements, then for any
  $a\in\bpi0\gamma(E)$ and $a\not=0$ there exists $b\leq a$
  with $r(b)\leq\gamma$.
\end{lemma}
\proof
To see this let $E=\{e_n:n\in\omega\}$
and let $\name{Y}$ be a name for the set in the generic extension
 $$Y=\{n\in\omega: e_n\in G\}.$$
Note that $e_n=\val{n\in \name{Y}}$.
For elements $b$ of $\bool$ in the complete subalgebra generated by
$E$ let us associate sentences $\theta_b$
of the infinitary propositional logic $L_\infty(P_n:n\in\omega)$
as follows:
$$\theta_{e_n}=P_n$$
$$\theta_{-b}=\neg\theta_b$$
$$\theta_{\prod R}=\Conj_{r\in R}\theta_r$$
Note that $\val {Y\models\theta_b}=b$ and if $b\in \bpi0\gamma(E)$ then
$\theta_b$ is a $\Pi_\gamma$-sentence.
The Rank and Forcing Lemma \ref{forcing} gives us
(by translating
$p\forces Y\models\theta_b$ into $p\leq \val{Y\models\theta_b}=b$) that:
\begin{quote}
For any $\gamma\geq 1$ and $p\leq b\in\bpi0\gamma(E)$ there exists a $\hat{p}$
compatible with $p$ such that $\hat{p}\leq b$ and
$r(\hat{p})\leq\gamma$.
\end{quote}
\qed

Now we use the lemmas to see that $\ord(\bool) > \alpha$.

Given any countable
$E\subseteq\bool$,
let $Q\subseteq X$ be countable so that
for any $e\in E$ there exists $H\subseteq\poset$ countable
so that $e=\sum H$ and for every $p\in H$
we have $\rank(p,Q)=0$.  Let $x\in X\setminus Q$ be arbitrary;
then we claim:
$$\val{x\in U_{\emptyseq}}\notin\bsig0{\alpha}(E).$$

We have chosen $Q$ so that $r(p)=\rank(p,Q)=0$ for any  $p\in E$
so the hypothesis of
Lemma \ref{rankbool2} is satisfied.
Suppose for contradiction
that $\val{x\in U_{\emptyseq}}=b\in\bsig0\alpha(E)$.
Let
$b=\sum_{n\in\omega}b_n$ where each $b_n$ is $\bpi0{\gamma_n}(C)$
for some $\gamma_n<\alpha$. For some n and $p\in\poset$ we would
have $p\leq b_n$.  By Lemma \ref{rankbool2} we have that there
exists
$\hat{p}$ with $\hat{p}\leq b_n\leq b=\val{x\in U_{\emptyseq}}$ and
$\rank(\hat{p},Q)\leq\gamma_n$. But by the definition of
$\rank(\hat{p},Q)$ the pair $({\emptyseq},x)$ is not in $F_{\hat{p}}$,
but this contradicts $\hat{p}\leq b_n\leq b=\val{x\in U_{\emptyseq}}=
\langle \emptyset, \{({\emptyseq},x)\}\rangle$.

\bigskip

This takes care of all countable successor ordinals.
(We leave the case of $\alpha=0,1$ for the reader to contemplate.)
For $\lambda$ a
limit ordinal take $\alpha_n$ increasing to $\lambda$ and let
$\poset=\sum_{n<\omega}\poset_{\alpha_n}$ be the direct sum,
 where $\poset_{\alpha_n}$ is $\alpha_n$-forcing.  Another way to
describe essentially the same thing is as follows:
Let $\poset_{\lambda}$ be $\lambda$-forcing. Then take $\poset$ to be
the subposet
of $\poset_{\lambda}$ such that $\emptyseq$ doesn't occur, i.e.,
$$\poset=\{p\in\poset_{\lambda}: \neg\exists x\in X\;\;
(\emptyseq,x)\in F_p\}.$$
Now if $\poset$ is dense in the cBa $\bool$, then $\ord(\bool)=\lambda$.
This is easy to see, because for each $p\in\poset$ there exists
$\beta<\lambda$ with $p\in\bpi0\beta(C)$.  Consequently,
$\poset\subseteq\Union_{\beta<\lambda}\bpi0\beta(C)$ and
so since $\bool=\bsig01(\poset)$ we get $\bool=\bsig0\lambda(C)$.
Similarly to the other argument we see that for any countable
$E$ we can choose a countable $Q\subseteq X$ such for any
$s\in T$ with $2\leq r(s)=\beta<\lambda$ (so $s\not=\emptyseq$) we have that
$\val{x\in U_s}$ is not $\bsig0\beta(E)$.  Hence
$\ord(\bool)=\lambda$.

For $\ord(\bool)=\omega_1$ we postpone until section \ref{ord1}.
\qed

\sector{Borel order of a field of sets}

In this section we use the Sikorski-Loomis representation theorem
to transfer the abstract Borel hierarchy on a complete boolean
algebra into a field of sets.

A family $F\subseteq P(X)$ is a \dex{$\sigma$-field} iff it contains
the empty set and  is closed under countable unions and complements
in $X$. $I\subseteq F$ is a \dex{$\sigma$-ideal} in $F$ iff
\begin{enumerate}
\item $I$ contains the empty set,
\item $I$ is closed under countable unions,
\item $A\subseteq B\in I$ and $A\in F$ implies $A\in I$, and
\item $X\notin I$.
\end{enumerate}

 $F/I$ is the countably complete boolean algebra
formed by taking $F$ and modding out by $I$, i.e. $A\eq B$ iff
$A\Delta B\in I$. \sdex{$F/I$}
For $A\in F$ we use $[A]$ or $[A]_I$ to denote the equivalence
class of $A$ modulo $I$. \sdex{$[A]_I$}

\begin{theorem}(Sikorski,Loomis, see \site{sikor} section 29)
\label{sikorthm}
For any countably complete boolean algebra $B$ there exists
a $\sigma$-field $F$ and a $\sigma$-ideal $I$ such that
$B$ is isomorphic to $F/I$.
\end{theorem}
\proof
Recall that the Stone space of $B$, $\st(B)$, is the space of ultrafilters
$u$ on $B$ with the topology generated by the clopen sets
of the form: \sdex{$\st(B)$}
$$[b]=\{u\in\st(B): b\in u\}.$$
This space is a compact Hausdorff space in which the field of clopen sets
exactly corresponds to $B$.   $B$ is countably complete means that for
for any sequence $\{ b_n:n<\omega\}$ in $B$ there exists
$b\in B$ such that $b=\sum_{n\in\omega} b_n$. This translates to
the fact that given any countable family of clopen sets $\{C_n:n\in\omega\}$
in $\st(B)$ there exists a clopen set $C$ such that
$\Union_{n\in\omega}C_n\subseteq C$ and the closed set
$C\setminus \Union_{n\in\omega}C_n$ cannot contain a clopen set, hence
it has no interior, so it is nowhere dense. Let $F$ be the $\sigma$-field
generated by the clopen subsets of $\st(B)$.  Let $I$ be the $\sigma$-ideal
generated by the closed nowhere dense subsets of $F$
(i.e. the ideal of meager sets). The Baire
category theorem implies that no nonempty open subset of a compact
Hausdorff space is meager, so
$st(B)\notin I$ and the same holds for any nonempty clopen subset of
$\st(B)$.  Since the countable union of clopen sets is equivalent to
a clopen set modulo $I$ it follows that the map $C\mapsto [C]$ is
an isomorphism taking the clopen algebra of $\st(B)$ onto $F/I$.
\qed

Shortly after I gave a talk about my boolean algebra result
(Theorem \ref{boolord}),
Kunen pointed out the following result.

\begin{theorem}
  (Kunen see \site{borelhier}) For every $\alpha\leq\omega_1$
   there exists a field
  of sets $H$ such that $\ord(H)=\alpha$.
\end{theorem}
\proof
Clearly we only have to worry about $\alpha$ with $2<\alpha<\omega_1$.
Let $\bool$ be the complete boolean algebra given by Theorem~\ref{boolord}.
Let $\bool\isom F/I$ where $F$ is a $\sigma$-field of sets and
$I$ a $\sigma$-ideal.  Let $C\subseteq F/I$ be a countable
set of generators. Define
$$H=\{A\in F: [A]_I\in C\}.$$
By induction on $\beta$ it is easy to prove that for any $Q\in F$:
$$Q\in \bsig0\beta(H) \rmiff [Q]_I\in \bsig0\beta(C).$$
From which it follows that  $\ord(H)=\alpha$.
\qed

Note that there is no claim that the family $H$ is countable.
In fact, it is consistent (Miller \site{borelhier}) that
for every countable $H$ either $\ord(H)\leq 2$ or $\ord(H)=\omega_1$.

\sector{CH and orders of separable metric spaces}

In this section we prove that assuming CH that there
exists countable field of sets of all possible Borel orders,
which we know is equivalent to existence of separable
metric spaces of all possible orders.
We will need a sharper form of the representation theorem.

\begin{theorem}(Sikorski, see \site{sikor} section 31)\label{sikorccc}
$\bool$ is a countably generated ccc cBa iff there exists
a ccc $\sigma$-ideal $I$ in $\borel(\cantorspace)$ such that
$\bool\isom\borel(\cantorspace)/I$.  Furthermore if
$\bool$ is generated by the countable set $C\subseteq\bool$, then
this isomorphism can be taken so as to map the
clopen sets mod $I$ onto $C$.
\end{theorem}
\proof
To see that $\borel(\cantorspace)/I$ is countably generated is
trivial since the clopen sets modulo $I$ generate it.
A general theorem of Tarski is
that any $\kappa$-complete $\kappa$-cc boolean algebra
is complete.

For the other direction, we may assume by using the Sikorski-Loomis
Theorem, that $\bool$ is $F/J$ where $F$ is a $\sigma$-field
and $J$ a $\sigma$-ideal in $F$.  Since $\bool$ is countably
generated there exists $C_n\in F$ for $n\in\omega$ such that
$\{[C_n]:n\in\omega\}$ generates $F/J$ where
$[C]$ denotes the equivalence class of $C$ modulo $J$. Now
let $h:X\arrow\cantorspace$ be defined by
$$h(x)(n)=\left\{\begin{array}{ll}
                     1 & \mbox{if $x\in C_n$}\\
                     0 & \mbox{if $x\notin C_n$}
                 \end{array}\right.$$
and define $\phi:\borel(\cantorspace)\arrow F$ by
$$\phi(A)=h^{-1}(A).$$
Define $I=\{A\in\borel(\cantorspace):\phi(A)\in J\}$.
Finally, we claim that
$$\hat{\phi}:\borel(\cantorspace)/I\arrow F/I
\mbox{ defined by }\hat{\phi}([A]_I)=[\phi(A)]_J$$
is an isomorphism of the two boolean algebras.
\qed

For $I$ a $\sigma$-ideal in $\borel(\cantorspace)$ we say that
$X\subseteq\cantorspace$ is an \dex{$I$-Luzin set} iff for every
$A\in I$ we have that $X\intersect A$ is countable.  We say that
$X$ is \dex{super-$I$-Luzin} iff $X$ is $I$-Luzin and for every
$B\in\borel(\cantorspace)\setminus I$ we have that
$B\intersect X\not=\emptyset$.  The following Theorem was
first proved by Mahlo \site{mahlo} and later by Luzin \site{luz} for
the ideal of meager subsets of the real line.
Apparently, Mahlo's paper was overlooked and hence these
kinds of sets have always been referred to as Luzin sets.

\begin{theorem} (Mahlo \site{mahlo}) CH.  \label{mah}
  Suppose $I$ is a $\sigma$-ideal in $\borel(\cantorspace)$ containing
  all the singletons.  Then there exists a super-$I$-Luzin set.
\end{theorem}
\proof
Let $$I=\{A_\alpha:\alpha<\omega_1\}$$ and let
$$\borel(\cantorspace)\setminus I=\{B_\alpha:\alpha<\omega_1\}.$$
Inductively choose $x_\alpha\in\cantorspace$ so that
$$x_\alpha\in B_\alpha\setminus (\{x_\beta:\beta<\alpha\}\union
\Union_{\beta<\alpha}A_\alpha).$$
Then $X=\{x_\alpha:\alpha<\omega_1\}$ is  a super-$I$-Luzin set.
\qed

\begin{theorem} (Kunen see \site{borelhier})
   Suppose $\bool=\borel(\cantorspace)/I$ is a cBa, $C\subseteq\bool$
   are the clopen mod $I$ sets,
   $ord(C)=\alpha>2$, and $X$ is
   super-$I$-Luzin.   Then $\ord(X)=\alpha$.
\end{theorem}

\proof
Note that the $\ord(X)$ is the minimum $\alpha$ such that
for every $B\in\borel(\cantorspace)$ there exists
$A\in\bpi0\alpha(\cantorspace)$ with $A\intersect X=B\intersect X$.

Since $\ord(C)=\alpha$ we know that given any Borel set $B$
there exists a $\bpi0\alpha$ set
$A$ such that $A\Delta B\in I$.  Since $X$ is Luzin we know
that $X\intersect (A\Delta B)$ is countable.  Hence there exist countable
sets $F_0,F_1$ such that
$$X\intersect B= X\intersect ((A\setminus F_0)\union F_1).$$
But since $\alpha>2$ we have that $((A\setminus F_0)\union F_1)$ is
also $\bpi0\alpha$ and hence $\ord(X)\leq\alpha$.

On the other hand for any $\beta<\alpha$ we know there exists
a Borel set $B$ such that for every $\bpi0\beta$ set $A$ we have
$B\Delta A\notin I$ (since $\ord(C)>\beta$).
But since $X$ is super-$I$-Luzin we have that for every
$\bpi0\beta$ set $A$ that $X\intersect (B\Delta A)\not=\emptyset$ and
hence $X\intersect B\not= X\intersect A$.  Consequently,
$\ord(X)>\beta$.\qed

\begin{corollary}
  (CH) For every $\alpha\leq\omega_1$ there exists a separable
  metric space $X$ such that $\ord(X)=\omega_1$.
\end{corollary}

While a graduate student at Berkeley I had obtained the result
that it was consistent with
any cardinal arithmetic to assume that for every $\alpha\leq\omega_1$
there exists a separable metric space $X$ such that $\ord(X)=\alpha$.
It never occured to me at the time to ask what CH implied. In fact,
my way of thinking at the time was that
proving something from CH is practically the same as just showing
it is consistent.  I found out in the real world
(outside of Berkeley) that they are considered very differently.

In Miller \site{borelhier} it is shown that
for every $\alpha<\omega_1$ it is consistent there
exists a separable metric space of order $\beta$ iff
$\alpha<\beta\leq\omega_1$.  But the general question
is open.

\begin{question}
  For what $C\subseteq\omega_1$ is it consistent that
  $$C=\{\ord(X): X \mbox{ separable metric } \}?$$
\end{question}

\sector{Martin-Solovay Theorem}

In this section we the theorem below.  The technique of proof
will be used in the next section to produce a boolean algebra
of order $\omega_1$. \sdex{Martin-Solovay Theorem}

\begin{theorem} (Martin-Solovay \site{marsol}) \label{marsolthm}
  The following are equivalent for an infinite cardinal $\kappa$:
  \begin{enumerate}
        \item $\MA_\kappa$, i.e., for any poset $\poset$ which is ccc
        and family $\dense$ of dense subsets of $\poset$ with
        $|\dense|<\kappa$ there exists a $\poset$-filter $G$ with
        $G\intersect D\not=\emptyset$ for all $D\in\dense$
        \item For any ccc $\sigma$-ideal $I$ in $\borel(\cantorset)$
        and $\ideal\subset I$ with $|\ideal|<\kappa$ we have that
        $$\cantorspace\setminus \Union\ideal\not=\emptyset.$$
  \end{enumerate}
\end{theorem} \sdex{$\MA_\kappa$}

\begin{lemma}  \label{martinlem}
  Let $\bool=\borel(\cantorspace)/I$ for some ccc $\sigma$-ideal $I$ and
  let $\poset=\bool\setminus\{0\}$.   The following are equivalent
  for an infinite cardinal $\kappa$:
 \begin{enumerate}
  \item for any family $\dense$ of dense subsets of $\poset$ with
        $|\dense|<\kappa$ there exists a $\poset$-filter $G$ with
        $G\intersect D\not=\emptyset$ for all $D\in\dense$
  \item for any family $\famly\subseteq\bool^\omega$ with
        $|\famly|<\kappa$ there exists an ultrafilter $U$ on
        $\bool$ which is $\famly$-complete, i.e., for every
        $\langle b_n:n\in\omega\rangle\in \famly$
        $$\sum_{n\in\omega}b_n\in U\rmiff \exists n \; b_n\in U$$
  \item for any $\ideal\subset I$ with $|\ideal|<\kappa$
        $$\cantorspace\setminus \Union\ideal\not=\emptyset$$
\end{enumerate}
\end{lemma}
\proof
To see that (1) implies (2) note that for any
$\langle b_n:n\in\omega\rangle\in\bool^\omega$ the set
$$D=\{p\in\poset: p\leq -\sum_n b_n \rmor\exists n\; p\leq b_n\}$$
is dense. Note also that any filter extends to an ultrafilter.

To see that (2) implies (3) do as follows.
Let $M$ be an elementary substructure of
H$_{\gamma}$ for sufficiently large $\gamma$
with $|M|<\kappa$, $I\in M$, $\ideal\subseteq M$.

Let $\famly$ be all the $\omega$-sequences of Borel sets which
are in $M$. Since $|\famly|<\kappa$ we know there exists
$U$ an $\famly$-complete ultrafilter on $\bool$.
Define $x\in\cantorset$ by the
rule:
$$ x(n)=i \rmiff [\{y\in\cantorspace: y(n)=i\}]\in U.$$

\claim  For every Borel set $B\in M$:
$$x\in B\rmiff [B]\in U.$$
\proof
This is true for subbasic clopen sets by definition.
Inductive steps just use that $U$ is an M-complete ultrafilter.
\qed

\bigskip

To see that (3) implies (1), let $M$ be an elementary substructure of
H$_{\gamma}$ for sufficiently large $\gamma$
with $|M|<\kappa$, $I\in M$, $\dense\subseteq M$. Let
$$\ideal=M\intersect I.$$
By (3) there exists
$$x\in \cantorset\setminus\Union\ideal.$$
Let $\bool_M=\bool\intersect M$. Then define
$$G=\{[B]\in\bool_M: x\in B\}.$$
Check $G$ is a $\poset$ filter which meets every $D\in\dense$.
\qed

This proves Lemma \ref{martinlem}.

\bigskip

To prove the theorem it necessary to do a \dex{two step iteration}.
Let $\poset$ be a poset and $\name{\posetq}\in V^\poset$ be the
$\poset$-name of a poset, i.e.,
$$\forces_\poset\posetq \mbox{ is a poset.}$$
Then we form the poset
$$\poset * \name{\posetq}= \{(p,\name{q}):
p\forces \name{q}\in\name{\posetq}\}$$
ordered by $(\hat{p},\hat{q})\leq (p,q)$ iff $\hat{p}\leq p$ and
$\hat{p}\forces \hat{q}\leq q$.\sdex{$?poset * ?name{?posetq}$}
In general there are two problems
with this.  First, $\poset * \name{\posetq}$ is
a class.  Second, it does not satisfy antisymmetry:
$x\leq y$ and $y\leq x$ implies $x=y$.   These can be solved
by cutting down to a sufficiently large set of nice names and
moding out by the appropriate equivalence relation.  Three of the main
theorems are:

\begin{theorem}
  If $G$ is $\poset$-generic over $V$ and $H$ is
  $\posetq^G$-generic over $V[G]$, then
  $$G*H=\{(p,q)\in \poset * \name{\posetq}: p\in G, q^G\in H\}.$$
  is a $\poset * \name{\posetq}$ filter generic over $V$.
\end{theorem}

\begin{theorem}
  If $K$ is a $\poset * \name{\posetq}$-filter generic over $V$, then
  $$G=\{p:\exists q\;(p,q)\in K\}$$
   is $\poset$-generic over $V$ and
  $$H=\{q^G:\exists p \; (p,q)\in K\}$$ is $\posetq^G$-generic over $V[G]$.
\end{theorem}

\begin{theorem} (Solovay-Tennenbaum  \site{solten})
  If $\poset$ is ccc and $\forces_\poset $``$\name{\posetq}$ is ccc'', then
   $\poset * \name{\posetq}$ is ccc.
\end{theorem}

Finally we prove Theorem \ref{marsolthm}.  We may assume that
the ccc poset $\poset$ has
cardinality less than $\kappa$ and $\kappa\leq c$.
Choose $X=\{x_p:p\in\poset\}\subseteq\cantorspace$ distinct elements
of $\cantorset$.
If $G$ is $\poset$-filter generic over $V$ let
$\posetq$ be Silver's forcing for forcing a $G_\delta$-set,
$\Intersect_{n\in\omega}U_n$, in $X$ such that
$$G=\{p\in\poset: x_p\in \Intersect_{n\in\omega}U_n\}.$$
Let $\base\in V$ be a countable base for $X$.
A simple description of $\poset * \name{\posetq}$ can
be given by:
$$(p,q)\in \poset * \name{\posetq}$$
iff
$p\in\poset$ and
$q\in V$ is a finite set of consistent sentences of the form:
\begin{enumerate}
  \item  $``x \notin \name{U}_n\mbox{''}$ where $x\in X$ or
  \item $``B\subseteq \name{U}_n\mbox{''}$ where $B\in\base$ and $n\in\omega$.
\end{enumerate}
with the additional requirement that whenever the sentence
$``x \notin \name{U}_n\mbox{''}$ is in $q$ and $x=x_r$, then
$p$ and $r$ are incompatible (so $p\forces r\notin G$).

Note that if $D\subseteq \poset$ is dense in $\poset$, then
$D$ is predense in $\poset * \name{\posetq}$, i.e., every
$r\in \poset * \name{\posetq}$ is compatible with an element of $D$.
Consequently, it is enough to find sufficiently generic filters
for $\poset * \name{\posetq}$.   By Lemma \ref{martinlem}
and Sikorski's Theorem~\ref{sikorccc} it is enough to see that
if $\poset * \name{\posetq}\subseteq \bool$ is dense in the ccc cBa algebra
$\bool$, then $\bool$ is countably generated. Let
$$C=\{\val{B\subseteq U_n}:B\in\base, n\in\omega\}.$$
We claim that $C$ generates $\bool$.  To see this, note that
for each $p\in\poset$
$$\val{x_p\in \intersect_n U_n}=\prod_{n\in\omega}\val{x_p\in U_n}$$
$$\val{x_p\in U_n}=\sum_{B\in\base, x_p\in B}\val{B\subseteq U_n}$$
furthermore
$$(p,\emptyset)=\val{x_p\in \intersect_n U_n}$$
and so it follows that every element of $\poset * \name{\posetq}$
is in the boolean algebra generated by $C$ and so since
$\poset * \name{\posetq}$ is dense in $\bool$ it follows that
$C$ generates $\bool$.
\qed

Define $X\subseteq\cantorspace$ to be
a generalized $I$-Luzin set for an ideal $I$ in the Borel sets iff
$|X|=\continuum$ and $|X\intersect A|<\continuum$
for every $A\in I$.
It follows from the Martin-Solovay Theorem \ref{marsolthm} that
(assuming that the continuum is regular)

MA is equivalent to

for every ccc ideal $I$ in the Borel subsets of $\cantorspace$ there
exists a generalized $I$-Luzin set.

Miller and Prikry \site{prikmil} show that it necessary to assume
the continuum is regular in the above observation.

\sector{Boolean algebra of order $\omega_1$}

Now we use the Martin-Solovay technique to produce a countably \label{ord1}
generated ccc cBa with order $\omega_1$.   Before doing so we
introduce a countable version of $\alpha$-forcing which will be
useful for other results also.  It is similar
to one used in Miller \site{cat} to give a simple
proof about generating sets in the category algebra.

Let
$T$ be a nice tree of rank $\alpha$  ($2\leq\alpha<\omega_1$).
Define
$$\poset_\alpha=\{p:D\arrow\omega: D\in [\omega]^{<\omega},
\forall s,s\concat n\in D\;\; p(s)\not=p(s\concat n)\}.$$
This is ordered by $p\leq q$ iff $p\supseteq q$. \sdex{$?poset_?alpha$}
For $p\in\poset_\alpha$ define
$$\rank(p)=\max\{r_T(s):s\in\dom(p)\}$$
where $r_T$ is the rank function on $T$. \sdex{$\rank(p)$}

\begin{lemma} \label{littlerank} $\rank:\poset_\alpha\arrow \alpha+1$
satisfies the
Rank Lemma \ref{rank}, i.e, for every $p\in\poset_\alpha$ and
$\beta\geq 1$ there exists $\hat{p}\in\poset_\alpha$ such that
\begin{enumerate}
\item $\hat{p}$ is compatible with $p$,
\item $\rank(\hat{p})\leq\beta$, and
\item for any $q\in\poset_\alpha$ if $\rank(q)<\beta$ and
$\hat{p}$ and $q$ are compatible, then $p$ and $q$ are compatible.
\end{enumerate}
\end{lemma}
\proof
First let $p_0\leq p$ be such that for every $s\in\dom(p)$ and $n\in\omega$
if $r_T(s\concat n)<\beta<\lambda=r_T(s)$ then there exists
$m\in\omega$ with $p_0(s\concat n\concat m)=p(s)$.
Note that $r_T(s\concat n)<\beta<\lambda=r_T(s)$ can happen only when
$\lambda$ is a limit ordinal and for any such $s$ there can be at most
finitely many $n$ (because $T$ is a nice tree).

Now let
$$E=\{s\in\dom(p_0): r_T(s)\leq\beta\}$$
and define
$\hat{p}= p_0\res E$. It is compatible with $p$ since $p_0$ is stronger
than both. From its definition it has rank $\leq\beta$.  So let
$q\in\poset_\alpha$ have $\rank(q)<\beta$ and be incompatible with
$p$.  We need to show it is incompatible with $\hat{p}$. There are
only three ways for $q$ and $p$ to be incompatible:
\begin{enumerate}
  \item $\exists s\in \dom(p)\intersect\dom(q)\;\; p(s)\not=q(s)$,
  \item $\exists s\in\dom(q)\;\exists s\concat n\in \dom(p)
         \;\; q(s)=p(s^\concat n)$, or
  \item $\exists s\in\dom(p)\;\exists s\concat n\in \dom(q)
         \;\; p(s)=q(s^\concat n)$.
\end{enumerate}
For (1) since $\rank(q)<\beta$ we know $r_T(s)<\beta$ and hence
by construction $s\in\dom(\hat{p})$ and so $q$ and $\hat{p}$ are
incompatible. For (2) since $r_T(s\concat n)<r_T(s)<\beta$
we get the same conclusion.  For (3) since $s\concat n\in \dom(q)$ we
know $r_T(s\concat n)<\beta$. If $r_T(s)=\beta$, then
$s\in \dom(\hat{p})$ and so $q$ and $\hat{p}$ are
incompatible.  Otherwise since $T$ is a nice tree,
$r_T(s\concat n)<\beta<r_T(s)=\lambda$ a limit ordinal. In this
case we have arranged $\hat{p}$ so that there exists $m$ with
$p(s)=\hat{p}(s\concat n\concat m)$ and so again
$q$ and $\hat{p}$ are
incompatible. \qed

\begin{lemma} \label{denseset}
  There exists a countable family $\dense$ of dense subsets of
  $\poset_\alpha$ such that for every $G$ a  $\poset_\alpha$-filter
  which meets each dense set in $\dense$ the filter $G$ determines
  a map $x:T\arrow\omega$ by $p\in G$ iff $p\subseteq x$. This
  map has the property that for every $s\in\tone$ the value
  of $x(s)$ is the unique element of $\omega$ not in
  $\{x(s\concat n):n\in\omega\}$.
\end{lemma}
\proof
For each $s\in T$ the set  $$D_s=\{p : s\in\dom(p)\}$$ is dense.
Also for each $s\in\tone$ and $k\in\omega$ the set
$$E^k_s=\{p : p(s)=k \rmor \exists n \;\;p(s\concat n)=k\}$$
is dense.
\qed

The poset $\poset_\alpha$ is separative, since if $p\not\leq q$ then
either $p$ and $q$ are incompatible or there exists
$s\in\dom(q)\setminus\dom(p)$ in which case we can find
$\hat{p}\leq p$ with $\hat{p}(s)\not=q(s)$.

Now if $\poset_\alpha\subseteq\bool$ is dense in the cBa $\bool$, it
follows that for each $p\in \poset_\alpha$
$$p=\val{p\subseteq x}$$
and for any $s\in \tone$ and $k$
$$ \val{x(s)=k}=\prod_{m\in\omega}\val{x(s^\concat m)\not= k}.$$
Consequently if
$$C=\{p\in\poset_\alpha: \dom(p)\subseteq T^0\}$$
then $C\subseteq\bool$ has the property that $\ord(C)=\alpha+1$.

Now let
$\sum_{\alpha<\omega_1}\poset_\alpha$ be the \dex{direct sum},
i.e., $p=\langle p_\alpha:\alpha<\omega_1\rangle$
with $p_\alpha\in\poset_\alpha$ and $p_\alpha={\bf 1}_\alpha=\emptyset$
for all but finitely many $\alpha$.
\sdex{$\sum_{\alpha<\omega_1}?poset_\alpha$}
  This forcing is equivalent
to adding $\omega_1$ Cohen reals, so the usual delta-lemma
argument shows that it is ccc.
Let
$$X=\{x_{\alpha,s,n}\in\cantorspace:\alpha<\omega_1,s\in\tzero_\alpha,
n\in\omega\}$$
be distinct elements of $\cantorspace$. For
$G=\langle G_\alpha:\alpha<\omega_1\rangle$ which is
$\sum_{\alpha<\omega_1}\poset_\alpha$-generic over $V$, use
$X$ and Silver forcing to code the rank zero parts of each $G_\alpha$, i.e.,
define $(\sum_{\alpha<\omega_1}\poset_\alpha) *\name{\posetq}$ by
$(p,q)\in (\sum_{\alpha<\omega_1}\poset_\alpha) *\name{\posetq}$
\sdex{$(\sum_{\alpha<\omega_1}?poset_\alpha) *?name{?posetq}$}

iff

$p\in \sum_{\alpha<\omega_1}\poset_\alpha $ and  $q$ is a finite
set of consistent sentences of the form:
\begin{enumerate}
  \item  $``x \notin \name{U}_n\mbox{''}$ where $x\in X$ or
  \item $``B\subseteq \name{U}_n\mbox{''}$ where $B$ is clopen
   and $n\in\omega$.
\end{enumerate}
with the additional proviso that whenever
$``x_{\alpha,s,n} \notin \name{U}_n\mbox{''}\in q$ then
$s\in\dom({p_\alpha})$ and $p_\alpha(s)\not=n$.  This is a little
stronger than saying $p\forces \check{q}\in\posetq$, but would be true for
a dense set of conditions.

The rank function
$$\rank:(\sum_{\alpha<\omega_1}\poset_\alpha)
 *\name{\posetq}:\arrow\omega_1$$
is defined by
$$\rank(\langle p_\alpha:\alpha<\omega_1\rangle,q)=
\max\{\rank(p_\alpha):\alpha<\omega_1\}$$ which means we ignore
$q$ entirely.

\begin{lemma} \label{rank2}
For every $p\in (\sum_{\alpha<\omega_1}\poset_\alpha) *\name{\posetq}$ and
$\beta\geq 1$ there exists $\hat{p}\in
(\sum_{\alpha<\omega_1}\poset_\alpha) *\name{\posetq}$ such that
\begin{enumerate}
\item $\hat{p}$ is compatible with $p$,
\item $\rank(\hat{p})\leq\beta$, and
\item for any $q\in (\sum_{\alpha<\omega_1}\poset_\alpha) *\name{\posetq}$
if $\rank(q)<\beta$ and
$\hat{p}$ and $q$ are compatible, then $p$ and $q$ are compatible.
\end{enumerate}
\end{lemma}
\proof
Apply Lemma \ref{littlerank} to each $p_\alpha$ to obtain
$\hat{p_\alpha}$ and then let
$$\hat{p}=(\langle \hat{p_\alpha}:\alpha<\omega_1\rangle,q).$$
  This is
still a condition because $\hat{p_\alpha}$ retains all the rank
zero part of $p_\alpha$ which is needed to force
$q\in\name{\posetq}$.
\qed

Let $(\sum_{\alpha<\omega_1}\poset_\alpha) *\name{\posetq}\subseteq\bool$
be a dense subset of the ccc cBa $\bool$.  We show that
$\bool$ is countably generated and $\ord(\bool)=\omega_1$.  A
strange thing about $\omega_1$ is that if one countable set
of generators has order $\omega_1$, then all countable sets of
generators have order $\omega_1$.  This is because
any countable set will be generated by a countable stage.

One set of generators for $\bool$ is
$$C=\{\val{\check{B}\subseteq \name{U_n}}: B\mbox{ clopen }, n\in\omega\}.$$
Note that
$$\val{x\in\intersect_{n\in\omega}U_n}=\prod_{n\in\omega}\val{x\in U_n}
=\prod_{n\in\omega}\sum\{\val{\check{B}\subseteq \name{U_n}}:x\in B\}$$
and also each $\poset_\alpha$ is generated by
$$\{p\in\poset_\alpha: \dom(p)\subseteq \tzero_\alpha\}.$$
We know that for each $\alpha<\omega_1$, $ s\in \tzero_\alpha$ and
$n\in\omega$ if $p=(\langle p_\alpha:\alpha<\omega_1\rangle,q)$
is the condition for which $p_\alpha$ is the function with domain
$\{s\}$, and $p_\alpha(s)=n$, and the rest of $p$ is the trivial condition,
then
$$p=\val{\check{x}_{\alpha,s,n}\in\Intersect_{n\in\omega}\name{U_n}}.$$
From these facts it follows that $C$ generates $\bool$.

It follows from  Lemma \ref{rankbool2} that the order of $C$ is
$\omega_1$.  For any $\beta<\omega_1$ let
$b=(\langle p_\alpha:\alpha<\omega_1\rangle, q)$ be the condition
all of whose components are trivial except for $p_{\beta}$,
and $p_{\beta}$ any the function with domain $\emptyseq$.
Then $b\notin\bsig0\beta(C)$.  Otherwise by Lemma \ref{rankbool2},
there would be some $a\leq b$ with $\rank(a,C)<beta$, but then
$p^a_\beta$ would not have $\emptyseq$ in its domain.

This proves the $\omega_1$ case of Theorem \ref{boolord}.

\sector{Luzin sets}
In this section we use Luzin sets and generalized $I$-Luzin sets
to construct separable metric spaces of various Borel orders.
Before doing so we review some standard material on
the property of Baire.

Given a topological space $X$ the $\sigma$-ideal of meager sets
is defined as follows. $Y\subseteq X$ is nowhere dense iff
the interior of the closure of $Y$ is empty,
i.e, $\interior(\cl(Y))=\emptyset$. A subset
of $X$ is meager iff it is the countable union of nowhere dense
sets.   The Baire category Theorem is the statement that
nonempty open subsets of
compact Hausdorff spaces or completely metrizable spaces are
not meager.  A subset $B$ of $X$ has the Baire property
iff there exists $U$ open such that $B\Delta U$ is meager.

\begin{theorem} \label{borelmod}
  (Baire) The family of sets with the Baire property forms
  a $\sigma$-field.
\end{theorem}
\proof
If $B\Delta U$ is meager
where $U$ is open, then
$$B\Delta\cl(U)=(B\setminus\cl(U))\union (\cl(U)\setminus B)$$
and $(B\setminus\cl(U))\subseteq B\setminus U$ is meager and
$(\cl(U)\setminus B)\subseteq (U\setminus B)\union (\cl(U)\setminus U)$
is meager because $\cl(U)\setminus U$ is nowhere dense.
Therefore,
$$\comp{B}\;\;\Delta\comp{\cl(U)}=B\Delta\cl(U)$$
is meager.

If $B_n\Delta U_n$ is meager for each $n$, then
$$(\Union_{n\in\omega}B_n)\Delta (\Union_{n\in\omega}U_n)\subseteq
\Union_{n\in\omega}B_n\Delta U_n$$
is meager.
\qed

Hence every Borel set has the property of Baire.
\begin{theorem}
  Suppose that every nonempty open subset of $X$ is nonmeager,
  then $\bool=\borel(X)/\meager(X)$ is a cBa.   \sdex{$\borel(X)/\meager(X)$}
\end{theorem}                                   \sdex{$\meager(X)$}
\proof
It is enough to show that it is complete.  Suppose
$\Gamma\subseteq\bool$ is arbitrary.
Let $\cal U$ be a family of open sets such
that
$$\Gamma=\{[U]_{\meager(X)}: U\in {\cal U}\}.$$
Let $V=\Union{\cal U}$ and we claim that
$[V]$ is the minimal upper bound of $\Gamma$ in
$\bool$.  Clearly it is an upper bound. Suppose
$[W]$ is any upper bound for $\Gamma$ with $W$ open.
So $U\setminus W$ is meager for every $U\in{\cal U}$.
We need to show that $V\setminus W$ is meager (so
$[V]\leq [W]$).
$V\setminus W\subseteq \cl(V)\setminus W$ and if the latter
is  not nowhere dense, then there exists $P$ a nonempty open
set with $P\subseteq \cl(V)\setminus W$. Since
$V=\Union{\cal U}$ we may assume that that there exists
$U\in{\cal U}$ with $P\subseteq U$. But $P$ is a nonempty open
set and $[P]\leq [W]$ so it is impossible for $P$ to be
disjoint from $W$.
\qed

We say that $X\subseteq\cantorset$ is a \dex{super Luzin set}
iff for every Borel  set $B$
the set $X\intersect B$ is countable iff B is meager.
It is easy to see that if $X$ is an ordinary
Luzin set, then in some basic clopen set $C$ it is a super Luzin set
relative to $C$.  Also since $\bairespace$ can be obtained by
deleting countably many points from $\cantorspace$ it is clear
that having a Luzin set for one is equivalent to having it for
the other.  With a little more work it can be seen that it is
equivalent to having one for any completely metrizable separable metric
space without isolated points.

The generic set
of Cohen reals in the Cohen real model is a Luzin set.
Let $\fin(\kappa,2)$ be the partial order of finite partial
functions from $\kappa$ into $2$. \sdex{$\fin(\kappa,2)$}
If $G$ is $\fin(\kappa,2)$-generic over V and
for each $\alpha<\kappa$ we define $x_\alpha$ by
$x_\alpha(n)=G(\omega*\alpha+n)$, then
$X=\{x_\alpha:\alpha<\kappa\}$ is a Luzin set in $V[G]$.

\begin{theorem} \label{luzinset}
  (Miller \site{borelhier})
  If there exists a Luzin set in $\bairespace$, then
  for every $\alpha$ with $3\leq\alpha<\omega_1$ there exists
  $Y\subseteq\bairespace$ with $\ord(Y)=\alpha$.
\end{theorem}
\proof

Let $T_\alpha$ \sdex{$T_\alpha$}
be the nice $\alpha$-tree used in the definition of
$$\poset_{\alpha}=\{p:\;\; p:D\arrow \omega, D\in [T_\alpha]^{<\omega},
\forall s,s\concat n\in D\;\;p(s)\not= p(s\concat n)\}.$$
Let \sdex{$Q_\alpha$} $Q_\alpha$ be the closed subspace of $\omega^{T_\alpha}$
$$Q_\alpha=\{x\in \omega^{T_\alpha}:
\forall s,s\concat n\in T_\alpha\;\;x(s)\not= x(s\concat n).$$
For each $p\in\poset_{\alpha}$ we get a basic clopen set
$$[p]=\{x\in Q_\alpha: p\subseteq x\}.$$
It easy to check that $Q_\alpha$ is homeomorphic to $\bairespace$.
Hence there exists
a super Luzin set $X\subseteq Q_\alpha$.  Consider the map
$r:Q_\alpha\arrow \omega^{\tzero_\alpha}$ defined by restriction, i.e.,
$r(x)=x\res \tzero_\alpha$.  Note that by Lemma \ref{denseset}
there exists a countable sequence of dense  open subsets of
$Q_\alpha$, $\langle D_n:n\in\omega\rangle$, such that
$r$ is one-to-one on $\Intersect_{n\in\omega}D_n$.
Since $\Intersect_{n\in\omega}D_n$ is a comeager set in $Q_\alpha$ and
$X$ is Luzin we may assume that
$$X\subseteq\Intersect_{n\in\omega}D_n.$$

$Y$ is just the image of $X$ under $r$.  (So, in fact,
$Y$ is the one-to-one  continuous image of  a Luzin set.)
An equivalent way to view
$Y$ is just to imagine $X$ with the topology given by
$$\base=\{[p]: p\in\poset_\alpha, \dom(p)\subseteq \tzero_\alpha\}.$$
We know by Lemma \ref{rankbool2} that
$$\ord\{[B]:B\in\base\}=\alpha+1$$
as a subset of
$\borel(Q_\alpha)/\meager(Q_\alpha)$ which means that:

$\alpha+1$ is
minimal  such that for every
$B\in\borel(Q_\alpha)$ there exists a $\bsig0{\alpha+1}(\base)$ set $A$
such that such that $B\Delta A$ is meager in $Q_\alpha$.

This translates (since $X$ is super-Luzin) to:

$\alpha+1$ is
minimal  such that for every
$B\in\borel(Q_\alpha)$ there exists a $\bsig0{\alpha+1}(\base)$ set $A$
such that such that $(B\Delta A)\intersect X$ is countable.

Which means for $Y$ that:

$\alpha+1$ is
minimal  such that for every
$B\in\borel(Y)$ there exists a $\bsig0{\alpha+1}(Y)$ set $A$
such that such that $B\Delta A$ is countable.

But since countable subsets of $Y$ are $\bsig02$ and $\alpha>2$,
this means $\ord(Y)=\alpha+1$.

To get $Y$ of order $\lambda$ for
a limit $\lambda<\omega_1$ just take a clopen separated
union of sets whose order increases to $\lambda$.

\qed

Now we clean up a loose end from Miller \site{borelhier}.
In that paper we had shown that assuming MA for every $\alpha<\omega_1$
there exists a separable metric space $X$ with
$\alpha\leq \ord(X)\leq\alpha+2$ or something silly like that.
Shortly afterwards, Fremlin supplied the missing arguments to
show the following.

\begin{theorem} \label{frem}
  (Fremlin \site{fremlin}) MA implies that for every
  $\alpha$ with $2\leq\alpha\leq\omega_1$ there exists
  a second countable Hausdorff space $X$
  with $\ord(X)=\alpha$.
\end{theorem}
\proof

Since the union of less than continuum many meager sets is meager, the
Mahlo construction \ref{mah} gives us a set $X\subseteq Q_\alpha$
of cardinality $\continuum$ such that for every Borel set
$B\in\borel(Q_\alpha)$ we have that $B$ is meager iff
$B\intersect X$ has cardinality less than $\continuum$.

Letting $\base$ be defined as in the proof of Theorem
\ref{luzinset} we see that:

$\alpha+1$ is
minimal such that for every
$B\in\borel(Q_\alpha)$ there exists a $\bsig0{\alpha+1}(\base)$ set $A$
such that such that $(B\Delta A)\intersect X$ has cardinality less
than $\continuum$.

What we need to see to complete the proof is that:

for
every $Z\subseteq X$ of cardinality less than $\continuum$ there
exists a $\bsig02(\base)$ set $F$ such that $F\intersect X=Z$.

\begin{lemma} (MA)
  For any $Z\subseteq Q_\alpha$ of cardinality less than $\continuum$,
  there exists $\langle D_n:n\in\omega\rangle$ such that:
  \begin{enumerate}
    \item $D_n$ is predense in $\poset_\alpha$,
    \item $p\in D_n$ implies $\dom(p)\subseteq\tzero_\alpha$, and
    \item $Z\intersect \Intersect_{n\in\omega}\Union_{s\in D_n}[s]=\emptyset$.
  \end{enumerate}
\end{lemma}

\proof
Force with the following poset
$$P=\{(F,\langle p_n:n<N\rangle): F\in [Z]^{<\omega},\; N<\omega,\;
\dom(p)\in [\tzero_\alpha]^{<\omega}\}$$
where $(F,\langle p_n:n<N\rangle)\leq (H,\langle q_n:n<M\rangle)$
iff $F\supseteq H$, $N\geq M$, $p_n=q_n$ for $n<M$, and for each
$x\in H$ and $M\leq n<N$ we have $x\notin [p_n]$.
Since this forcing is ccc we can apply MA with the appropriate
choice of family of dense sets to get $D_n=\{p_m:m>n\}$
to do the job.
\qed

By applying the Lemma we get that for
every $Z\subseteq X$ of cardinality less than $\continuum$ there
exists a $\bsig02(\base)$ set $F$ which is meager in
$Q_\alpha$ and such that $Z\subseteq F\intersect X$.  But since
$F$ is meager we know
$F\intersect X$ has cardinality less than $\continuum$.
By Theorem \ref{silverlem} every subset of $r(F\intersect X)$
is a relative $\bsig02$ in $Y$, so
there exists an $F_0$ a $\bsig02(\base)$ set such that
 $Z=(F\intersect X)\intersect F_0$.
This proves Theorem \ref{frem}.
\qed

\sector{Cohen real model}

I have long wondered if there exists an uncountable \label{cohen}
separable metric space of order $2$ in the Cohen real model.
I thought there weren't any.
We already
know from Theorem \ref{luzinset} that since there is an uncountable
Luzin set in Cohen real model that for every $\alpha$ with
$3\leq\alpha\leq\omega_1$ there is an
uncountable separable metric space $X$ with $\ord(X)=\alpha$.

\begin{theorem} \label{cohenmodel}
Suppose $G$ is $\fin(\kappa,2)$-generic over $V$
where $\kappa\geq\omega_1$.
Then in
$V[G]$ there is a
separable metric space $X$ of cardinality $\omega_1$
with $\ord(X)=2$.
\end{theorem}
\proof
We may assume that $\kappa=\omega_1$. This is because
$\fin(\kappa,2)\cross\fin(\omega_1,2)$ is isomorphic
to $\fin(\kappa,2)$ and so by the product lemma we may
replace $V$ by $V[H]$ where $(H,G)$ is
$\fin(\kappa,2)\cross\fin(\omega_1,2)$-generic over $V$.

We are going to use the fact that forcing with
$\fin(\omega_1,2)$ is equivalent to any finite support
$\omega_1$ iteration of countable posets.
The main idea of the proof is to construct an Aronszajn tree of perfect sets,
a technique first used by Todorcevic (see Galvin and Miller \site{galmil}).
We construct an Aronszajn tree $(A,\fleq)$ and a family of
\sdex{Aronszajn tree}
perfect sets $([T_s]:s\in A)$ such that $\supseteq$ is the
same order as $\fleq$.  We will then show that if $X=\{x_s:s\in A\}$
is such that $x_s\in [T_s]$, then the order of $X$ is $2$.

In order to insure the construction can keep going at limit ordinals
we will need to use a fusion argument. Recall that a perfect set
corresponds to the infinite branches $[T]$  of a \dex{perfect tree}
$T\subseteq 2^{<\omega}$, i.e., a tree with the property that
for every $s\in T$ there exist a $t\in T$ such that both
$t\concat 0\in T$ and $t\concat 1\in T$.  Such a $T$ is called
a \dex{splitting node} of $T$. There is a natural correspondence of
the splitting nodes of a perfect tree $T$ and $2^{<\omega}$.

Given two perfect trees $T$ and $T^\prime$ and $n\in\omega$ define
$T\leq_n T^\prime$ \sdex{$T\leq_n T^\prime$}
iff $T\subseteq T^\prime$ and the first $2^{<n}$ splitting nodes
of $T$ remain in  $T^\prime$.

\begin{lemma}\label{fusion}
  (\dex{Fusion}) Suppose $(T_n:n\in\omega)$ is a sequence of
  perfect sets such that $T_{n+1}\leq_n T_n$ for every $n\in\omega$.
  Then $T=\Intersect_{n\in\omega}T_n$ is a perfect tree and
  $T\leq_n T_n$ for every $n\in\omega$.
\end{lemma}
\proof
If $T=\Intersect_{n<\omega}T_n$, then $T$ is a perfect
tree because the first $2^{<n}$ splitting nodes of $T_n$ are in
$T_m$ for every $m>n$ and thus in $T$.
\qed

By identifying $\fin(\omega_1,2)$ with
$\sum_{\alpha<\omega_1}\fin(\omega,2)$ we may assume that
$$G=\langle G_\alpha:\alpha<\omega_1\rangle$$
where $G_\beta$ is
$\fin(\omega,2)$-generic over $V[G_\alpha:\alpha<\beta]$
for each $\beta<\omega_1$.

Given an Aronszajn tree $A$ we let $A_\alpha$ \sdex{$A_\alpha$}
\sdex{$A_{<\alpha}$}  be the nodes of $A$ at
level $\alpha$, i.e.
$$A_\alpha=\{s\in A: \{t\in A: t\fless s\}\mbox{ has order type } \alpha\}$$
and
$$A_{<\alpha}=\Union_{\beta<\alpha} A_\beta.$$
We use $\langle G_\alpha:\alpha<\omega_1\rangle$ to construct
an Aronszajn tree $(A,\fleq)$ and a family of
perfect sets $([T_s]:s\in A)$ such that
\begin{enumerate}
\item $s\fleq t$ implies $T_s\supseteq T_t$,
\item if $s$ and $t$ are distinct elements of $A_\alpha$, then
$[T_s]$ and $[T_t]$ are disjoint,
\item every $s\in A_\alpha$ has infinitely many distinct extensions
 in $A_{\alpha+1}$,
\item for each $s\in A_{<\alpha}$ and
$n<\omega$ there exists $t\in A_\alpha$ such that $T_t\leq_n T_s$,
\item for each $s\in A_\alpha$ and $t\in A_{\alpha+1}$ with
$s\fless t$, we have that
$[T_t]$ is a generic perfect subset of $[T_s]$ obtained by using $G_\alpha$
(explained below in Case 2), and
\item $\{T_s:s\in A_{<\alpha}\}\in V[G_\beta:\beta<\alpha]$.
\end{enumerate}

The first three items simply say that $\{[T_s] : s\in A\}$
and its ordering by $\subseteq$ determines
$(A,\sleq)$, so what we really have here is an Aronszajn tree of perfect sets.
Item (4) is there in order to allow the construction to proceed at
limits levels.

Item (5) is what we do a successor levels and guarantees the set we
are building has order $2$.  Item (6) is a consequence of the construction
and would be true for a closed unbounded set of ordinals no matter what we
did anyway.

Here are the details of our construction.

\bigskip
\par\noindent Case 1. $\alpha$ a limit ordinal.

The construction is done uniformly enough so that we already have
that $\{T_s:s\in A_{<\alpha}\}\in V[G_\beta:\beta<\alpha]$.
Working in $V[G_\beta:\beta<\alpha]$ choose a sequence $\alpha_n$
for $n\in\omega$ which strictly increases to $\alpha$. Given
any $s_n\in A_{\alpha_n}$ we can choose by inductive hypothesis
a sequence $s_m\in A_{\alpha_m}$ for $m\geq n$ such that
$$T_{s_{m+1}}\leq_m T_{s_m}.$$
If $T=\Intersect_{m>n} T_{s_m}$, then by Lemma \ref{fusion} we
have that $T\leq_nT_{s_n}$. Now let $\{T_t:t\in A_\alpha\}$
be a countable collection of perfect trees so that for
every $n$ and $s\in A_{\alpha_n}$ there exists $t\in A_\alpha$
with $T_t\leq_n T_s$.  This implies item (4) because
for any $s\in A_{<\alpha}$ and $n<\omega$ there exists some
$m\geq n$ with $s\in A_{<\alpha_m}$ hence by inductive hypothesis
there exists $\hat{s}\in A_{\alpha_m}$ with $T_{\hat{s}}\leq_n T_s$
and by construction there exists $t\in A_\alpha$
with $T_t\leq_m T_{\hat{s}}$ and so $T_t\leq_n T_s$ as desired.

\bigskip
\par\noindent Case 2. Successor stages.

Suppose we already have constructed
$$\{T_s:s\in A_{<\alpha+1}\}\in V[G_\beta:\beta<\alpha+1].$$
Given a perfect tree $T\subseteq 2^{<\omega}$ define the countable
partial order $\poset(T)$ as follows. $p\in\poset(T)$ iff $p$
\sdex{$?poset(T)$}
is a finite subtree of $T$ and $p\leq q$ iff $p\supseteq q$ and
$p$ is an end extension of $q$, i.e., every new node of $p$ extends
a terminal node of $q$.  It is easy to see that if $G$ is $\poset(T)$-generic
over a model $M$, then
$$T_G=\Union\{p: p\in G\}$$
is a perfect subtree of $T$.  Furthermore, for any $D\subseteq [T]$
dense open in $[T]$ and
coded in $M$,  $[T_G]\subseteq D$.
i.e., the branches of $T_G$ are Cohen reals (relative to $T$) over $M$.
This means that for any Borel set $B\subseteq [T]$ coded in $M$,
there exists an clopen set $C\in M$ such that
$$C\intersect [T_G]=B\intersect [T_G].$$
To see why this is true let $p\in\poset(T)$ and $B$ Borel. Since
$B$ has the Baire property relative to $[T]$
by extending each terminal node of $p$,
if necessary, we can obtain $q\geq p$ such that for every terminal
node $s$ of $q$ either  $[s]\intersect B$ is meager in $[T]$
or $[s]\intersect {B}$
is comeager in $[T]\intersect [s]$.
If we let
$C$ be union of all $[s]$ for $s$ a terminal node of $q$
such that $[s]\intersect {B}$
is comeager in $[T]\intersect [s]$, then
$$q\forces B\intersect T_G= C\intersect T_G.$$
To get $T_G\leq_n T$ we could instead force with
$$\poset(T,n)=\{p\in\poset(T):
p \mbox{ end extends the first $2^{<n}$ splitting nodes of }
T\}.$$

Finally to determine $A_{\alpha+1}$ consider
$$\sum\{\poset(T_s,m):s\in A_{\alpha}, m\in\omega\}.$$
This poset is countable and hence $G_{\alpha+1}$ determines
a sequence
      $$\langle T_{s,m}: s\in A_\alpha, m\in\omega\rangle$$
of generic perfect trees such that $T_{s,m}\leq_m T_s$.
Note that genericity also guarantees that corresponding
perfect sets will be disjoint.  We define $A_{\alpha+1}$
to be this set of generic trees.

\medskip

This ends the construction.

\medskip

By taking generic perfect sets at successor steps we have guaranteed the
following. For any Borel set $B$ coded in
$V[G_\beta:\beta<\alpha+1]$ and $T_t$ for $t\in A_{\alpha+1}$
there exists a clopen set $C_t$ such that
$$C_t\intersect [T_t]= B\intersect [T_t].$$
Suppose $X=\{x_s: s\in A\}$ is such that
$x_s\in [T_s]$ for every $s\in A$.  Then $X$ has order $2$.
To verify this, let $B\subseteq\cantorspace$ be any Borel set.
By ccc there exists a countable $\alpha$ such that $B$ is coded
in $V[G_\beta:\beta<\alpha+1]$.  Hence,
$$B\intersect \Union_{t\in A_{\alpha+1}} [T_t]=
  \Union_{t\in A_{\alpha+1}} (C_t\intersect [T_t]).$$
Hence $B\intersect X$ is equal to a $\bsig02$ set intersected
$X$:
$$ X\intersect \Union_{t\in A_{\alpha+1}} (C_t\intersect [T_t])$$
 union a countable set:
$$(B\intersect X)\setminus \Union_{t\in A_{\alpha+1}} [T_t]$$
and therefore $B\intersect X$ is $\bsig02$ in $X$.
\qed

Another way to get a space of order $2$ is to use the following
argument. If the ground model satisfies CH, then there exists
a Sierpinski set.  Such a set has order $2$
(see Theorem \ref{pop}) in $V$ and therefore by the next
theorem it has order 2 in  $V[G]$.  It also follows from the
next theorem that
if $X=\cantorspace\intersect V$,
then $X$ has order $\omega_1$ in $V[G]$.  Consequently, in what
I think of as ``the Cohen real model'', i.e. the
model obtained by adding $\omega_2$ Cohen reals to a model of
CH,  there are separable metric spaces of cardinality $\omega_1$ and
order $\alpha$ for every $\alpha$ with $2\leq\alpha\leq\omega_1$.

\begin{theorem} \label{cohonord}
Suppose $G$ is $\fin(\kappa,2)$-generic over $V$
where $\kappa\geq\omega$ and
$V\models ``\ord(X)=\alpha$''.
Then $V[G]\models ``\ord(X)=\alpha$''.
\end{theorem}

By the usual ccc arguments it is clearly
enough to prove the Theorem  for $\fin(\omega,2)$.
To prove it we will need the following
lemma.

\begin{lemma} \label{cohlem}
  (Kunen, see \site{kunmil}) Suppose $p\in \fin(\omega,2)$,
  $X$ is a second countable Hausdorff space in $V$, and
  $\name{B}$ is a name such that
  $$p\forces \name{B}\subseteq\check{X}\mbox{ is a $\bpi0\alpha$-set}. $$
  Then the set
  $$\{x\in X:p\forces \check{x}\in\name{B}\}$$
  is a $\bpi0\alpha$-set in $X$.
\end{lemma}
\proof
This is proved by induction on $\alpha$.

For $\alpha=1$ let
$\base\in V$ be a countable base for the closed subsets of $X$, i.e.,
every closed set is the intersection of elements of $\base$.
Suppose $p\forces$``$ \name{B}$ is a closed set in $\check{X}$''.
Then for every $x\in X$
$p\forces$``$ \check{x}\in\name{B}$'' iff
for every $q\leq p$ and
for every $C\in\base$ if $q\forces$``$\name{B}\subseteq\check{C}$'',
then $x\in C$. But
$$\{x\in X: \forall q\leq p\;\forall B\in \base\;(q\forces
\mbox{``}\name{B}\subseteq\check{C}\mbox{''}\arrow x\in C)\}$$
is closed in  $X$.

Now suppose $\alpha>1$ and $p\forces\name{B}\in \bpi0\alpha(X)$.
Let $\beta_n$ be a sequence which is either constantly $\alpha-1$
if $\alpha$ is a successor or which is unbounded in $\alpha$ if
$\alpha$ is a limit.
By the usual forcing facts there exists a sequence of names
$\langle B_n:n\in\omega\rangle$ such that
for each $n$, $$p\forces B_n\in\bpi0{\beta_n},$$ and
$$p\forces B=\Intersect_{n<\omega}\comp{B_n}.$$
Then for every $x\in X$
$$p\forces \check{x}\in B$$ iff
$$\forall n\in\omega\;\; p\forces x\in\comp{B_n}$$ iff
$$\forall  n\in\omega\;\forall q\leq p\;\; q\not\forces x\in{B_n}.$$
Consequently,
$$\{x\in X: p\forces x\in \name{B}\}=\Intersect_{n\in\omega}
\Intersect_{q\leq p} \comp{\{x: q\forces \check{x}\in\name{B_m}\}}.$$
\qed

Now let us prove the Theorem.  Suppose
$V\models$``$\ord(X)=\alpha$''.  Then in $V[G]$ for any Borel
set $B\in\borel(X)$
$$B=\Union_{p\in G}\{x\in X: p\forces \check{x}\in \name{B}\}.$$
By the lemma, each of the sets
$\{x\in X: p\forces \check{x}\in \name{B}\}$ is a Borel set in $V$, and since
$\ord(X)=\alpha$, it is a $\bsig0\alpha$ set.  Hence, it
follows that $B$ is a $\bsig0\alpha$ set.
So,  $V[G]\models\ord(X)\leq\alpha$. To see that
$\ord(X)\geq\alpha$ let $\beta<\alpha$ and suppose
in $V$ the set $A\subseteq X$ is $\bsig0\beta$ but not
$\bpi0\beta$. This must remain true in $V[G]$ otherwise
there exists a $p\in G$ such
that
$$p\forces\mbox{``$\check{A}$ is $\bpi0\beta$''}$$
but by the lemma
$$\{x\in X:p\forces \check{x}\in \check{A}\}=A$$
is $\bpi0\beta$ which is a contradiction.
\qed

Part of this argument is similar to one used by Judah and Shelah
\site{qsetzero}
who showed that it is consistent to have a Q-set which does not
have strong measure zero.

It is natural to ask if there are spaces of order $2$
of higher cardinality.

\begin{theorem} \label{cohenmodel2}
Suppose $G$ is $\fin(\kappa,2)$-generic over $V$
where $V$ is a model of CH and $\kappa\geq\omega_2$.
Then in
$V[G]$  for every separable metric space $X$ with
  $|X|>\omega_1$, we have $\ord(X)\geq 3$.
\end{theorem}
\proof

This will follow easily from the next lemma.

\begin{lemma} (Miller \site{maponto})  \label{mapping}
Suppose $G$ is $\fin(\kappa,2)$-generic over $V$
where $V$ is a model of CH and $\kappa\geq\omega_2$.
Then  $V[G]$ models that for every
$X\subseteq\cantorspace$ with $|X|=\omega_2$ there
exists a Luzin set $Y\subseteq\cantorspace$ and a
one-to-one continuous function $f:Y\arrow X$.
\end{lemma}
\proof
Let $\langle \tau_\alpha:\alpha<\omega_2\rangle$ be a sequence
of names
for distinct elements of $X$.  For each
$\alpha$ and $n$ choose a maximal antichain $A_n^\alpha\union B_n^\alpha$
such that
$$p\forces \tau_\alpha(n)=0 \mbox{ for each $p\in A_n^\alpha$ and}$$
$$p\forces \tau_\alpha(n)=1 \mbox{ for each $p\in B_n^\alpha$}.$$
Let $X_\alpha\subseteq\kappa$ be union of domains of elements
from $\Union_{n\in\omega}A^\alpha_n\union B^\alpha_n$.
Since each $X_\alpha$ is
countable we may as well assume that the
$X_\alpha$'s form a $\Delta$-system, i.e. there exists $R$ such
that $X_\alpha\intersect X_\beta=R$ for every $\alpha\not=\beta$.
We can assume that $R$ is the empty set.  The
reason is we can just replace $A^\alpha_n$ by
$$\hat{A}^\alpha_n=\{p\res (\comp{R}):p\in A^\alpha_n \rmand
p\res R\in G\}$$
and similarly for $B^\alpha_n$. Then let $V[G\res R]$ be the new
ground model.

Let
$$\langle j_\alpha:X_\alpha\arrow \omega:\alpha<\omega_2\rangle$$
be a sequence of bijections in the ground model.
Each
$j_\alpha$ extends to an order preserving map from $\fin(X_\alpha,2)$ to
$\fin(\omega,2)$.  By CH, we may as well assume that
there exists a  single sequence, $\langle (A_n,B_n):n\in\omega\rangle$,
such that every  $j_\alpha$ maps
$\langle A^\alpha_n,B^\alpha_n:n\in\omega\rangle$
to $\langle (A_n,B_n):n\in\omega\rangle$.

The Luzin set is
$Y=\{y_\alpha:\alpha<\omega_2\}$
where $y_\alpha(n)=G(j_\alpha^{-1}(n))$.  The continuous
function, $f$, is the map determined by
$\langle (A_n,B_n):n\in\omega\rangle$:
$$f(x)(n)=0 \rmiff \exists m \; \; x\res m\in A_n.$$
This proves the Lemma.
\qed

If  $f:Y\arrow X$ is one-to-one continuous function from a Luzin
set $Y$, then $\ord(X)\geq 3$.  To see this assume that
$Y$ is dense and let
$D\subseteq Y$ be a countable dense subset of $Y$.  Then
$D$ is not $G_\delta$ in $Y$.  This is because any $G_\delta$
set containing $D$ is comeager and therefore must meet $Y$ in an uncountable
set.  But note that $f(D)$ is a countable set which cannot be
$G_\delta$ in $X$, because $f^{-1}(f(D))$ would be
$G_\delta$ in $Y$ and since $f$ is one-to-one we have $D=f^{-1}(f(D))$.
This proves the Theorem.
\qed

It is natural to ask about the cardinalities of sets of various orders
in this model.  But note that there is a trivial way to get a large
set of order $\beta$.  Take a clopen separated union of a large
Luzin set (which has order 3) and a set of size $\omega_1$ with
order $\beta$.  One possible way to strengthen the notion of
order is to say that a space $X$ of cardinality
$\kappa$ has essential order $\beta$ iff every nonempty
open subset of $X$ has order $\beta$ and cardinality $\kappa$.
But this is also open to a simple trick of combining a small set
of order $\beta$ with a large set of small order.  For example,
let $X\subseteq\cantorspace$ be a clopen separated
union of a Luzin set of cardinality $\kappa$ and set of cardinality
$\omega_1$ of order $\beta\geq 3$.
Let $\langle P_n:n\in\omega\rangle$ be a sequence of disjoint nowhere
dense perfect subsets of $\cantorspace$ with the property that for
every $s\in 2^{<\omega}$ there exists $n$ with $P_n\subseteq [s]$.
Let $X_n\subseteq P_n$ be a homeomorphic copy of $X$ for each $n<\omega$.
Then $\Union_{n\in\omega} X_n$ is a set of cardinality $\kappa$ which
has essential order $\beta$.

With this cheat in mind let us define a stronger notion of order.
A separable metric space $X$ has \dex{hereditary order} $\beta$ iff
every uncountable $Y\subseteq X$ has order $\beta$.  We begin
with a stronger version of Theorem \ref{luzinset}.

\begin{theorem}\label{cohenluz}
  If there exists a Luzin set $X$ of cardinality $\kappa$, then for
  every $\alpha$ with $2<\alpha<\omega_1$ there exists a separable
  metric space $Y$ of cardinality $\kappa$ which is hereditarily of
  order $\alpha$.
\end{theorem}
\proof
This is a slight modification of the proof of Theorem \ref{luzinset}.
Let $\posetq_\alpha$ \sdex{$?posetq_\alpha$} be the following partial
order. Let
$\langle \alpha_n:n\in\omega\rangle$ be a sequence such that
if $\alpha$ is a limit ordinal, then $\alpha_n$ is a cofinal increasing
sequence in $\alpha$ and if $\alpha=\beta+1$ then $\alpha_n=\beta$
for every $n$.

The rest of the proof is same except we use $\posetq_{\alpha+1}$ instead
of $\poset_\alpha$ for successors and $\posetq_\alpha$ for limit
$\alpha$ instead of taking a clopen separated union.  By using the direct
sum (even in the successor case) we get a stronger property for
the order.  Let
$$\hat{Q}_\alpha=\prod{Q_{\alpha_n}}$$ be the closed subspace of
\sdex{$?hat{Q}_\alpha$}
$$\prod_{n\in\omega}\omega^{T_{\alpha_n}}$$
and let $\base$ be the collection of clopen subsets of $Q_\alpha$ which
are given by rank zero conditions of $\posetq(\alpha)$, i.e.,
all rectangles of the form $\prod_{n\in\omega}[p_n]$ such
that
$p_n\in \posetq_{\alpha_n}$ with $\dom(p)\subseteq\tzero_\alpha$
and $p_n$ the trivial condition for all but finitely many $n$.

As in the proof of Theorem \ref{luzinset} we get that
 the order of $\{[B]:B\in\base\}$ as a subset
of $\borel(\hat{Q}_\alpha)/\meager(\hat{Q}_\alpha)$
is $\alpha$.  Because
we took the direct sum we get the
stronger property that for any nonempty clopen set $C$ in $\hat{Q}_\alpha$
the order of $\{[B\intersect C]:B\in\base\}$ is $\alpha$.

But know given $X$ a Luzin set in $\hat{Q}_\alpha$ we know that for any
uncountable $Y\subseteq X$ there is a nonempty clopen set
$C\subseteq\hat{Q}_\alpha$ such that $Y\intersect C$ is a super-Luzin
set relative to $C$.  (The \dex{accumulation points} of $Y$, the set of
all points all of whose neighborhoods contain uncountably many points of
$Y$, is closed and uncountable, therefore must have nonempty interior.)
If $C$ is a nonempty clopen set in the interior of the accumulation points
of $Y$, then since $\{[B\intersect C]:B\in\base\}$ is $\alpha$,
we have by the proof of Theorem \ref{luzinset}, that the order of $Y$ is
$\alpha$.
\qed

\begin{theorem} \label{cohen1}
Suppose that in $V$ there is a separable metric space,
$X$, with hereditary order $\beta$ for some $\beta\leq\omega_1$.
Let
$G$ be $\fin(\kappa,2)$-generic over $V$ for any
$\kappa\geq\omega$.
Then in $V[G]$ the space $X$ has hereditary order $\beta$.
\end{theorem}
\proof
In $V[G]$ let $Y\subseteq X$ be uncountable.
For contradiction, suppose that
$$p\forces \ord(\name{Y})\leq\alpha\rmand |\name{Y}|=\omega_1$$
for some $p\in \fin(\kappa,2)$ and $\alpha<\beta$.
Working in $V$ by the usual
$\Delta$-system argument  we can get $q\leq p$ and
$$\langle p_x : x\in A\rangle$$
for some $A\in [X]^{\omega_1}$ such that
and $p_x\leq q$ and
$$p_x\forces \check{x}\in \name{Y}$$
for each $x\in A$ and
$$dom(p_x)\intersect dom(p_y)=dom(q)$$
for  distinct $x$ and $y$ in $A$.
Since $A$ is uncountable we know that in $V$ the order of $A$ is
$\omega_1$.  Consequently, there exists $R\subseteq A$ which is
$\bsig0\alpha(A)$ but not $\bpi0\alpha(A)$. We claim that in $V[G]$
the set $R\intersect Y$ is not $\bpi0\alpha(Y)$.  If not, there
exists $r\leq q$ and $\name{S}$ such that
$$r\forces \mbox{``}\name{Y}\intersect R= \name{Y}\intersect \name{S}
\rmand \name{S}\in \bpi0\alpha(A)\mbox{''}.$$
Since Borel sets are coded by reals there exists
$\Gamma\in[\kappa]^\omega\intersect V$
such that for any $x\in A$ the statement ``$\check{x}\in \name{S}$''
is decided by conditions in $\fin(\Gamma,2)$ and also
let $\Gamma$ be large enough to contain the domain of $r$.

Define
$$T=\{x\in A: q\forces \check{x}\in \name{S}\}.$$
According to Lemma \ref{cohlem} the set $T$ is $\bpi0\alpha(A)$.
Consequently, (assuming $\alpha\geq 3$) there are uncountably
many $x\in A$ with  $x \in R\Delta T$.  Choose such an $x$ which
also has the property that $dom(p_x)\setminus dom(q)$ is disjoint
from $\Gamma$.  This can be done since the $p_x$ form a
$\Delta$ system.  But now, if $x\in T\setminus R$, then
$$r\union p_x\forces \mbox{``}\check{x}\in \name{Y}\intersect \name{S}
\rmand x\notin \name{Y}\intersect \check{R}\mbox{''}.$$
On the other hand, if $x\in R\setminus T$, then
there exists $\hat{r}\leq r$ in $\fin(\Gamma,2)$ such that
$$\hat{r}\forces \check{x}\notin \name{S}$$
and consequently,
$$\hat{r}\union p_x\forces
\mbox{``}\check{x}\notin \name{Y}\intersect \name{S}
\rmand x\in \name{Y}\intersect \check{R}\mbox{''}.$$
Either way we get a contradiction and the result is proved.
\qed

\begin{theorem}\label{cohen2}
(CH) There exists $X\subseteq\cantorspace$ such that
$X$ has hereditary order $\omega_1$.
\end{theorem}
\proof
By Theorem \ref{boolord} there exists a countably generated
ccc cBa $\bool$ which has order $\omega_1$. For any $b\in\bool$
with $b\not=0$ let $\ord(b)$ be the order of the boolean algebra
you get by looking only at $\{c\in\bool: c\leq b\}$.  Note that
in fact $\bool$ has the property that for every $b\in\bool$
we have $\ord(b)=\omega_1$.  Alternatively, it easy to show that
any ccc cBa of order $\omega_1$ would have to contain
an element $b$ such that every $c\leq b$ has order $\omega_1$.

By the proof
of the Sikorski-Loomis Theorem \ref{sikorthm} we know that
$\bool$ is isomorphic to $\borel(Q)/\meager(Q)$ where
$Q$ is a ccc compact Hausdorff space with a basis of
cardinality continuum.

Since $Q$ has ccc, every open dense set
contains an open dense set which is a countable union of
basic open sets. Consequently, by using CH, there exists a
family $\famly$ of meager sets
with $|\famly|=\omega_1$ such that
every meager set is a subset of one in $\famly$.
Also note that for any nonmeager Borel set $B$ in $Q$ there
exists a basic open set $C$ and $F\in\famly$ with
$C\setminus F\subseteq B$.
Hence by Mahlo's construction (Theorem \ref{mah}) there
exists a set $X\subseteq Q$ with the property that for
any Borel subset $B$ of $Q$
$$|B\intersect X|\leq\omega \rmiff B \mbox{ meager.} $$

Let $\base$ be a countable field of clopen subsets of
$Q$ such that
$$\{[B]_{\meager(Q)}:B\in\base\}$$
generates $\bool$.  Let
$$R=\{X\intersect B: B\in\base\}.$$
If $\tilde{X}\subseteq \cantorspace$ is the image of $X$ under
the characteristic
function of the sequence $\base$ (see Theorem \ref{charseq}),
then $\tilde{X}$ has hereditary order $\omega_1$.  Of course
$\tilde{X}$ is really just the same as $X$ but retopologized using
$\base$ as a family of basic open sets.  Let $Y\in[X]^{\omega_1}$.
Since $\ord(p)=\omega_1$ for any basic clopen set the following
claim shows that
the order of $Y$ (or rather the image of $Y$ under the
characteristic function of the sequence $\base$) is $\omega_1$.

\claim There exists a basic clopen $p$ in $Q$ such that for every
Borel $B\subseteq p$,
$$|B\intersect Y|\leq\omega \rmiff B \mbox{ meager.} $$

\proof
Let $p$ and $q$ stand for nonempty basic clopen sets. Obviously
if $B$ is meager then $B\intersect Y$ is countable, since
$B\intersect X$ is countable.  To prove the other direction,
suppose for contradiction that for every $p$ there exists
$q\leq p$ and Borel $B_q\subseteq q$ such that $B_q$ is comeager
in $q$ and $B_q\intersect Y$ is countable.  By using ccc there exists
a countable dense family $\Sigma$ and $B_q$ for $q\in \Sigma$
with
$B_q\subseteq q$ Borel and comeager
in $q$ such that  $B_q\intersect Y$ is countable.  But
$$B=\Union \{B_q:q\in \Sigma\}$$
is a comeager Borel set which meets $Y$ in a countable set.  This
implies that $Y$ is countable since $X$ is contained in $B$ except
for countable many points.
\qed

\begin{theorem} \label{cohenmod2}
Suppose $G$ is $\fin(\kappa,2)$-generic over $V$
where $V$ is a model of CH and $\kappa\geq\omega$.
Then
in $V[G]$
there exists a separable metric space $X$ with
$|X|=\omega_1$ and hereditarily of order $\omega_1$.
\end{theorem}
\proof
Immediate from Theorem \ref{cohen1} and \ref{cohen2}.
\qed

Finally, we show that there are no large spaces of hereditary order
$\omega_1$ in the Cohen real model.

\begin{theorem} \label{cohenmodel3}
Suppose $G$ is $\fin(\kappa,2)$-generic over $V$
where $V$ is a model of CH and $\kappa\geq\omega_2$.
Then
in $V[G]$
for every separable metric space $X$ with
$|X|=\omega_2$ there exists $Y\in [X]^{\omega_2}$ with
 $\ord(Y)<\omega_1$.
\end{theorem}
\proof
By the argument used in the proof of Lemma \ref{mapping}
we can find
$$\langle G_\alpha:\alpha<\omega_2\rangle\in V[G]$$
 which is
$\sum_{\alpha<\omega_2}\fin(\omega,2)$-generic over
$V$ and
a $\fin(\omega,2)$-name $\tau$ for an element of $\cantorspace$ such
that $Y=\{\tau^{G_\alpha}:\alpha<\omega_2\}$ is subset of
$X$.  We claim that $\ord(Y)<\omega_1$.
Let $${\famly}=\{\val{\tau\in C}: C\subseteq\cantorset \mbox{ clopen }\}$$
where boolean values are in the
the unique complete boolean algebra $\bool$ in which $\fin(\omega,2)$
is dense.
Let $\fbool$ be the complete subalgebra of $\bool$ which is generated
by $\famly$.  First note that the order of $\famly$ in $\fbool$ is
less than $\omega_1$.  This is because $\fbool$ contains a countable dense
set:
$$D=\{\prod\{c\in {\fbool}:p\leq c\}: p\in \fin(\omega,2)\}.$$
Since $D$ is countable and $\bsig01(D)={\fbool}$, it follows that
the order of $\famly$ is countable.

I claim that the order of $Y$ is essentially less than or equal
to the order of $\famly$ in ${\fbool}$.

\begin{lemma}\label{boolborel}
 Let $\bool$ be a cBa, $\tau$ a $\bool$-name for an element of
 $\cantorset$, and
 $\famly=\{\val{\tau\in C}: C\subseteq\cantorset \mbox{ clopen }\}$.
 Then for each $B\subseteq\cantorspace$ a $\bpi0\alpha$ set coded in $V$
  the set
  $
  \val{\tau\in\check{B}}
  $
  is $\bpi0\alpha({\famly})$ and conversely, for every
  $c\in \bpi0\alpha({\famly})$ there exists a
  $B\subseteq\cantorspace$ a $\bpi0\alpha$ set coded in $V$ such
  that
  $
  c = \val{\tau\in\check{B}}
  .$
\end{lemma}
\proof
Both directions of the lemma are simple inductions.
\qed

Now suppose the order of $\famly$ in ${\fbool}$ is $\alpha$.
Let $B\subseteq\cantorspace$ be any  Borel set coded in $V[G]$.
By ccc there exists $H=G\res\Sigma$ where $\Sigma\subseteq\kappa$
is countable set in $V$ such that $B$ is coded in $V[H]$.  Consequently,
since we could replace $V$ with $V[H]$ and delete countably many elements
of $Y$ we may as well assume that $B$ is coded in the ground model.
Since the order of $\famly$ is $\alpha$ we have by the lemma
that there exists a $\bpi0\alpha$ set $A$ such that
$$\val{\tau\in \check{A}} = \val{\tau\in \check{B}}.$$
It follows that
$$Y\intersect A= Y\intersect B$$
and hence order of $Y$ is less than or equal to $\alpha$
(or three since we neglected countably many elements of $Y$).

\qed

\sector{The random real model}

In this section we consider the question of \label{random}
Borel orders in the random real model. We conclude with a few
remarks about perfect set forcing.

A set $X\subseteq\cantorspace$ is a \dex{Sierpinski set} iff
$X$ is uncountable and for every measure zero set
$Z$ we have $X\intersect Z$ countable.  Note that
by Mahlo's Theorem \ref{mah} we know that under CH
Sierpinski sets exists. Also it is easy to see that in the
random real model, the set of reals given by the generic
filter is a Sierpinski set.

\begin{theorem}
  (Poprougenko \site{pop}) \label{pop}
  If $X$ is Sierpinski, then $\ord(X)=2$.
\end{theorem}
\proof
For any Borel set $B\subseteq\cantorspace$ there exists
an $F_\sigma$ set with $F\subseteq B$ and $B\setminus F$
measure zero.
Since $X$ is Sierpinski $(B\setminus F)\intersect X=F_0$ is
countable, hence $F_\sigma$.  So
$$B\intersect X= (F\union F_0)\intersect X.$$
\qed

I had been rather hoping that every uncountable separable metric
space in the random real model has order either 2 or $\omega_1$.
The following result shows that  this is definitely not the case.

\begin{theorem} \label{rand}
  Suppose $X\in V$ is a subspace of $\cantorset$ of order $\alpha$
  and $G$ is measure algebra $2^\kappa$-generic over $V$, i.e.
  adjoin $\kappa$ many random reals.

  Then  $V[G]\models \alpha\leq\ord(X)\leq\alpha+1$.
\end{theorem}

The result will easily follow from the next two lemmas.

Presumably, $\ord(X)=\alpha$ in $V[G]$,
but I haven't been able to prove this.  Fremlin's proof
(Theorem \ref{frem}) having filled up one such missing gap,
 leaving this gap here restores a certain
cosmic balance of ignorance.\footnote{All things I thought I knew;
but now confess,
the more I know I know, I know the less.- John Owen (1560-1622)}

Clearly, by the usual ccc arguments, we may assume that
that $\kappa=\omega$ and $G$ is just a random real.
In the following lemmas boolean values $\val{\theta}$ will
be computed in the measure algebra $\bool$ on $\cantorspace$.
Let $\mu$ \sdex{$\mu$} be the usual product measure on $\cantorspace$.

\begin{lemma} \label{random1}
Suppose $\epsilon$ a real,  $b\in\bool$, and
$\name{U}$ the name of a $\bpi0\alpha$
subset of $\cantorspace$ in $V[G]$.   Then the set
$$\{x\in\cantorspace: \mu(b\conj\val{\check{x}\in\name{U}})\geq\epsilon\}$$
is $\bpi0\alpha$ in $V$.
\end{lemma}
\proof
The proof is by induction on $\alpha$.

\bigskip
 Case $\alpha=1$.
\par\noindent
If $\name{U}$ is a name for a closed set, then
$$\val{\check{x}\in \name{U}}=\prod_{n\in\omega}\val{[x\res n]\intersect
\name{U}\not=\emptyset}.$$
Consequently,
$$\mu(b\conj\val{\check{x}\in\name{U}})\geq\epsilon$$
iff
$$\forall n\in\omega\;\;\mu(b\conj\val{[x\res n]\intersect
\name{U}\not=\emptyset})\geq\epsilon$$
and the set is closed.

\bigskip
Case $\alpha>1$.
\par\noindent  Suppose $\name{U}=\Intersect_{n\in\omega}\comp{\name{U_n}}$
where
each $\name{U_n}$ is a name for a $\bpi0{\alpha_n}$ set
for some $\alpha_n<\alpha$.  We can assume that the sequence
$\comp{U_n}$ is descending.
Consequently,

$$\mu(b\conj\val{\check{x}\in\check{U}})\geq\epsilon$$
iff
$$\mu(b\conj\val{\check{x}\in\Intersect_{n\in\omega}\comp{\name{U_n}} })
  \geq\epsilon $$
iff
$$ \forall n \in \omega \;\;
\mu(b\conj\val{\check{x}\in\comp{\name{U_n}} })
  \geq\epsilon $$
iff
$$ \forall n \in \omega \;\;
\mbox{ not }\mu(b\conj\val{\check{x}\in{\name{U_n}} })
  >\mu(b)-\epsilon. $$
iff
$$ \forall n \in \omega \;\;
\mbox{ not }\exists m\in \omega\;\;\mu(b\conj\val{\check{x}\in{\name{U_n}} })
  \geq\mu(b)-\epsilon+1/m $$
By induction, each of the sets
$$\{x\in\cantorspace: \mu(b\conj\val{\check{x}\in{\name{U_n}} })
  \geq\mu(b)-\epsilon+1/m\}$$
 is $\bpi0{\alpha_n}$ and so the result is proved.
\qed

It follows from this lemma that if $X\subseteq\cantorspace$ and
$V\models$``$\ord(X)>\alpha$'', then
$V[G]\models$``$\ord(X)>\alpha$''.  For suppose $F\subseteq\cantorspace$
is $\bsig0\alpha$ such that for every $H\subseteq\cantorspace$ which
is $\bpi0\alpha$ we have $F\intersect X\not= H\intersect X$.
Suppose for contradiction that

$$b\forces \mbox{``}\name{U}\intersect \check{X}=\check{F}\intersect \check{X} \
rmand
\name{U}\mbox{ is } \bpi0\alpha\mbox{''}.$$
But then
$$\{x\in\cantorspace: \mu(b\conj\val{\check{x}\in\check{U}})=\mu(b)\}$$
is a $\bpi0\alpha$ set which must be equal to $F$ on $X$, which is
a contradiction.

To prove the other direction of the inequality we will use the following
lemma.
\begin{lemma} \label{random2}
  Let $G$ be $\bool$-generic
  (where $\bool$ is the measure algebra on $\cantorspace$)
  and $r\in\cantorspace$ is the associated
  random real.  Then for any $b\in\bool$
  $$b\in G \rmiff
   \forall^\infty n\; \mu([r\res n]\conj b)\geq {3\over 4}\mu([r\res n]).$$
\end{lemma}
\proof
Since $G$ is an ultrafilter it is enough to show that
$b\in G$ implies
$$\forall^\infty n\; \mu([r\res n]\conj b)\geq {3\over 4}\mu([r\res n]).$$
Let $\bool^+$ be the nonzero elements of $\bool$.
To prove this it suffices to show:
\claim For any $b\in\bool^+$ and
for every $d\leq b$ in $\bool^+$ there exists a tree
$T\subseteq 2^{<\omega}$ with $[T]$ of positive measure, $[T]\leq d$,
and
$$\mu([s]\intersect b)\geq {3\over 4}\mu([s])$$
for all but finitely many $s\in T$.
\proof
Without loss we may assume that $d$ is a closed set and
let $T_d$ be a tree such that $d=[T_d]$.  Let $t_0\in T_d$ be such
that $$\mu([t_0]\intersect [T_d])\geq {9\over 10}\mu([t_0]).$$
Define a subtree $T\subseteq T_d$ by $r\in T$ iff $r\subseteq t_0$ or
$t_0\subseteq r$ and
$$\forall t\;( t_0\subseteq t\subseteq r\mbox{ implies } \mu([t]\intersect b)
\geq{3\over 4}\mu([t])\;).$$
So we only need to see that $[T]$ has positive measure.  So suppose
for contradiction that $\mu([T])=0$.  Then for some sufficiently large
$N$
$$\mu(\Union_{s\in T\intersect 2^N} [s])\leq {1\over 10}\mu([t_0]).$$
For every $s\in T_d\intersect 2^N$ with $t_0\subseteq s$, if
$s\notin T$ then there exists $t$ with $t_0\subseteq t\subseteq s$ and
$\mu([t]\intersect b)<{3\over 4}\mu([t])$.  Let $\Sigma$ be a maximal
antichain of $t$ like this.  But note that
$$[t_0]\intersect [T_d]\subseteq \Union_{s\in 2^N\intersect T}[s]\union
\Union_{t\in\Sigma}([t]\intersect b).$$
By choice of $\Sigma$
$$\mu(\Union_{s\in\Sigma}[s]\intersect b)\leq {3\over 4}\mu([t_0])$$
and by choice of $N$
$$\mu(\Union_{s\in 2^N\intersect T}[s])\leq {1\over 10}\mu([t_0])$$
which contradicts the choice of $t_0$:
$$\mu([t_0]\intersect [T_d])\leq ({1\over 10}+{3\over 4})
\mu([t_0])=\footnote{Trust me on this, I have been teaching a lot
of Math 99 ``College Fractions''.}
\;\;{17\over 20}\mu([t_0])<{9\over 10} \mu([t_0]).$$
This proves the claim and the lemma.
\qed

So now suppose that the order of $X$ in $V$ is  $\leq \alpha$.
We show that it is $\leq \alpha+1$ in $V[G]$.  Let $\name{U}$ be
any name for a Borel subset of $X$ in the extension.  Then
we know that $x\in U^G$ iff $\val{\check{x}\in \name{U}}\in G$.
By Lemma \ref{random1} we know that
for any  $s\in 2^{<\omega}$ the set
$$B_s=\{x\in X:
    \mu([s]\intersect \val{\check{x}\in \name{U}})\geq {3\over 4}\mu([s])\}$$
is a Borel subset of $X$ in the ground model and hence is
$\bpi0\alpha(X)$.  By Lemma \ref{random2} we have
that for any $x\in X$
$$x\in U \rmiff \forall^\infty n \; \; x\in B_{r\res n}$$
and so $U$ is $\bsig0{\alpha+1}(X)$ in V[G].

This concludes the proof of Theorem \ref{rand}.
\qed

Note that this result does allow us to get sets of order
$\lambda$ for any countable limit ordinal $\lambda$ by taking a
clopen separated union of a sequence of sets whose order goes up $\lambda$.

Also a Luzin set $X$ from the ground model has order $3$ in the random
real extension.  Since $(\ord(X)=3)^V$ we know that
$(3\leq \ord(X)\leq 4)^{V[G]}$.  To see that $(\ord(X)\leq 3)^{V[G]}$ suppose
that $B\subseteq X$ is Borel in $V[G]$. The above proof shows that there
there exists Borel sets $B_n$ each coded in $V$ (but the sequence
may not be in $V$) such that
$$x\in B \rmiff \forall^{\infty} n \;\; x\in B_n.$$
For each $B_n$ there exists an open set $U_n\subseteq X$
such that $B_n\Delta U_n$ is countable. If we let
$$C=\Union_{n\in\omega}\Intersect_{m>n} U_m$$
then $C$ is $\bsig03(X)$ and $B\Delta C$ is countable.  Since
subtracting and adding a countable set from a $\bsig03(X)$ is
still $\bsig03(X)$ we have that $B$ is $\bsig03(X)$ and so
the order of $X$ is $\leq 3$ in $V[G]$.

\begin{theorem} \label{randlast}
  Suppose $V$ models CH and $G$ is measure algebra on $2^\kappa$-generic
  over $V$ for some $\kappa\geq \omega_2$.  Then in $V[G]$ for
  every $X\subseteq \cantorspace$ of cardinality $\omega_2$
  there exists $Y\in [X]^{\omega_2}$ with $\ord(Y)=2$.
\end{theorem}
\proof
Using the same argument as in the proof of Theorem \ref{cohenmodel3}
we can get a Sierpinski set $S\subseteq\cantorspace$ of cardinality
$\omega_2$ and a term $\tau$ for any element of $\cantorspace$ such
that $Y=\{\tau^r:r\in S\}$ is a set of distinct elements of
$X$.  This Sierpinski set has two additional properties:
every element of it is random over the ground model and
it meets every set of positive measure.

We will show that the order of $Y$ is $2$.

\begin{lemma}
  Let $\famly\subseteq\bool$ be any subset
  of a measure algebra $\bool$ closed under finite conjunctions.
  Then $\bpi02(\famly)=\bsig02(\famly)$, i.e. $\famly$ has order $\leq 2$.
\end{lemma}
\proof
Let $\mu$ be the measure on $\bool$.

\medskip
\noindent (1) For any $b\in\bpi01(\famly)$ and real $\epsilon>0$
there exists $a\in\famly$ with $b\leq a$ and $\mu(a-b)<\epsilon$.

pf:\footnote{Pronounced `puff'.} $\;\; b=\prod_{n\in\omega} a_n$.
Let $a=\prod_{n<N}a_n$ for
some sufficiently large $N$.

\medskip
\noindent (2) For any $b\in\bsig02(\famly)$ and real $\epsilon>0$
there exists $a\in\bsig01(\famly)$ with $b\leq a$ and $\mu(a-b)<\epsilon$.

pf: $b=\sum_{n<\omega}b_n$.  Applying (1) we get $a_n\in\famly$ with
$b_n\leq a_n$ and $$\mu(a_n-b_n)<{\epsilon\over{2^n}}.$$
Then let $a=\sum_{n\in\omega} a_n$.

\bigskip
Now suppose $b\in\bsig02(\famly)$. Then by applying (2) there exists
$a_n\in\bsig01(\famly)$ with $b\leq a_n$ and $\mu(a_n-b)<1/n$.
Consequently, if $a=\prod_{n\in\omega}a_n$, then $b\leq a$ and
$\mu(a-b)=0$ and so $a=b$.
\qed
Let $${\famly}=\{\val{\tau\in C}: C\subseteq\cantorset \mbox{ clopen }\}$$
where boolean values are in the measure algebra $\bool$ on
$\cantorspace$.
Let $\fbool$ be the complete subalgebra of $\fbool$ which is generated
by $\famly$.

Since the order of $\famly$ is $2$, by Lemma \ref{boolborel}
we have
that for any Borel set $B\subseteq Y$ there exists a $\bsig02(Y)$ set
$F$ such that $y\in B$ iff $y\in F$ for all but countably many $y\in Y$.
Thus we see that the order of $Y$ is $\leq 3$.  To get it down
to $2$ we use the following lemma.  If $B=(F\setminus F_0)\union F_1$ where
$F_0$ and $F_1$ are countable and $F$ is $\bsig02$, then by the lemma
$F_0$ would be $\bpi02$ and thus $B$ would be $\bsig02$.

\begin{lemma}
 Every countable subset of $Y$ is $\bpi02(Y)$.
\end{lemma}
\proof

It suffices to show that every countable subset of $Y$ can be covered
by a countable $\bpi02(Y)$ since one can always subtract a countable
set from a $\bpi02(Y)$ and remain $\bpi02(Y)$.

For any $s\in 2^{<\omega}$ define
$$b_s=\val{s\subseteq \tau}.$$
Working in the ground model let $B_s$ be a Borel
set with $[B_s]_{\bool}=b_s$.  Since the
Sierpinski set consists only of reals random over the ground
model we know that for every $r\in S$
$$r\in B_s\rmiff s\subseteq \tau^r.$$
Also since the Sierpinski set meets every Borel set of positive measure
we know that for any $z\in Y$ the set $\Intersect_{n<\omega}B_{z\res n}$
has measure zero.   Now let $Z=\{z_n:n<\omega\}\subseteq Y$ be arbitrary
but listed with infinitely many repetitions.  For each $n$ choose $m$ so
that
if $s_n=z_n\res m$, then $\mu(B_{s_n})<1/{2^n}$.
Now for every $r\in S$  we have that
$$r\in\Intersect_{n<\omega}\Union_{m>n} B_{s_m} \rmiff
  \tau^r\in\Intersect_{n<\omega}\Union_{m>n} [s_m]. $$
The set $H=\Intersect_{n<\omega}\Union_{m>n} [s_m]$ covers $Z$ and
is $\bpi02$.
It has countable intersection with $Y$ because
the set $\Intersect_{n<\omega}\Union_{m>n} B_{s_m}$ has
measure zero.

This proves the Lemma and Theorem \ref{randlast}.
\qed

%%%%%%%%%%%%%%% sacks model %%%%%%%%%%%%%%%

\newpage

\begin{center}
  Perfect Set Forcing
\end{center}

 \sdex{perfect set forcing}

In the  iterated Sack's real model the continuum is
$\omega_2$ and every set $X\subseteq\cantorspace$
of cardinality $\omega_2$ can
be mapped continuously onto $\cantorspace$ (Miller \site{maponto}).
It follows from Rec{\l}aw's
Theorem \ref{rec}
that in this model every separable metric space of cardinality
$\omega_2$ has order $\omega_1$.  On the other hand this
forcing (and any other with the Sack's property) has the property
that every meager set in the extension is covered by a meager set
in the ground model and every measure set in the extension is covered
by a measure zero set in the ground model (see Miller \site{someprop}).
Consequently, in this model there are Sierpinski sets and Luzin sets
of cardinality $\omega_1$.  Therefore in the iterated Sacks real
model there are separable metric spaces of cardinality $\omega_1$
of every order $\alpha$ with $2\leq\alpha<\omega_1$.  I do not know
if there is an uncountable separable metric space which is hereditarily
of order $\omega_1$ in this model.

Another way to obtain the same orders is to use the construction
of Theorem 22 of Miller \site{borelhier}.  What was done there implies
the following:
\begin{quote}
For any model $V$ there
exists a ccc extension $V[G]$ in which every uncountable separable
metric space has order $\omega_1$.
\end{quote}

If we apply this result $\omega_1$ times with a finite support extension,

we get a model, $V[G_\alpha:\alpha<\omega_1]$,  where there are
separable metric spaces of all orders of cardinality $\omega_1$,
but every separable metric space of cardinality $\omega_2$ has
order $\omega_1$.

To see the first fact note that  $\omega_1$ length finite
support iteration always adds
a Luzin set.  Consequently, by Theorem \ref{cohenluz},
for each $\alpha$ with
$2<\alpha<\omega_1$ there
exists a separable metric space of cardinality $\omega_1$
which is hereditarily of order $\alpha$.   Also there is such an
$X$ of order $2$ by the argument used in Theorem \ref{cohenmodel}.

On the other hand if $X$ has cardinality $\omega_2$ in
$V[G_\alpha:\alpha<\omega_1]$, then for some $\beta<\omega_1$ there
exists and uncountable $Y\in V[G_\alpha:\alpha<\beta]$
with $Y\subseteq X$. Hence $Y$ will have order $\omega_1$ in
$V[G_\alpha:\alpha<\beta+1]$ and
by examining the proof it is easily seen that it remains of order
$\omega_1$ in $V[G_\alpha:\alpha<\omega_1]$.

\sector{Covering number of an ideal}

This section is a small diversion.\footnote{All men's gains are the fruit
of venturing.  Herodotus BC 484-425.}
It is motivated
Theorem \ref{marsolthm} of  Martin and Solovay.

Define for any ideal $I$ in $\borel(\cantorspace)$
$$\cov(I)=\min\{|\ideal|: \ideal\subseteq I, \; \Union\ideal=\cantorset\}.$$
The following theorem is well-known. \sdex{$\cov(I)$}
\begin{theorem} \label{count}
  For any cardinal $\kappa$ the following are equivalent:
  \begin{enumerate}
    \item \dex{MA$_\kappa$(ctbl)}, i.e. for any countable poset, $\poset$,
    and family $\dense$ of dense subsets of $\poset$ with $|\dense|<\kappa$
    there exists a $\poset$-filter $G$ with $G\intersect D\not=\emptyset$
    for every $D\in\dense$, and
    \item $\cov(\meager(\cantorset))\geq\kappa$.
  \end{enumerate}
\end{theorem} \sdex{$\cov(\meager(\cantorset))$}
\proof
MA$_\kappa$(ctbl) implies $\cov(\meager(\cantorset))\geq\kappa$, is easy
because if $U\subseteq\cantorspace$ is a dense open set, then
$$D=\{s\in 2^{<\omega}: [s]\subseteq U\}$$
is dense in $2^{<\omega}$.

$\cov(\meager(\cantorset))\geq\kappa$ implies MA$_\kappa$(ctbl) follows
from the fact that any countable poset, $\poset$, either contains
a dense copy of $2^{<\omega}$ or contains a $p$ such that
every two extensions of $p$ are compatible.
\qed

\begin{theorem} \label{cofcov}
  (Miller \site{ctlbcof}) $\cof(\cov(\meager(\cantorset)))>\omega$, e.g.,
  it is impossible to have $\cov(\meager(\cantorset))=\aleph_\omega$.
\end{theorem}
\proof
Suppose for contradiction that $\kappa=\cov(\meager(\cantorset))$ has
countable cofinality and let $\kappa_n$ for $n\in\omega$ be a cofinal
sequence in $\kappa$. Let $\langle C_\alpha:\alpha<\kappa\rangle$ be
a family of closed nowhere dense sets which cover $\cantorset$.
We will construct a sequence $P_n\subseteq\cantorset$
of perfect sets with the properties that
\begin{enumerate}
  \item $P_{n+1}\subseteq P_n$,
  \item $P_n\intersect \Union\{C_\alpha:\alpha<\kappa_n\}=\emptyset$, and
  \item $\forall \alpha<\kappa\;\;\;C_\alpha\intersect P_n$ is nowhere dense
  in $C_\alpha$.
\end{enumerate}
This easily gives a contradiction, since $\Intersect_{n<\omega}P_n$ is
nonempty and disjoint from all $C_\alpha$, contradicting the
fact that the $C_\alpha$'s cover $\cantorspace$.

We show how to obtain $P_0$, since the argument easily relativizes to show
how to obtain $P_{n+1}$ given $P_n$. Since
$\cov(\meager(\cantorset))>\kappa_n$ there exists a countable sequence
$$D=\{x_n:n\in\omega\}\subseteq\cantorspace$$
such that $D$ is dense and for every $n$
$$x_n\notin\Union_{\alpha<\kappa_n}C_\alpha.$$

Consider the following forcing notion $\poset$.
$$\poset=\{(H,n): n\in\omega \rmand H\in [D]^{<\omega}\}$$
This is ordered by $(H,n)\leq (K,m)$ iff
\begin{enumerate}
\item $H\supseteq K$,
\item $n\geq m$, and
\item for every $x\in H$ there exists $y\in K$ with $x\res m= y\res m$.
\end{enumerate}
Note that $\poset$ is countable.

For each $n\in\omega$ define $E_n\subseteq\poset$ by $(H,m)\in E_n$
iff
\begin{enumerate}
\item $m>n$ and
\item $\forall x\in H \exists y\in H\;
x\res n=y\res n$ but $x\res m\not=y\res m$.
\end{enumerate}
and for each $\alpha<\kappa_0$ let
$$F_\alpha=\{(H,m)\in\poset: \forall x\in H\;\;
[x\res m]\intersect C_\alpha=\emptyset\}.$$
For $G$ a $\poset$-filter, define $X\subseteq D$ by
$$X=\Union\{H: \exists n\; (H,n)\in G\}$$
and let $P=cl(X)$.
It easy to check that the $E_n$'s are dense and if $G$ meets each
one of them, then $P$ is perfect (i.e. has no isolated points).
The $F_\alpha$ for $\alpha<\kappa_0$ are
dense in $\poset$. This is because $D\intersect C_\alpha=\emptyset$
so given $(H,n)\in\poset$ there exists $m\geq n$ such that
for every $x\in H$ we have $[x\res m]\intersect C_\alpha=\emptyset$
and thus $(H,m)\in F_\alpha$.  Note that if
$G\intersect F_\alpha\not=\emptyset$, then
$P\intersect C_\alpha=\emptyset$. Consequently, by Theorem \ref{count},
there exists a $\poset$-filter $G$ such that $G$ meets
each $E_n$ and all $F_\alpha$ for $\alpha<\kappa_0$.
Hence $P=cl(X)$ is a perfect set which is disjoint from
each $C_\alpha$ for $\alpha<\kappa_0$. Note also that
for every $\alpha<\kappa$ we have that $C_\alpha\intersect D$ is
finite and hence $C_\alpha\intersect X$ is finite and
therefore $C_\alpha\intersect P$ is nowhere dense in $P$.
This ends the construction of $P=P_0$ and since the $P_n$ can
be obtained with a similar argument, this proves the Theorem.
\qed

\begin{question}
  (Fremlin) Is the same true for the measure zero ideal in place of
  the ideal of meager sets?
\end{question}

Some partial results are known (see Bartoszynski, Judah, Shelah
\site{barto}\site{barto2}\site{barto3}).

\begin{theorem}
  (Miller \site{ctlbcof}) It is consistent that
  $\cov(\meager(2^{\omega_1}))=\aleph_\omega$.
\end{theorem}
\proof
In fact, this holds in the model obtained by forcing with
$\fin(\aleph_\omega,2)$ over a model of GCH.

\bigskip

$\cov(\meager(2^{\omega_1}))\geq\aleph_\omega$: Suppose
for contradiction that
$$\{C_\alpha:\alpha<\omega_n\}\in V[G]$$
is a family
of closed nowhere dense sets covering $2^{\omega_1}$.
Define
$$E_\alpha=\{s\in \fin(\omega_1,2): [s]\intersect C_\alpha=\emptyset\}.$$
Using ccc, there exists $\Sigma\in[\aleph_\omega]^{\omega_n}$
in $V$ with
$$\{E_\alpha:\alpha<\omega_n\}\in V[G\res\Sigma].$$
Let $X\subseteq \aleph_\omega$ be a set in $V$ of cardinality $\omega_1$
which is disjoint from $\Sigma$.  By the product lemma
$G\res X$ is $\fin(X,2)$-generic over $V[G\res \Sigma]$.  Consequently,
if $H:\omega_1\arrow 2$ corresponds to $G$ via an isomorphism
of $X$ and $\omega_1$, then $H\notin C_\alpha$ for every
$\alpha<\omega_n$.

\bigskip

$\cov(\meager(2^{\omega_1}))\leq\aleph_\omega$:  Note that for
every uncountable $X\subseteq\omega_1$ with $X\in V[G]$ there exists
$n\in\omega$ a $Z\in[\omega_1]^{\omega_1}\intersect V[G\res\omega_n]$ with
$Z\subseteq X$. To see this note that for every $\alpha\in X$ there
exists $p\in G$ such that
$p\forces \alpha\in X$
and $p\in\fin(\omega_n,2)$ for some $n\in\omega$.  Consequently, by ccc,
some $n$ works for uncountably many $\alpha$.

Consider the family
of all closed nowhere dense sets $C\subseteq 2^{\omega_1}$ which
are coded in some $V[G\res\omega_n]$ for some $n$.  We claim that
these cover $2^{\omega_1}$.  This follows from above,
because for any $Z\subseteq\omega_1$ which is infinite the
set
$$C=\{x\in 2^{\omega_1}: \forall \alpha\in Z \;\; x(\alpha)=1\}$$
is nowhere dense.

\qed

\begin{theorem}
  (Miller \site{ctlbcof}) It is consistent that
  there exists a ccc $\sigma$-ideal $I$ in
  $\borel(\cantorset)$ such that $\cov(I)=\aleph_\omega$.
\end{theorem}

\proof
Let $\poset=\fin(\omega_1,2)*\name{\posetq}$ where $\name{\posetq}$
is a name for the Silver forcing which codes up generic filter
for $\fin(\omega_1,2)$
just like in the proof of Theorem \ref{marsolthm}.  Let
$\prod_{\alpha<\aleph_{\omega}}\poset$ be the direct sum (i.e. finite
support product) of $\aleph_{\omega}$ copies of $\poset$.  Forcing
with the direct sum adds a
filter $G=\langle G_\alpha:\alpha<\aleph_{\omega} \rangle$
where each $G_\alpha$ is $\poset$-generic.  In general,
a direct sum is ccc iff every finite subproduct is ccc.  This follows by
a delta-system argument.  Every finite product of $\poset$ has
ccc, because $\poset$ is $\sigma$-centered, i.e., it is the countable
union of centered sets.

Let $V$ be a model of GCH and
$G=\langle G_\alpha:\alpha<\aleph_{\omega} \rangle$
be $\prod_{\alpha<\aleph_{\omega}}\poset$ generic over $V$.
We claim that in $V[G]$ if $I$ is the $\sigma$-ideal given by Sikorski's
Theorem \ref{sikorthm} such that $\prod_{\alpha<\aleph_{\omega}}\poset$ is
densely embedded into $\borel(\cantorset)/I$ then
$\cov(I)=\aleph_\omega$.

First define,
$m_\poset$, to be the cardinality of the minimal failure
of MA for $\poset$, i.e., the least $\kappa$ such that
there exists a family $|\dense|=\kappa$ of dense subsets of
$\poset$ such that there is no $\poset$-filter meeting all
the $D\in\dense$. \sdex{$m_?poset$}

\begin{lemma}
  In $V[\langle G_\alpha:\alpha<\aleph_{\omega}\rangle]$ we have that
  $m_\poset=\aleph_\omega$.
\end{lemma}
\proof
Note that for any set $D\subset\poset$ there exists a set
$\Sigma\in[\aleph_\omega]^{\omega_1}$ in $V$ with
$D\in V[\langle G_\alpha:\alpha\in\Sigma\rangle]$.  So if
$|\dense|=\omega_n$ then there exists
$\Sigma\in[\aleph_\omega]^{\omega_n}$ in $V$ with
$\dense\in V[\langle G_\alpha:\alpha\in\Sigma\rangle]$.
Letting $\alpha\in\aleph_\omega\setminus \Sigma$ we get
$G_\alpha$ a $\poset$-filter meeting every $D\in\dense$.
Hence   $m_\poset\geq\aleph_\omega$.

On the other hand:
\claim For every
$X\in[\omega_1]^\omega_1\intersect
V[\langle G_\alpha:\alpha<\aleph_{\omega}\rangle]$ there exists $n\in\omega$
and $Y\in [\omega_1]^\omega_1\intersect
V[\langle G_\alpha:\alpha<\aleph_{n}\rangle]$ with $Y\subseteq X$.
\proof
For every $\alpha\in X$ there exist $p\in G$ and $n<\omega$ such that
$p\forces \check{\alpha}\in \name{X}$ and $\dom(p)\subseteq\aleph_n$.
Since $X$ is uncountable there is one $n$ which works for uncountably
many $\alpha\in X$.
\qed
It follows from the Claim that there is no $H$ which is
$\fin(\omega_1,2)$ generic over all the models
$V[\langle G_\alpha:\alpha<\aleph_{n}\rangle]$, but forcing with $\poset$
would add such an $H$ and so   $m_\poset\leq\aleph_\omega$ and
the Lemma is proved. \qed

\begin{lemma}
  If $\poset$ is ccc and dense in the cBa $\borel(\cantorset)/I$,
   then $m_\poset=\cov(I)$.
\end{lemma}
\proof
This is the same as Lemma \ref{martinlem} equivalence of (1) and (3),
except you have to check that $m$ is the same for both
$\poset$ and $\borel(\cantorset)/I$.
\qed

Kunen \site{kunenma} showed that least cardinal for which MA fails
can be a singular cardinal of cofinality $\omega_1$, although
it is impossible for it to have cofinality $\omega$ (see Fremlin
\site{frembook}).
It is still
open whether it can be a singular cardinal of cofinality
greater than $\omega_1$ (see Landver \site{landver}).
Landver \site{landver2} generalizes Theorem \ref{cofcov}
to the space $2^\kappa$ with basic clopen sets of
the form $[s]$ for $s\in 2^{<\kappa}$.  He uses a generalization
of a characterization
of $\cov(\meager(\cantorspace))$ due to Bartoszynski \site{barto4}
and Miller \site{char}.

\addtocontents{toc}{\protect\newpage}
\newpage
\part{Analytic sets}
\section{Analytic sets}

Analytic sets were discovered by Souslin when he encountered a mistake
of Lebesgue.  Lebesgue had erroneously proved that the Borel sets
were closed under projection.  I think the mistake  he made was
to think that the countable intersection commuted with projection.
A good reference is the volume devoted to analytic sets edited by
Rogers \site{anal}.   For the more classical viewpoint
of operation-A, see Kuratowski \site{kur}.  For the whole
area of descriptive set theory and its history, see
Moschovakis \site{mosch}.

\bigskip

Definition. A set $A\subseteq\bairespace$ is $\Sigma^1_1$ iff
\sdex{$\Sigma^1_1$}
there exists a recursive
$$R\subseteq \Union_{n\in\omega}(\omega^n\cross\omega^n)$$
such that for all $x\in\bairespace$
$$x\in A \rmiff \exists y \in\bairespace\;\forall n\in\omega\;\;
R(x\res n, y\res n).$$
A similar definition applies for $A\subseteq \omega$ and also
for $A\subseteq \omega\cross\bairespace$ and so forth.  For
example, $A\subseteq\omega$ is $\Sigma^1_1$ iff there
exists a recursive $R\subseteq \omega\cross \finseq$ such that
for all $m\in\omega$
$$m\in A \rmiff \exists y \in\bairespace\;\forall n \in\omega\;\;
R(m, y\res n).$$
A set $C\subseteq\bairespace\cross\bairespace$ is $\Pi_1^0$ iff
\sdex{$\Pi_1^0$}
there exists a recursive predicate
$$R\subseteq \Union_{n\in\omega}(\omega^n\cross\omega^n)$$ such that
$$C=\{(x,y):\forall n \; R(x\res n, y\res n)\}.$$  That means basically
that $C$ is a recursive closed set.

The $\Pi$ classes are the complements of the $\Sigma$'s and
 $\Delta$ is the class of sets which are both $\Pi$ and $\Sigma$.
The relativized classes, e.g. $\Sigma^1_1(x)$ are obtained by
\sdex{$\Sigma^1_1(x)$}
allowing $R$ to be recursive in $x$, i.e., $R\leq_T x$.  The boldface
classes, e.g., $\bsig11$, $\bpi11$, are obtained by taking arbitrary $R$'s.
\sdex{$?bsig11$} \sdex{$?bpi11$}

\begin{lemma}
$A\subseteq\bairespace$ is $\Sigma^1_1$ iff there exists
set $C\subseteq\bairespace\cross\bairespace$ which is $\Pi_1^0$
and $$A=\{x\in\bairespace: \exists y \in \bairespace\;\; (x,y)\in C\}.$$
\end{lemma}

\begin{lemma}   \label{catchall} The following are all true:
   \begin{enumerate}
        \item For every $s\in\finseq$ the basic
     clopen set $[s]=\{x\in\bairespace: s\subseteq x\}$ is $\Sigma^1_1$,
        \item if $A\subseteq \bairespace\cross\bairespace$ is $\Sigma^1_1$,
     then so is
     $$B=\{x\in\bairespace: \exists y\in\omega \;(x,y)\in A\},$$
        \item if $A\subseteq \omega\cross\bairespace$ is $\Sigma^1_1$,
     then so is
     $$B=\{x\in\bairespace: \exists n\in\omega \;(n,x)\in A\},$$
        \item if $A\subseteq \omega\cross\bairespace$ is $\Sigma^1_1$,
     then so is
     $$B=\{x\in\bairespace: \forall n\in\omega \;(n,x)\in A\},$$
        \item if $\langle A_n: n\in\omega\rangle$ is sequence of
     $\Sigma_1^1$ sets given by the recursive predicates
     $R_n$ and $\langle R_n: n\in\omega\rangle$ is (uniformly)
      recursive, then both
        $$\Union_{n\in\omega} A_n  \rmand  \Intersect_{n\in\omega} A_n
         \mbox{ are } \Sigma_1^1.$$
        \item if $f:\bairespace\arrow\bairespace$ is a function
      whose graph is $\Sigma^1_1$ and
      $A\subseteq \bairespace$ is $\Sigma^1_1$, then  so is $f^{-1}(A)$.
   \end{enumerate}
\end{lemma}

Of course, the above lemma is true with $\omega$ or $\omega\cross\bairespace$,
etc., in place of $\bairespace$.  It also relativizes to
any class $\Sigma^1_1(x)$.  It follows
from the Lemma that every Borel subset of $\bairespace$ is $\bsig11$
and that the continuous pre-image of $\bsig11$ set is $\bsig11$.

\begin{theorem}
  There exists a $\Sigma^1_1$ set $U\subseteq \bairespace\cross\bairespace$
  which is universal for all $\bsig11$ sets, i.e., for every
  $\bsig11$ set $A\subseteq\bairespace$ there exists $x\in\bairespace$
  with $$A=\{y:(x,y)\in U\}.$$ \sdex{universal for $?bsig11$ sets}
\end{theorem}
\proof
There exists $C\subseteq\bairespace\cross\bairespace\cross\bairespace$
a $\Pi_1^0$ set which is universal for $\bpi01$ subsets
of $\bairespace\cross\bairespace$.  Let $U$ be the projection of
$C$ on its second coordinate.
\qed

Similarly we can get $\Sigma^1_1$ sets contained in  $\omega\cross\omega$
(or $\omega\cross\bairespace$) which are universal for $\Sigma^1_1$ subsets
of $\omega$ (or $\bairespace$).

The usual diagonal argument shows that
there are $\Sigma^1_1$ subsets of $\bairespace$ which are not $\bpi11$ and
$\Sigma^1_1$ subsets of $\omega$ which are not $\Pi^1_1$.

\begin{theorem}\label{normal} \sdex{Normal form for $\Sigma^1_1$}
  (Normal form) A set $A\subseteq\bairespace$
  is $\Sigma^1_1$ iff there exists a recursive map
    $$\bairespace\arrow 2^{\finseq}\;\; x\mapsto T_x$$
  such that $T_x\subseteq\finseq$ is a tree for every $x\in\bairespace$,
  and $x\in A$ iff $T_x$ is ill-founded. By recursive map we
  mean that there is a Turing machine $\{e\}$ such that
  for $x\in\bairespace$ the machine $e$ computing with an oracle
  for $x$, $\{e\}^x$ computes the characteristic function of $T_x$.
\end{theorem}
\proof
Suppose
$$x\in A \rmiff \exists y \in\bairespace\;\forall n\in\omega\;\;
R(x\res n, y\res n).$$

Define
$$T_x=\{s\in \finseq: \forall i \leq |s|\;\; R(x\res i, s\res i)\}.$$
\qed
A similar thing is true for $A\subseteq\omega$, i.e.,
$A$ is $\Sigma^1_1$ iff there is a uniformly recursive list
of recursive trees $\langle T_n:n<\omega\rangle$ such that
$n\in A$ iff $T_n$ is ill-founded.

The connection between $\Sigma^1_1$ and well-founded trees, gives
us the following:

\begin{theorem} \label{mostow}
  (Mostowski's Absoluteness) Suppose $M\subseteq N$ are
  transitive models of ZFC$^*$ and $\theta$ is $\bsig11$ sentence
  with parameters in $M$.   Then \sdex{Mostowski's Absoluteness}
  $M\models\theta$ iff $N\models\theta$.
\end{theorem}
\proof
ZFC$^*$ is a nice enough finite fragment of ZFC to know that
trees are well-founded iff they have rank functions (Theorem \ref{tree}).
$\theta$ is $\bsig11$ sentence \sdex{ZFC$^*$}
  with parameters in $M$ means there exists $R$ in $M$ such that
  $$\theta=\exists x\in\bairespace \forall n \; R(x\res n).$$
This means that for some tree $T\subseteq\finseq$ in $M$
$\theta$ is equivalent to ``$T$ has an infinite branch''.
So if $M\models\theta$ then $N\models\theta$ since a branch
$T$ exists in $M$.  On the other hand if
$M\models\neg\theta$, then
$$M\models\exists r:T\arrow\ordinals\mbox{ a rank function''}$$
and then for this same $r\in M$
$$N\models r:T\arrow\ordinals\mbox{ is a rank function''}$$
and so $N\models\neg\theta$. \qed

\sector{Constructible well-orderings}

G\"odel proved the axiom of choice relatively consistent with ZF
by producing a definable well-order of the constructible universe.
He announced in G\"odel \site{god} that if V=L,
then there exists an uncountable $\Pi^1_1$ set without perfect subsets.
Kuratowski wrote down a
proof of the theorem below but the manuscript was lost during
World War II (see Addison \site{addison2}).

A set is $\Sigma^1_2$ iff it is the projection of a $\Pi_1^1$ set.
\sdex{$\Sigma^1_2$}

\begin{theorem} \label{delluz}
  [V=L] There exists a $\Delta^1_2$
  well-ordering of $\bairespace$.     \sdex{$\Delta^1_2$ well-ordering}
\end{theorem} \sdex{V=L}
\proof

Recall the definition of G\"odel's Constructible sets $L$.
$L_0=\emptyset$, $L_\lambda=\Union_{\alpha<\lambda}L_\alpha$ for
$\lambda$ a limit ordinal, and $L_{\alpha+1}$
is the definable subsets of $L_\alpha$.  Definable means
with parameters from $L_\alpha$.  $L=\Union_{\alpha\in\ordinals}L_\alpha$.

The set $x$ is constructed before $y$, \sdex{$x<_cy$}
($x<_cy$) iff the least $\alpha$ such that $x\in L_\alpha$ is less
than the least $\beta$ such that $y\in L_\beta$, or $\alpha=\beta$
and the ``least'' defining formula for $x$ is less than the one
for $y$.   Here ``least'' basically boils down to lexicographical order.
Whatever the exact formulation of $x<_c y$ is it satisfies:
$$x<_cy\rmiff L_\alpha\models x<_cy$$
where $x,y\in L_\alpha$ and $L_\alpha\models$ZFC$^*$ where
ZFC$^*$ is a sufficiently large finite fragment of ZFC.  (Actually,
it is probably enough for $\alpha$ to a limit ordinal.)
Assuming $V=L$,  for $x,y\in\bairespace$ we have that $x<_cy$ iff
there exists $E\subseteq\omega\cross\omega$ and $\name{x},\name{y}\in\omega$
such that letting $M=(\omega,E)$ then
\begin{enumerate}
  \item $E$ is extensional and well-founded,
  \item $M\models$ZFC$^*$+ V=L
  \item $M\models \name{x}<_c\name{y}$,
  \item for all $n,m\in\omega\;\;\;(x(n)=m$
  iff $M\models \name{x}(\name{n})=\name{m}$), and
 \item for all $n,m\in\omega$ ($y(n)=m$ iff
 $M\models \name{y}(\name{n})=\name{m}$).
\end{enumerate}
The first clause guarantees (by the Mostowski collapsing lemma) that
$M$ is isomorphic to a transitive set. The second, that this transitive
set will be of the form $L_\alpha$.  The last two clauses guarantee
that the image under the collapse of $\name{x}$ is $x$ and
$\name{y}$ is $y$.

Well-foundedness of $E$ is $\Pi^1_1$.  The remaining clauses are
all $\Pi^0_n$ for some $n\in\omega$.  Hence, we have given a
$\Sigma^1_2$ definition of $<_c$.  But
$x\not <_c y$ iff $y=x$ or $y<_cx$.  It follows that $<_c$ is
also $\Pi^1_2$ and hence $\Delta^1_2$.

\qed

\sector{Hereditarily countable sets}

HC is the set of all hereditarily countable sets.  There is a close
\sdex{HC}\sdex{hereditarily countable sets}
connection between the projective hierarchy above level 2 and
a natural hierarchy on the subsets of HC.  A formula of
set theory is $\Delta_0$ iff it is in the smallest family
\sdex{$\Delta_0$-formulas}
of formulas containing the atomic formulas of the form ``$x\in y$'' or
``$x=y$'', and closed under conjunction, $\theta\conj\phi$,
negation, $\neg\theta$, and bounded quantification,
$\forall x\in y$ or $\exists x\in y$.  A formula $\theta$ of set
theory is $\Sigma_1$ iff it of the form $\exists u_1,\ldots,u_n\;\; \psi$
where $\psi$ is $\Delta_0$. \sdex{$\Sigma_1$-formula}

\begin{theorem}
  A set $A\subseteq\bairespace$ is $\Sigma^1_2$ iff
  there exists a $\Sigma_1$ formula $\theta(.)$ of set theory
  such that $$A=\{x\in\bairespace: HC\models \theta(x)\}.$$
\end{theorem}
\proof
We note that $\Delta_0$ formulas are absolute between transitive
sets, i.e., if $\psi(\ldots)$ is $\Delta_0$ formula, $M$ a transitive
set and $\overline{y}$ a finite sequence of elements of $M$, then
$M\models\psi(\overline{y})$ iff  $V\models\psi(\overline{y})$.
Suppose that $\theta(.)$ is a $\Sigma_1$ formula of set theory.
Then for every $x\in\bairespace$ we have that
$HC\models\theta(x)$ iff there exists a countable transitive
set $M\in HC$ with $x\in M$ such that $M\models\theta(x)$.
Hence, $HC\models\theta(x)$ iff
there exists $E\subseteq\omega\cross\omega$ and $\name{x}\in\omega$
such that letting $M=(\omega,E)$ then
\begin{enumerate}
  \item $E$ is extensional and well-founded,
  \item $M\models$ZFC$^*$, (or just that $\omega$ exists)
  \item $M\models \theta(\name{x})$,
  \item for all $n,m\in\omega\;\;\;x(n)=m$
  iff $M\models \name{x}(\name{n})=\name{m}$.
\end{enumerate}
Therefore, $\{x\in\bairespace: HC\models \theta(x)\}$ is a $\Sigma^1_2$
set.  On the other hand given a $\Sigma^1_2$ set $A$ there
exists a $\Pi_1^1$ formula $\theta(x,y)$  such that
$A=\{x:\exists y \theta(x,y)\}$.  But then by Mostowski absoluteness
(Theorem \ref{mostow}) we have that
$x\in A$ iff there exists a countable transitive set $M$ with
$x\in M$ and there exists $y\in M$ such that $M\models$ZFC$^*$
and $M\models\theta(x,y)$.  But this is a $\Sigma_1$ formula
for HC.\qed

The theorem says that    $\Sigma^1_{2}=\Sigma^{HC}_1$.
\sdex{$\Sigma^1_{2}=\Sigma^{HC}_1$}
Similarly, $\Sigma^1_{n+1}=\Sigma^{HC}_n$.
Let us illustrate this with an example construction.

\begin{theorem}
  If V=L, then there exists a $\Delta^1_2$ Luzin set $X\subseteq\bairespace$.
\end{theorem}
\proof
Let $\{T_\alpha:\alpha<\omega_1\}$ (ordered by $<_c$) be all subtrees $T$ of
$\finseq$ whose branches $[T]$ are a closed nowhere dense
subset of $\bairespace$. Define $x_\alpha$ to be the least
constructed ($<_c$) element
of $\bairespace$ which is not in
$$\Union_{\beta<\alpha}[T_\beta]\union\{x_\beta:\beta<\alpha\}.$$
Define $X=\{x_\alpha:\alpha<\omega_1\}$.  So $X$ is a Luzin set.

To see that $X$ is $\Sigma_1^{HC}$ note that $x\in X$ iff
there exists a transitive countable $M$ which models ZFC$^*$+V=L such
that $M\models$``$x\in X$''(i.e. $M$ models the first paragraph
of this proof).

To see that $X$ is $\Pi_1^{HC}$ note that $x\in X$ iff
for all $M$ if $M$ is a transitive countable model of ZFC$^*$+V=L
with $x\in M$ and
 $M\models$``$\exists y\in X \;\; x<_cy$'', then
$M\models$``$x\in X$''.
This is true because the nature of the construction is such that
if you put a real into $X$ which is constructed after $x$, then
$x$ will never get put into $X$ after this.  So $x$ will be in
$X$ iff it is already in $X$.

\qed

\sector{Shoenfield Absoluteness}

 For a tree $T\subseteq \Union_{n<\omega}\kappa^n\cross\omega^n$
 define
 $$p[T]=\{y\in\bairespace: \exists x\in \kappa^\omega\;\;\forall n
 (x\res n, y\res n)\in T\}.$$
 A set defined this way is called $\kappa$-Souslin.  Thus
 \sdex{$\kappa$-Souslin}
 $\bsig11$ sets are precisely the $\omega$-Souslin sets.  Note that
 if $A\subseteq\bairespace\cross\bairespace$ and $A=p[T]$ then the projection
 of $A$, $\{y:\exists x\in\bairespace \;(x,y)\in A\}$ is
 $\kappa$-Souslin.  To see this let $<,>:\kappa\cross\omega\arrow\kappa$
 be a pairing function.  For $s\in \kappa^n$ let $s_0\in\kappa^n$
 and $s_1\in\omega^n$ be defined by $s(i)=<s_0(i),s_1(i)>$.
 Let $T^*$ be the tree defined by
$$T^*=\Union_{n\in\omega}\{(s,t)\in \kappa^n\cross\omega^n:
 (s_0,s_1,t)\in T\}.$$
Then  $p[T^*]=\{y:\exists x\in\bairespace \;(x,y)\in A\}$.

\begin{theorem} \label{shoentree}
  (Shoenfield \site{shoen}) If $A$ is a $\Sigma^1_2$ set, then $A$ is
  $\omega_1$-Souslin set coded in L, i.e. $A=p[T]$
  where $T\in L$.
\end{theorem}
\proof
From the construction of $T^*$ it is clear that is enough to
see this for $A$ which is $\Pi^1_1$.

We know that a  countable tree is well-founded
iff there exists a rank function $r:T\arrow\omega_1$.
Suppose
$$x\in A\rmiff \forall y \exists n \;\;(x\res n,y\res n)\notin T$$
where $T$ is a recursive tree. So defining
$T_x=\{t: (x\res |t|,t)\in T\}$ we have that
$x\in A$ iff $T_x$ is well-founded (Theorem \ref{normal}).

The $\omega_1$ tree $\hat{T}$ is just the tree of partial rank functions.
Let $\{s_n:n\in\omega\}$ be  a recursive listing of $\finseq$ with
$|s_n|\leq n$. Then for every $N<\omega$,
and $(r,s)\in \omega_1^N\cross\omega^N$ we have $(r,t)\in\hat{T}$ iff
$$ \forall n,m<N\;
[(t,s_n),(t,s_m)\in T\rmand s_n\subset s_m] \mbox{ implies } r(n)>r(m) .$$
Then $A=p[\hat{T}]$.\qed

\begin{theorem} \label{shoenabs} \sdex{Shoenfield Absoluteness}
  (Shoenfield Absoluteness \site{shoen}) If $M\subseteq N$ are transitive models
  of ZFC$^*$ and $(\omega_1)^N\subseteq M$, then for any $\Sigma^1_2(x)$
  sentence $\theta$ with parameter $x\in M$.
  $$M\models\theta \rmiff N\models\theta.$$
\end{theorem}
\proof
If $M\models\theta$, then $N\models\theta$, because $\bpi11$ sentences
are absolute.  On the other hand suppose $N\models\theta$.  Working
in $N$ using the proof of Theorem \ref{shoentree} we get
a tree $T\subseteq\omega_1^{<\omega}$ with
$T\in L[x]$ such that $T$ is ill-founded, i.e., there exists
$r\in [T]$.  Note that $r$ codes a witness to a $\Pi^1_1(x)$
predicate and a rank function showing the tree corresponding to
this predicate is well-founded. Since for some
$\alpha<\omega_1$, $r\in\alpha^\omega$ we see that
$$T_\alpha=T\intersect \alpha^{<\omega}$$
is ill-founded. But $T_\alpha\in M$ (since by assumption
$(\omega_1)^N\subseteq M$) and so by the absoluteness of well-founded
trees, $M$ thinks that $T_\alpha$ is ill-founded.  But a branch
thru $[T]$ gives a witness and a rank function showing that $\theta$ is
true, and consequently, $M\models\theta$.\qed

\sector{Mansfield-Solovay Theorem}

\begin{theorem} \label{mansol}
  (Mansfield \site{mansf}, Solovay \site{solov})
  If $A\subseteq\bairespace$ is a $\bsig12$ set with constructible
  parameter which contains a nonconstructible element of $\bairespace$,
  then $A$ contains a perfect set which is coded in $L$.
\end{theorem} \sdex{Mansfield-Solovay Theorem}

\proof
By Shoenfield's Theorem \ref{shoentree}, we may assume
$A=p[T]$ where $T\in L$ and
$T\subseteq \Union_{n<\omega}\omega_1^n\cross\omega^n$.
Working in $L$ define the following decreasing sequence of
subtrees as follows.

$T_0=T$,

$T_\lambda=\Intersect_{\beta<\lambda}T_\beta$,  if $\lambda$ a limit ordinal, and

$T_{\alpha+1}=\{(r,s)\in T_\alpha: \exists (r_0,s_0), (r_1,s_1)\in T_\alpha$
such that
$(r_0,s_0), (r_1,s_1)$ extend $(r,s)$, and $s_0$ and $s_1$ are incompatible$\}$.

Each $T_\alpha$ is tree, and for $\alpha<\beta$ we have
$T_\beta\subseteq T_\alpha$.  Thus there exists some $\alpha_0$ such
that $T_{\alpha_0+1}=T_{\alpha_0}$.

\claim  $[T_{\alpha_0}]$ is nonempty.

\proof
Let $(x,y)\in [T]$ be any pair with $y$ not constructible.  Since
$A=p[T]$ and $A$ is not a subset of $L$, such a pair must exist.
Prove by induction on $\alpha$ that $(x,y)\in [T_\alpha]$.
This is easy for $\alpha$ a limit ordinal. So suppose
$(x,y)\in [T_\alpha]$ but $(x,y)\notin [T_{\alpha+1}]$. By the
definition it must be that there exists $n<\omega$ such
that $(x\res n,y\res n)=(r,s) \notin T_{\alpha+1}$. But
in $L$ we can define the tree:
$$T_{\alpha}^{(r,s)}=\{(\hat{r},\hat{s})\in T_{\alpha}:
(\hat{r},\hat{s})\subseteq
(r,s) \rmor (r,s) \subseteq (\hat{r},\hat{s})\}$$
which has the property that $p[T_{\alpha}^{(r,s)}]=\{y\}$. But by
absoluteness of well-founded trees, it must be that
there exists $(u,y_0)\in [T_{\alpha}^{(r,s)}]$ with $(u,y_0)\in L$.
But then $y_0=y\in L$ which is a contradiction.
This proves the claim. \qed

Since $T_{\alpha_0+1}=T_{\alpha_0}$, it follows that
for every $(r,s)\in T_{\alpha_0}$ there exist
$$(r_0,s_0), (r_1,s_1)\in T_{\alpha_0}$$ such that
$(r_0,s_0), (r_1,s_1)$ extend $(r,s)$ and $s_0$ and $s_1$ are incompatible.
This allows us to build by induction (working in $L$):
$$\langle (r_\sigma,s_\sigma):\sigma\in 2^{<\omega}\rangle$$
with $(r_\sigma,s_\sigma)\in T_{\alpha_0}$ and for each
$\sigma\in 2^{<\omega}$
$(r_{\sigma_0},s_{\sigma_0}), (r_{\sigma_1},s_{\sigma_1})$ extend
$(r_{\sigma},s_{\sigma})$ and $s_{\sigma_0}$ and $s_{\sigma_1}$ are
incompatible.  For any $q\in\cantorspace$ define
$$x_q=\Union_{n<\omega}r_{q\res n}
\rmand y_q=\Union_{n<\omega}s_{q\res n}.$$
Then we have that $(x_q,y_q)\in [T_{\alpha_0}]$ and
therefore $P=\{y_q:q\in \cantorspace\}$ is a perfect set such that
$$P\subseteq p[T_{\alpha_0}]\subseteq p[T]=A$$ and $P$ is coded
in $L$.\qed

This proof is due to Mansfield.  Solovay's
proof used forcing.  Thus we have
departed\footnote{``Consistency is the hobgoblin of little minds.
 With consistency
a great soul has simply nothing to do.''  Ralph Waldo Emerson.}
from our theme of giving forcing proofs.

\sector{Uniformity and Scales}

Given $R\subseteq X\cross Y$ we say that $S\subseteq X\cross Y$ uniformizes
$R$ iff
\begin{enumerate}
  \item $S\subseteq R$,
  \item for all $x\in X$ if there exists $y\in Y$ such that
  $R(x,y)$, then there exists $y\in Y$ such that $S(x,y)$, and
  \item for all $x\in X$ and $y,z\in Y$ if $S(x,y)$ and $S(x,z)$, then
  $y=z$.
\end{enumerate}
\sdex{uniformization property}
Another way to say the same thing is that $S$ is a subset of $R$ which
is the graph of a function whose domain is the same as $R$'s.

\begin{theorem} \label{uniform}
  (Kondo \site{kondo})
  Every $\Pi^1_1$ set $R$ can be uniformized
  by a $\Pi^1_1$ set $S$.
\end{theorem}

Here, $X$ and $Y$ can be taken to be either $\omega$ or $\bairespace$ or
even a singleton $\{0\}$.  In this last case, this amounts to
saying for any nonempty $\Pi_1^1$ set $A\subseteq\bairespace$ there
exists a $\Pi^1_1$ set $B\subseteq A$ such that $B$ is a singleton,
i.e., $|B|=1$.  The proof of this Theorem is to use a property which
has become known as the Scale property.

\begin{lemma}\label{scale} \sdex{scale property}
  (Scale property) For any $\Pi^1_1$ set $A$ there exists
  $\langle \phi_i:i<\omega\rangle$ such that
  \begin{enumerate}
    \item each $\phi_i:A\arrow\ordinals$,
    \item for all $i$ and $x,y\in A$ if $\phi_{i+1}(x)\leq\phi_{i+1}(y)$, then
          $\phi_{i}(x)\leq\phi_{i}(y)$,
    \item for every $x,y\in A$ if $\forall i \;\;\phi_i(x)=\phi_i(y)$, then $x=y$,
    \item for all $\langle x_n:n<\omega\rangle\in A^{\omega}$ and
              $\langle \alpha_i:i<\omega\rangle\in \ordinals^{\omega}$ if
              for every $i$ and for all but finitely many
              $n\;\;\phi_i(x_n)=\alpha_i$, then there exists $x\in A$ such
              that $\lim_{n\arrow\infty}x_n=x$ and for each $i$
              $\phi_i(x)\leq\alpha_i$,
     \item  there exists $P$ a $\Pi^1_1$ set such that
             for all $x,y\in A$ and $i$
             $$P(i,x,y)\rmiff \phi_i(x)\leq \phi_i(y)$$
             and
             for all $x\in A$, $y\notin A$, $i\in\omega$ $P(i,x,y)$, and
     \item  there exists $S$  a $\Sigma^1_1$ set such that
             for all $x,y\in A$ and $i$
             $$S(i,x,y)\rmiff \phi_i(x)\leq \phi_i(y)$$
             and
             for all $y\in A$, $x$, $i\in\omega$
             if $S(i,x,y)$, then $x\in A$.
  \end{enumerate}
\end{lemma}

Another way to view a scale is from the point of view of
the relations on $A$ defined by $x\leq_i y$ iff $\phi_i(x)\leq\phi)i(y)$.
These are called \sdex{prewellorderings} prewellorderings.
They are well orderings
if we mod out by $x\equiv_iy$ which is defined by
$$x\equiv_iy\rmiff x\leq_iy \rmand y\leq_ix . $$
The second item says that
these relations get finer and finer as $i$ increases.
The third item says that in the ``limit'' we get a linear order.
The fourth item is some sort of continuity condition.
And the last two items are the definability properties of the
scale.

Before proving the lemma, let us deduce uniformity from it.
We do not use the last item in the lemma.  First let us show that
for any nonempty $\Pi_1^1$ set $A\subseteq\bairespace$ there
exists a $\Pi^1_1$ singleton  $B\subseteq A$.  \sdex{$\Pi^1_1$ singleton}
Define
$$x\in B\rmiff x\in A\rmand\forall n \forall y\; P(n,x,y).$$
Since $P$ is $\Pi^1_1$ the set $B$ is $\Pi_1^1$.  Clearly
$B\subseteq A$, and also by item (2) of the lemma,  $B$ can have at
most one element.  So it remains to show that $B$ is nonempty.
Define $\alpha_i=\min\{\phi_i(x):x\in A\}$.  For each
$i$ choose $x_i\in A$ such that $\phi_i(x_i)=\alpha_i$.

\claim If $n>i$ then $\phi_i(x_n)=\alpha_i$.
\proof
By choice of $x_n$ for every $y\in A$ we have $\phi_n(x_n)\leq \phi_n(y)$.
By item (2) in the lemma, for every $y\in A$ we have that
$\phi_i(x_n)\leq \phi_i(y)$ and hence $\phi_i(x_n)=\alpha_i$.
\qed

By item (4) there exists $x\in A$ such that $\lim_{n\arrow\infty}x_n=x$ and
$\phi_i(x)\leq\alpha_i$ all $i$. By the minimality of $\alpha_i$ it
must be that $\phi_i(x)=\alpha_i$.  So $x\in B$ and we are done.

Now to prove  a more general case of uniformity suppose that
$R\subseteq\bairespace\cross\bairespace$ is $\Pi^1_1$. Let
$\phi_i:R\arrow \ordinals$ be scale given by
the lemma and
$$P\subseteq \omega\cross (\bairespace\cross\bairespace)\cross
(\bairespace\cross\bairespace)$$ be the $\Pi_1^1$ predicate given by
item (5).  Then define the $\Pi^1_1$ set
$S\subseteq \bairespace\cross\bairespace$ by
$$(x,y)\in S\rmiff
(x,y)\in R\rmand\forall z \forall n\; P(n,(x,y),(x,z)).$$
The same proof shows that $S$ uniformizes $R$.

The proof of the Lemma will need the following two elementary facts
about well-founded trees.  For $T,\hat{T}$ subtrees of $Q^{<\omega}$ we
say that $\sigma:T\arrow \hat{T}$ is a \dex{tree embedding} iff
for all $s,t\in T$ if $s\subset t$ then $\sigma(s)\subset \sigma(t)$.
Note that $s\subset t$ means that $s$ is a proper initial segment of $t$.
\sdex{$s\subset t$}
Also note that tree embeddings need not be one-to-one. We write
$T\embed\hat{T}$ iff there exists a tree embedding from $T$ into $\hat{T}$.
\sdex{$T?embed?hat{T}$}
We write $T\strictembed\hat{T}$ iff there is a tree embedding which
takes the root node of $T$ to a nonroot node of $\hat{T}$.
\sdex{$T?strictembed\hat{T}$}
Recall that $r:T\arrow\ordinals$ is a rank function iff
for all $s,t\in T$ if $s\subset t$ then $r(s)>r(t)$.  Also
the rank of $T$ is the minimal ordinal $\alpha$ such that
there exists a rank function $r:T\arrow\alpha+1$.

\begin{lemma} \label{emb1}
Suppose $T\embed\hat{T}$ and $\hat{T}$ is well-founded,
i.e., $[\hat{T}]=\emptyset$, then $T$ is well-founded and rank of
$T$ is less than or equal to rank of $\hat{T}$.
\end{lemma}
\proof
Let $\sigma:T\arrow\hat{T}$ be a tree embedding and
$r:\hat{T}\arrow\ordinals$ a rank function.  Then
$r\circ \sigma$ is a rank function on $T$.
\qed

\begin{lemma}\label{emb2}
Suppose $T$ and $\hat{T}$ are well founded trees and
rank of $T$ is less than or equal rank of $\hat{T}$, then
$T\embed \hat{T}$.
\end{lemma}
Let $r_T$ and $r_{\hat{T}}$ be the canonical rank functions on
$T$ and $\hat{T}$ (see Theorem \ref{tree}).  Inductively
define $\sigma:T\intersect Q^{n}\arrow\hat{T}\intersect Q^{n}$,
so as to satisfy $r_T(s)\leq r_{\hat{T}}(\sigma(s))$.\qed

Now we a ready to prove the existence of scales (Lemma \ref{scale}).
Let $$\omega_1^{-}=\{-1\}\union\omega_1$$ be well-ordered in the obvious
way.  Given a well-founded tree $T\subseteq\finseq$ with rank function
$r_T$ extend $r_T$ to all of $\finseq$ by defining $r_T(s)=-1$ if
$s\notin T$.  Now suppose $A\subseteq\bairespace$ is $\Pi_1^1$ and
$x\in A$ iff $T_x$ is well-founded (see Theorem \ref{normal}).
Let $\{s_n:n<\omega\}$ be a recursive listing of $\finseq$ with
$s_0=\emptyseq$. For each $n<\omega$ define
$\psi_n:A\arrow \omega_1^-\cross\omega\cross\cdots\omega_1^-\cross\omega$
by
$$\psi_n(x)=\langle r_{T_x}(s_0), x(0), r_{T_x}(s_1), x(1),\ldots,
r_{T_x}(s_{n}), x(n)\rangle.$$
The set $\omega_1^-\cross\omega\cross\cdots\omega_1^-\cross\omega$ is
well-ordered by the lexicographical order.  The scale
$\phi_i$ is just obtained by mapping the range of $\psi_i$ order
isomorphically
to the ordinals.  (Remark: by choosing $s_0=\emptyseq$, we guarantee that
 the first
coordinate is always the largest coordinate, and so
the range of $\psi_i$ is less than or equal to $\omega_1$.)
Now we verify the properties.

For item (2):  if $\psi_{i+1}(x)\leq_{lex}\psi_{i+1}(y)$, then
$\psi_{i}(x)\leq_{lex}\psi_{i}(y)$. This is true because we
are just taking the lexicographical order of a longer sequence.

For item (3):  if $\forall i \;\;\psi_i(x)=\psi_i(y)$, then $x=y$.
This is true, because $\psi_i(x)=\psi_i(y)$ implies $x\res i=y\res i$.

For item (4):
Suppose $\langle x_n:n<\omega\rangle\in A^{\omega}$ and
for every $i$ and for all but finitely many
  $n\;\;\psi_i(x_n)=t_i$.  Then since $\psi_i(x_n)$
  contains $x_n\res i$ there must be
 $x\in\bairespace$ such that $\lim_{n\arrow\infty}x_n=x$. Note
 that since $\{s_n:n\in\omega\}$ lists every element of $\finseq$,
 we have that for every $s\in \finseq$ there exists $r(s)\in\ordinals$
 such that $r_{T_{x_n}}(s)=r(s)$ for all but finitely many $n$.  Using
 this and  $$\lim_{n\arrow\infty}T_{x_n}=T_x$$ it follows that
 $r$ is a rank function on $T_x$.  Consequently $x\in A$.  Now
 since $r_{T_x}(s)\leq r(s)$, it follows that
 $\psi_i(x)\leq_{lex}t_i$.

For item (5),(6): The following set is $\Sigma^1_1$:
$$\{(T,\hat{T}): T,\hat{T} \mbox{ are subtrees of }\finseq,
T\embed \hat{T}\}.$$
Consequently, assuming that $T,\hat{T}$ are well-founded, to
say that $r_T(s)\leq r_{\hat{T}}(s)$ is equivalent to
saying there exists a tree embedding which takes $s$ to $s$.
Note that this is
$\Sigma^1_1$.  This shows that it is possible to define a $\Sigma^1_1$
set $S\subseteq \omega\cross\bairespace \cross\bairespace$
such that for every
$x,y\in A$ we have $(n,x,y)\in S$ iff
$$\langle r_{T_x}(s_0), x(0), r_{T_x}(s_1), x(1),\ldots,
r_{T_x}(s_{n}), x(n)\rangle$$
is lexicographically less than or equal to
$$\langle r_{T_y}(s_0), y(0), r_{T_y}(s_1), y(1),\ldots,
r_{T_y}(s_{n}), y(n)\rangle.$$
Note that if $(n,x,y)\in S$ and $y\in A$ then $x\in A$, since
$T_y$ is a well-founded tree and $S$ implies $T_x\embed T_y$, so
$T_x$ is well-founded and so $x\in A$.

To get the $\Pi_1^1$ relation $P$ (item (5)), instead of
saying $T$ can be embedded into $\hat{T}$ we say that
$\hat{T}$ cannot be embedded properly into $T$,i.e.,
 $\hat{T}\not\strictembed T$ or in other words,
there
does not exists a tree embedding $\sigma:\hat{T}\arrow T$ such
that $\sigma(\emptyseq)\not=\emptyseq$.  This is a $\Pi^1_1$ statement.
For $T$ and $\hat{T}$ well-founded trees saying that
rank of $T$ is less than or equal to $\hat{T}$ is equivalent to saying
rank of $\hat{T}$ is not strictly smaller than the rank of $T$.
But by Lemma \ref{emb2} this is equivalent to the nonexistence
of such an embedding.  Note also that if $x\in A$ and
$y\notin A$, then we will have $P(n,x,y)$ for every $n$.
This is because $T_y$ is not well-founded and so
cannot be embedded into the well-founded tree $T_x$.

This finishes the proof of the Scale Lemma \ref{scale}.
\qed

The scale property was invented by Moschovakis \site{playful}
to show how determinacy
could be used to get
uniformity properties\footnote{I have yet to see any problem,
however complicated, which, when you looked
at it in the right way, did not become
still more complicated. Poul Anderson}
in the higher projective classes.
He was building on earlier ideas of Blackwell, Addison, and
Martin.  The 500 page book by Kuratowski and Mostowski \site{kurmos}
ends with a proof of the uniformization theorem.

\sector{Martin's axiom and Constructibility}

\begin{theorem}
  (G\"{o}del see Solovay \site{solov}) If V=L, there exists
  uncountable $\Pi_1^1$ set $A\subseteq\bairespace$ which
  contains no perfect subsets.
\end{theorem}
\proof
Let $X$ be any uncountable $\Sigma^1_2$ set containing no
perfect subsets.  For example, a $\Delta^1_2$ Luzin set
would do (Theorem \ref{delluz}).  Let $R\subset\bairespace\cross\bairespace$
be $\Pi_1^1$ such that $x\in X$ iff $\exists y\; R(x,y)$.
Use  $\Pi^1_1$ uniformization (Theorem \ref{uniform}) to get
$S\subseteq R$ with the property that
$X$ is the one-to-one image of $S$ via the projection map
$\pi(x,y)=x$.  Then $S$ is an uncountable $\Pi^1_1$ set which
contains no perfect subset.  This is because if $P\subseteq S$ is
perfect, then $\pi(P)$ is a perfect subset of $X$.\qed

Note that it is sufficient to assume that $\omega_1=(\omega_1)^L$.
Suppose  $A\in L$ is defined by the $\Pi^1_1$ formula $\theta$.
Then let $B$ be the set which is defined by $\theta$ in $V$.
So by $\Pi^1_1$ absoluteness $A=B\intersect L$. The set
$B$ cannot contain a perfect set since the sentence:
$$\exists T\mbox{ $T$ is a perfect tree and } \forall x\;
(x\in [T]\mbox{ implies } \theta(x))$$
is a $\Sigma^1_2$ and false in $L$ and so by Shoenfield absoluteness
(Theorem \ref{shoenabs}) must be false in $V$.  It follows then
by the Mansfield-Solovay Theorem \ref{mansol} that $B$ cannot contain
a nonconstructible real and so $A=B$.

Actually, by tracing thru the actual definition of $X$ one can
see that the elements of the uniformizing set $S$ (which is what $A$ is)
consist
of pairs $(x,y)$ where $y$ is isomorphic to some $L_\alpha$ and
$x\in L_\alpha$.  These pairs are reals which witness
their own constructibility,
so one can avoid using the Solovay-Mansfield Theorem.

\begin{corollary}
  If $\omega_1=\omega_1^L$, then there exists a $\Pi^1_1$ set
  of constructible reals which contains no perfect set.
\end{corollary}

\begin{theorem}
  (Martin-Solovay \site{marsol}) Suppose
  MA + $\neg$CH + $\omega_1=(\omega_1)^L$.
  Then every $A\subseteq\cantorspace$ of cardinality $\omega_1$ is
  $\bpi11$.
\end{theorem}
\proof
Let $A\subseteq\cantorspace$ be a uncountable $\Pi^1_1$ set of constructible
reals and let $B$ be an arbitrary subset of $\cantorset$ of cardinality
$\omega_1$.  Arbitrarily well-order the two sets,
$A=\{a_\alpha:\alpha<\omega_1\}$ and $B=\{b_\alpha:\alpha<\omega_1\}$.

By Theorem \ref{silverlem} there exists two sequences
of $G_\delta$ sets $\langle U_n:n<\omega\rangle$ and
$\langle V_n:n<\omega\rangle$ such that
for every $\alpha<\omega_1$ for every $n<\omega$
$$a_\alpha(n)=1\rmiff b_\alpha\in U_n$$
and
$$b_\alpha(n)=1\rmiff a_\alpha\in V_n.$$
This is because the set $\{a_\alpha: b_\alpha(n)=1\}$, although
it is an arbitrary subset of $A$, is relatively $G_\delta$ by
Theorem \ref{silverlem}.

But note that $b\in B$ iff $\forall a\in \cantorset$
$$ [\forall n \; (a(n)=1 \rmiff b\in U_n)]
                          \mbox{ implies }
                   [a\in A \rmand
                   \forall n \;(b(n)=1 \rmiff a\in U_n)].$$
Since $A$  is $\Pi^1_1$ this definition of $B$ has the form:
$$\forall a ([\bpi03]\mbox{ implies } [\Pi^1_1 \rmand \bpi03])$$
So $B$ is $\bpi11$. \qed

Note that if every set of reals of size $\omega_1$ is $\bpi11$ then
every $\omega_1$ union of Borel sets is $\bsig12$.  To see this
let $\langle B_\alpha:\alpha<\omega_1\rangle$ be any sequence of Borel sets.
Let $U$ be a universal $\bpi11$ set and let
$\langle x_\alpha:\alpha<\omega_1\rangle$ be a sequence such that
$B_\alpha=\{y:(x_\alpha,y)\in U\}$.  Then
$$y\in \Union_{\alpha<\omega_1}B_\alpha \rmiff
\exists x\;\; x \in\{x_\alpha:\alpha<\omega_1\}\conj (x,y)\in U.$$
But $\{x_\alpha:\alpha<\omega_1\}$ is $\bpi11$ and so the union
is $\bsig12$.

\sector{$\Sigma^1_2$ well-orderings}

\begin{theorem}
  (Mansfield \site{manfield}) If $(F,\lhd)$ is a $\Sigma^1_2$
  well-ordering, i.e.,
  $$F\subseteq\bairespace \rmand \lhd\subseteq F^2$$
  are both $\Sigma^1_2$, then $F$ is a
  subset of $L$.
\end{theorem}
\proof

We will use the following:

\begin{lemma}
  Assume there exists $z\in \cantorset$ such that $z\notin L$.
  Suppose $f:P\arrow F$ is a 1-1 continuous function from
  the perfect set $P$ and both $f$ and $P$ are coded in $L$,
  then there exists $Q\subseteq P$ perfect and
  $g:Q\arrow F$ 1-1 continuous so that both $g$ and $Q$ are
  coded in $L$ and for every $x\in Q$ we have
  $g(x)\lhd f(x)$.
\end{lemma}
\proof
(Kechris \site{kechris3}) First note that
there exists $\sigma:P\arrow P$ an autohomeomorphism
coded in $L$ such that for every $x\in P$ we have $\sigma(x)\not=x$ but
$\sigma^2(x)=x$.  To get this let $c:\cantorset\arrow \cantorset$ be
the complement function, i.e., $c(x)(n)=1-x(n)$ which just switches
0 and 1. Then $c(x)\not=x$ but $c^2(x)=x$.  Now if
$h:P\arrow\cantorset$ is a homeomorphism coded in $L$, then
$\sigma=h^{-1}\circ c \circ h$ works.

Now let $A=\{x\in P: f(\sigma(x))\lhd f(x)\}$.  The set $A$ is
a $\bsig12$ set with code in $L$.  Now since $P$ is coded in $L$ there
must be a $z\in P$ such that $z\notin L$.  Note that $\sigma(z)\notin L$
also.  But either
$$f(\sigma(z))\lhd f(z) \rmor f(z)=f(\sigma^2(z))\lhd f(\sigma(z))$$
and so either $z\in A$ or $\sigma(z)\in A$.
In either case $A$ has a nonconstructible member and so by the
Mansfield-Solovay Theorem \ref{mansol} the set $A$ contains a perfect
set $Q$ coded in $L$.
Let $g=f\circ \sigma$.\qed

Assume there exists $z\in F$ such that $z\notin L$.  By the
Mansfield-Solovay Theorem there exists
a perfect set $P$ coded in $L$ such that $P\subseteq F$.
Let $P_0=P$ and $f_0$ be the identity function. Repeatedly apply
the Lemma to obtain $f_n:P_n\arrow F$ so that for every
$n$ and $P_{n+1}\subseteq P_n$, for every  $x\in P_{n+1}$
$f_{n+1}(x)\lhd f_n(x)$.  But then
if $x\in\Intersect_{n<\omega}$ the sequence $\langle f_n(x):n<\omega\rangle$
is a descending $\lhd$ sequence with contradicts the fact that
$\lhd$ is a well-ordering.
\qed

Friedman \site{fried} proved the weaker result that if there
is a $\bsig12$ well-ordering of the real line, then
$\bairespace\subseteq L[g]$ for some $g\in \bairespace$.

\sector{Large $\Pi^1_2$ sets}

A set is $\bpi12$ iff it is the complement of a $\bsig12$ set.
\sdex{$?bpi12$}  Unlike $\bsig12$ sets which cannot have
size strictly in between $\omega_1$ and the continuum
(Theorem \ref{mansol}), $\bpi12$ sets can be
practically anything.\footnote{It's life Jim,  but not as we know it.-
             Spock of Vulcan}

\begin{theorem}\label{harring}
  (Harrington \site{harring}) Suppose $V$ is a model of
  set theory which satisfies $\omega_1=\omega_1^L$ and
  $B$ is arbitrary subset of $\bairespace$ in $V$.  Then
  there exists a ccc extension of $V$, $V[G]$, in which
  $B$ is a $\bpi12$ set.
\end{theorem}

\proof
Let $\poset_B$ be the following poset.  $p\in\poset_B$ iff
$p$ is a finite consistent set of sentences of the form:
\begin{enumerate}
  \item ``$[s]\intersect \name{C}_n =\emptyset$'', or
  \item ``$x\in \name{C}_n$, where $x\in B$.
\end{enumerate}
This partial order is isomorphic to Silver's view of
almost disjoint sets forcing (Theorem \ref{silverlem}).
So forcing with $\poset_B$ creates an $F_\sigma$ set
$\Union_{n\in\omega}C_n$ so that

$$\forall x\in \bairespace\intersect V
(x\in B\rmiff x\in \Union_{n<\omega}C_n).$$
Forcing with the direct sum of $\omega_1$ copies of
$\poset_B$,  $\prod_{\alpha<\omega_1}\poset_B$, we have that
$$\forall x\in \bairespace\intersect V[\langle G_\alpha:\alpha<\omega_1\rangle]
 (x\in B\rmiff
 x\in\Intersect_{\alpha<\omega_1}\union_{n<\omega}C_n^\alpha).$$
One way to see this is as follows.
Note that in any case
$$B\subseteq\Intersect_{\alpha<\omega_1}\union_{n<\omega}C_n^\alpha.$$
So it is the other implication which needs to be proved.
By ccc, for any $x\in V[\langle G_\alpha:\alpha<\omega_1\rangle]$ there
exists  $\beta<\omega_1$ with $x\in V[\langle G_\alpha:\alpha<\beta\rangle]$.
But considering $V[\langle G_\alpha:\alpha<\beta\rangle]$ as the new ground
model, then $G_\beta$ would be $\poset_B$- generic over
$V[\langle G_\alpha:\alpha<\beta\rangle]$ and hence if $x\notin B$ we
would have $x\notin\union_{n<\omega}C_n^\beta$.

Another argument will be given in the proof of the next lemma.

\begin{lemma}
Suppose $\langle c_\alpha:\alpha<\omega_1$ be a sequence in $V$
of elements of $\bairespace$ and
$\langle a_\alpha:\alpha<\omega_1\rangle$ is a sequence in
$V[\langle G_\alpha:\alpha<\omega_1\rangle]$ of elements of
$\cantorspace$.
Using Silver's forcing
add a sequence of $\bpi02$ sets $\langle U_n:n<\omega\rangle$
such that
$$\forall n\in\omega\forall\alpha<\omega_1
(a_\alpha(n)=1 \rmiff c_\alpha\in U_n).$$
Then
$$V[\langle G_\alpha:\alpha<\omega_1\rangle][\langle U_n:n<\omega\rangle]
\models
\forall x\in\bairespace\; (x\in B\rmiff
 x\in\Intersect_{\alpha<\omega_1}\union_{n<\omega}C_n^\alpha).$$
\end{lemma}
\proof
The lemma is not completely trivial, since adding the
$\langle U_n:n<\omega\rangle$ adds
new elements of $\bairespace$ which may somehow sneak into the $\omega_1$
intersection.

Working in $V$ define  $p\in \posetq$ iff
$p$ is a finite set of consistent sentences of the form:
\begin{enumerate}
  \item ``$[s]\subseteq U_{n,m}$'' where $s\in\finseq$, or
  \item ``$c_\alpha\in U_{n,m}$''.
\end{enumerate}
Here we intend that $U_n=\intersect_{m\in\omega}U_{n,m}$.
Since the $c$'s are in $V$ it is clear that the partial order $\posetq$
is too.
Define
$$\poset=\{(p,q)\in (\prod_{\alpha<\omega_1}\poset_B) \cross \posetq
:\mbox{ if ``$c_\alpha\in U_{n,m}$''$\in q$, then $p\forces a_\alpha(n)=1\}$}.$$
Note that $\poset$ is a semi-lower-lattice, i.e., if
$(p_0, q_0)$ and $(p_1, q_1)$ are compatible elements of $\poset$,
then $(p_0\union p_1, q_0\union q_1)$ is their greatest lower bound.
This is another way to view the iteration, i.e,
$\poset$ is dense in the usual iteration.
Not every iteration has
this property, one which Harrington calls ``innocuous''.

Now to prove the lemma, suppose for contradiction that
$$(p,q)\forces \name{x}\in \Intersect_{\alpha<\omega_1}\union_{n<\omega}
C_n^\alpha
\rmand \name{x}\notin B.$$
To simplify the notation, assume $(p,q)=(\emptyset,\emptyset)$.
Since $\poset$ has the ccc a sequence of
Working in $V$ let $\rangle A_n:n\in \omega\langle$ be a sequence of
 maximal antichains of $\poset$ which decide $\name{x}$, i.e.
for  $({p},{q})\in A_n$ there exists
$s\in \omega^n$ such that
$$({p},{q})\forces \name{x}\res n=\check{s}.$$
Since $\poset$ has the ccc, the $A_n$ are countable and we can find an
$\alpha<\omega$ which does not occur in the support of any ${p}$ for any
$({p},{q})$ in $\Union_{n\in\omega}A_n$.  Since $x$ is forced to be in
$\union_{n<\omega} C_n^\alpha$ there exists $(p,q)$ and $n\in\omega$
such that $$(p,q)\forces \name{x}\in C_n^\alpha.$$
Let ``$x_i\in C_n^\alpha$'' for $i<N$ be all the sentences of this type
which occur in $p(\alpha)$.  Since we are assuming $x$ is being forced not
in $B$ it must be different than
all the $x_i$, so there must be an $m$, $(\hat{p},\hat{q})\in A_m$, and
$s\in\omega^m$,
such that
\begin{enumerate}
  \item  $(\hat{p},\hat{q})$ and $(p,q)$ are compatible,
  \item $(\hat{p},\hat{q})\forces \name{x}\res m=\check{s}$, and
  \item $x_i\res m\not= s$ for every $i<N$.
\end{enumerate}
((To get $(\hat{p},\hat{q})$ and $s$ let $G$ be a generic filter containing
$(p,q)$, then since $x^G\not=x_i$ for every $i<N$ there must be
$m<\omega$ and $s\in\omega^m$ such that $x^G\res m=s$
and $s\not=x_i\res m$ for every $i<N$.
Let $(\hat{p},\hat{q})\in G\intersect A_m$.))

Now consider $(p\union\hat{p},q\union\hat{q})\in \poset$.  Since
$\alpha$ was not in the support of $\hat{p}$,
$(p\union\hat{p})(\alpha)=p(\alpha)$.  Since $s$ was chosen
so that $x_i\notin [s]$ for every $i<N$,
$$p(\alpha)\union \{[s]\intersect C^\alpha_n=\emptyset\}$$
is a consistent set of sentences, hence an element of $\poset_B$.
This is a contradiction, the condition
$$(p\union\hat{p}\union\{[s]\intersect
C^\alpha_n=\emptyset\} ,q\union\hat{q})$$
forces $x\in C^\alpha_n$ and also $x\notin C^\alpha_n$.

\qed

Let $F$ be a universal $\Sigma^0_2$ set coded in $V$ and
let $\langle a_\alpha\in\cantorspace:\alpha<\omega_1\rangle$ be
such that
$$F_{a_\alpha}=\Union_{n\in\omega }C^\alpha_n.$$
Let $C=\langle c_\alpha:\alpha<\omega_1\rangle$ be a $\Pi^1_1$ set
in $V$.  Such a set exists since $\omega_1=\omega_1^L$.

\begin{lemma}
In $V[\langle G_\alpha:\alpha<\omega_1\rangle][\langle U_n:n<\omega\rangle]$
the set $B$ is $\bpi12$.
\end{lemma}
\proof

$x\in B$ iff $x\in\Intersect_{\alpha<\omega_1}\Union_{n\in\omega }C^\alpha_n$
iff  $x\in\Intersect_{\alpha<\omega_1} F_{a_\alpha}$ iff

 $\forall a,c$ if $c\in C$ and $\forall n \;(a(n)=1 $ iff $c\in U_n$),
then $(a,x)\in F$, i.e. $(x\in F_a)$.

\bigskip

\noindent Note that
\begin{itemize}
\item ``$c\in C$'' is $\Pi^1_1$,
\item ``$\forall n \;(a(n)=1 $ iff $c\in U_n$)'' is Borel, and
\item ``$(a,x)\in F$'' is Borel,
\end{itemize}
 and  so this final definition for $B$ has the form:
$$\forall((\Pi^1_1\conj\borel))\arrow \borel)$$
Therefore $B$ is $\bpi12$. \qed

Harrington \site{harring} also shows how to choose
$B$ so that the generic extension
has a $\bdel13$ well-ordering of $\bairespace$.
He also shows how to take a further innocuous extensions
to make $B$ a $\Delta^1_3$  set and to get a $\Delta^1_3$  well-ordering.

\newpage\part{Classical Separation Theorems}
\section{Souslin-Luzin Separation Theorem}

Define $A\subseteq\bairespace$ to be $\kappa$-Souslin iff
\sdex{$\kappa$-Souslin}
\sdex{Souslin-Luzin Separation}
there exists a tree $T\subseteq \Union_{n<\omega}(\kappa^n\cross\omega^n)$
such that
$$y\in A\rmiff \exists x\in\kappa^\omega\;\forall n<\omega\;\;
(x\res n, y\res n)\in T.$$
In this case we write $A=p[T]$, the projection of the infinite branches
of the tree $T$.
Note that $\omega$-Souslin is the same as $\bsig11$.

Define the $\kappa$-Borel sets to be the smallest family of subsets
of $\bairespace$ containing \sdex{$\kappa$-Borel}
the usual Borel sets and closed under intersections or unions of size
$\kappa$ and complements.

\begin{theorem} \label{sep}
  Suppose $A$ and $B$ are disjoint $\kappa$-Souslin subsets of
  $\bairespace$.  Then there exists a $\kappa$-Borel set $C$ which separates
  $A$ and $B$, i.e., $A\subseteq C$ and $C\intersect B=\emptyset$.
\end{theorem}
\proof
Let $A=p[T_A]$ and $B=p[T_B]$.  Given a tree
$T\subseteq \Union_{n<\omega}(\kappa^n\cross\omega^n)$, and
$s\in \kappa^{<\omega}$, $t\in\omega^{<\omega}$ (possibly of
different lengths), define
$$T^{s,t}=\{(\hat{s},\hat{t})\in T:
(s\subseteq \hat{s}\rmor \hat{s}\subseteq {s})\rmand
(t\subseteq \hat{t}\rmor \hat{t}\subseteq {t})\}.$$

\begin{lemma}
  Suppose $p[T_A^{s,t}]$ cannot be separated from $p[T_B^{r,t}]$
  by a $\kappa$-Borel set.  Then for some $\alpha<\kappa$ the set

  $p[T_A^{s\concat\alpha,t}]$ cannot be separated from $p[T_B^{r,t}]$
  by a $\kappa$-Borel set.
\end{lemma}
\proof
Note that $p[T_A^{s,t}]=\Union_{\alpha<\kappa}p[T_A^{s\concat\alpha,t}]$.
If there were no such $\alpha$, then for every $\alpha$ we would have a
$\kappa$-Borel set $C_\alpha$ with
$$p[T_A^{s\concat\alpha,t}]\subseteq C_\alpha\rmand
C_\alpha\intersect p[T_B^{r,t}]=\emptyset.$$
But then $\Union_{\alpha<\kappa}C_\alpha$ is a $\kappa$-Borel
set separating $p[T_A^{s,t}]$ and $p[T_B^{r,t}]$.\qed

\begin{lemma}
  Suppose $p[T_A^{s,t}]$ cannot be separated from $p[T_B^{r,t}]$
  by a $\kappa$-Borel set.  Then for some $\beta<\kappa$

  $p[T_A^{s,t}]$ cannot be separated from $p[T_B^{r\concat\beta,t}]$
  by a $\kappa$-Borel set.
\end{lemma}
\proof
Since $p[T_B^{r,t}]=\Union_{\beta<\kappa}p[T_B^{r\concat\beta,t}]$,
if there were no such $\beta$ then for every $\beta$ we would have
$\kappa$-Borel set $C_\beta$ with
$$p[T_A^{s,t}]\subseteq C_\beta\rmand
C_\beta\intersect p[T_B^{r\concat\beta,t}]=\emptyset.$$
But then $\Intersect_{\beta<\kappa}C_\beta$ is a $\kappa$-Borel
set separating $p[T_A^{s,t}]$ and $p[T_B^{r,t}]$.\qed

\begin{lemma}
  Suppose $p[T_A^{s,t}]$ cannot be separated from $p[T_B^{r,t}]$
  by a $\kappa$-Borel set.  Then for some $n<\omega$

  $p[T_A^{s,t\concat n}]$ cannot be separated from $p[T_B^{r,t\concat n}]$
  by a $\kappa$-Borel set.
\end{lemma}
\proof
Note that
$$p[T_A^{s,t\concat n}]=p[T_A^{s,t}]\intersect [t\concat n]$$
 and
$$p[T_B^{r,t\concat n}]=p[T_B^{r,t}]\intersect [t\concat n].$$
Thus if $C_n\subseteq [t\concat n]$ were to separate
$p[T_A^{s,t\concat n}]$ and $p[T_B^{r,t\concat n}]$
for each $n$, then $\Union_{n<\omega}C_n$ would
separate $p[T_A^{s,t}]$ from  $p[T_B^{r,t}]$.\qed

To prove the separation theorem apply the lemmas iteratively in rotation to
obtain, $u,v\in\kappa^\omega$ and $x\in\bairespace$ so that
for every $n$, $p[T_A^{u\res n,x\res n}]$ cannot be separated from
$p[T_B^{v\res n,x\res n}]$.  But necessarily, for every $n$
$$(u\res n,x\res n)\in T_A \rmand (v\res n,x\res n)\in T_B$$
otherwise either  $p[T_A^{u\res n,x\res n}]=\emptyset$ or
$p[T_B^{v\res n,x\res n}]=\emptyset$ and they could be separated.
But this means that
$x\in p[T_A]=A$ and $x\in p[T_B]=B$ contradicting the fact that
they are disjoint.\qed

\sector{Kleene Separation Theorem}

We  begin by defining the hyperarithmetic subsets of $\bairespace$.
\sdex{hyperarithmetic sets}
We continue with our view of Borel sets as well-founded trees with little
dohickey's (basic clopen sets) attached to its terminal nodes.
\sdex{Kleene Separation Theorem}

A \dex{code for a hyperarithmetic set}
is a triple  $(T,p,q)$ where $T$ is a
recursive well-founded subtree of $\finseq$, $p:\tone\arrow 2$ is recursive,
 and $q:\tzero\arrow \base$
is a recursive map, where $\base$ is the set
of basic clopen subsets of $\bairespace$ including the empty set.
Given a code $(T,p,q)$ we define
$\langle C_s:s\in T\rangle$ as follows.

\begin{itemize}
  \item if $s$ is a terminal node of $T$, then
          $$C_s=q(s)$$
  \item if $s$ is a not a terminal node and $p(s)=0$, then
   $$C_s=\Union\{C_{s\concat n}:{s\concat n}\in T\},$$
   and
  \item if $s$ is a not a terminal node and $p(s)=1$, then
   $$C_s=\Intersect\{C_{s\concat n}:{s\concat n}\in T\}.$$
\end{itemize}
Here we are being a little more flexible by allowing unions and
intersections at various nodes.

Finally, the set $C$ coded by $(T,p,q)$ is the set $C_{\emptyseq}$.
A set $C\subseteq\bairespace$ is hyperarithmetic iff
it is coded by some recursive $(T,p,q)$.

\begin{theorem} (Kleene \site{kleene})
  Suppose $A$ and $B$ are disjoint $\Sigma^1_1$ subsets of
  of $\bairespace$.  Then there exists a
  hyperarithmetic set $C$ which separates them, i.e.,
  $A\subseteq C$ and $C\intersect B=\emptyset$.
\end{theorem}
\proof
This amounts basically to a constructive proof of the
classical Separation Theorem \ref{sep}.

Let $A=p[T_A]$ and $B=p[T_B]$ where $T_A$ and $T_B$ are
recursive subtrees of $\Union_{n\in\omega}(\omega^n\cross\omega^n)$,
and
$$p[T_A]=\{y:\exists x \forall n\;\;\; (x\res n,y\res n)\in T_A\}$$
and similarly for $p[T_B]$. Now
define the tree
$$T=\{(u,v,t): (u,t)\in T_A\rmand (v,t)\in T_B\}.$$
Notice that $T$ is recursive tree which is well-founded.  Any infinite
branch thru $T$ would give a point in the intersection of $A$ and
$B$ which would contradict the fact that they are disjoint.

Let $T^+$ be the tree of all nodes which are either ``in'' or ``just out'' of $T
$,
i.e.,
$(u,v,t)\in T^+$ iff $(u\res n,v\res n,t\res n)\in T$ where $|u|=|v|=|t|=n+1$.
Now we define the family of sets
$$\langle C_{(u,v,t)}:(u,v,t)\in T^+\rangle$$
as follows.

Suppose $(u,v,t)\in T^+$ is a terminal node of $T^+$.  Then since
$(u,v,t)\notin T$ either $(u,t)\notin T_A$ in which case we
define $C_{(u,v,t)}=\emptyset$ or $(u,t)\in T_A$ and $(v,t)\notin T_B$
in which case we define $C_{(u,v,t)}=[t]$.  Note that
in either case $C_{(u,v,t)}\subseteq [t]$ separates
$p[T_A^{u,t}]$ from $p[T_B^{v,t}]$.

\begin{lemma}
  Suppose $\langle A_n:n<\omega\rangle$, $\langle B_m:m<\omega\rangle$
  $\langle C_{nm}:n,m<\omega\rangle$
  are such that for every $n$ and $m$  $C_{nm}$ separates
  $A_n$ from $B_m$.  Then both
  $\Union_{n<\omega}\Intersect_{m<\omega}C_{nm}$ and
  $\Intersect_{m<\omega}\Union_{n<\omega}C_{nm}$ separate
  $\Union_{n<\omega}A_n$ from $\Union_{m<\omega}B_m$.
\end{lemma}
\proof
Left to reader.
\qed
It follows from the Lemma that if we let
$$C_{(u,v,t)}=\Union_{k<\omega}\Intersect_{m<\omega}\Union_{n<\omega}
C_{(u\concat n,v\concat m, t\concat k)}$$
(or any other permutation\footnote{Algebraic symbols are used
when you do not know what you are talking about.-
Philippe Schnoebelen}
of $\Intersect$ and $\Union$), then
by induction on rank of $(u,v,t)$ in $T^+$ that
$C_{(u,v,t)}\subseteq [t]$ separates $p[T_A^{u,t}]$ from $p[T_B^{v,t}]$.
Hence, $C=C_{(\emptyseq,\emptyseq,\emptyseq)}$ separates
$A=p[T_A]$ from $B=p[T_B]$.

To get a
hyperarithmetic code use the tree consisting of
all subsequences of sequences of the form,
$$\langle t(0),v(0),u(0),\ldots,t(n),v(n),u(n)\rangle$$ where
$(u,v,t)\in T^+$.  Details are left to the reader.
\qed

The theorem also holds for
$A$ and $B$ disjoint $\Sigma^1_1$ subsets
of $\omega$.  One way to see this is to identify $\omega$
with the constant functions in $\bairespace$.  The definition
of hyperarithmetic code $(T,p,q)$ is changed only by
letting $q$ map into the finite subsets of $\omega$.

\begin{theorem}\label{hypcode}
  If $C$ is a hyperarithmetic set, then $C$ is $\Delta_1^1$.
\end{theorem}
\proof
This is true whether $C$ is a subset of $\bairespace$ or $\omega$.
We just do the case $C\subseteq\bairespace$.
Let $(T,p,q)$ be a hyperarithmetic code for $C$.  Then
$x\in C$ iff

\noindent there exists a function $in:T\arrow \{0,1\}$ such that
\begin{enumerate}
  \item if $s$ a terminal node of $T$, then $in(s)=1$ iff $x\in q(s)$,
  \item if $s\in T$ and not terminal and $p(s)=0$, then
  $in(s)=1$ iff there exists $ n $ with $s^\concat n\in T$ and
  $in(s^\concat n)=1$,
  \item if $s\in T$ and not terminal and $p(s)=1$, then

  $in(s)=1$ iff for all $ n $ with $s^\concat n\in T$ we have $in(s^\concat n)=1
$,
  and finally,
  \item $in(\emptyseq)=1$.
\end{enumerate}
Note that (1) thru (4) are all $\Delta_1^1$ (being a terminal
node in a recursive tree is $\Pi^0_1$, etc).
It is clear
that $in$ is just coding up whether or not $x\in C_s$ for $s\in T$.
Consequently, $C$ is $\Sigma^1_1$.  To see that $\comp{C}$ is $\Sigma^1_1$
note
that $x\notin C$ iff

\noindent there exists $in:T\arrow \{0,1\}$ such that
(1), (2), (3), and (4)$^\prime$ where

$4^\prime \;\;\;\;in(\emptyseq)=0.$

\qed

\begin{corollary}
  A set is $\Delta^1_1$ iff it is hyperarithmetic.
\end{corollary}

\begin{corollary} \label{separation}
  If $A$ and $B$ are disjoint $\Sigma^1_1$ sets, then
  there exists a $\Delta^1_1$ set which separates them.
\end{corollary}

For more on the effective Borel hierarchy, see Hinman \site{hinman}.
See Barwise \site{barwise} for a model theoretic or admissible
sets approach to the hyperarithmetic hierarchy.

\sector{$\Pi^1_1$-Reduction}

We say that $A_0$,$B_0$ reduce $A$,$B$ iff
\begin{enumerate}
  \item $A_0\subseteq A$ and $B_0\subseteq B$,
  \item $A_0\union B_0=A \union B$, and
  \item $A_0\intersect B_0=\emptyset$.
\end{enumerate}

$\Pi^1_1$-reduction is the property that every pair of
\sdex{$\Pi^1_1$-Reduction}
$\Pi_1^1$ sets  can
be reduced by a pair of $\Pi^1_1$ sets.
The sets can be  either subsets of $\omega$ or of $\bairespace$.

\begin{theorem}
  $\Pi^1_1$-uniformity implies  $\Pi^1_1$-reduction.
\end{theorem}
\proof
Suppose $A,B\subseteq X$ are $\Pi^1_1$ where $X=\omega$ or
$X=\bairespace$.  Let $P=(A\cross \{0\})\union (B\cross\{1\})$.
Then $P$ is a $\Pi^1_1$ subset of $X\cross\bairespace$ and so
by $\Pi_1^1$-uniformity (Theorem \ref{uniform}) there exists
$Q\subseteq P$ which is $\Pi_1^1$ and for every $x\in X$,
if there exists $i\in\{0,1\}$ such that $(x,i)\in P$, then there
exists a unique $i\in\{0,1\}$ such that $(x,i)\in Q$.
Hence, letting
$$A_0=\{x\in X: (x,0)\in Q\}$$
and
$$B_0=\{x\in X: (x,1)\in Q\}$$
gives a pair of $\Pi^1_1$ sets which reduce $A$ and $B$.\qed

There is also a proof of reduction using the \sdex{prewellordering}
prewellordering
property, which is a weakening of the scale property used
in the proof of $\Pi_1^1$-uniformity.  So, for example, suppose
$A$ and $B$ are are $\Pi^1_1$ subsets of $\bairespace$.  Then
we know there are maps from $\bairespace$ to trees,
$$x\mapsto T^a_x \rmand y\mapsto T^b_y$$
which are ``recursive'' and

$x\in A$ iff $T_x^a$ is well-founded and

$y\in B$ iff $T_y^b$ is well-founded.

\noindent Now define
\begin{enumerate}
  \item $x\in A_0$ iff $x\in A$ and not ($T_x^b\strictembed T_x^a$), and
  \item $x\in B_0$ iff $x\in B$ and not ($T_x^a\embed T_x^b$).
\end{enumerate}
Since $\strictembed$ and $\embed$ are both $\Sigma^1_1$ it is clear,
that $A_0$ and $B_0$ are $\Pi_1^1$ subsets of $A$ and $B$ respectively.
If $x\in A$ and $x\notin B$, then $T_x^a$ is well-founded and
$T_x^b$ is ill-founded and so not ($T_x^b\strictembed T_x^a$) and
$a\in A_0$.  Similarly, if $x\in B$ and $x\notin A$, then $x\in B_0$.
If $x\in A\intersect B$, then both $T_x^b$ and $T_x^a$ are well-founded
and either $T_x^a\embed T_x^b$, in which case $x\in A_0$ and $x\notin B_0$,
or $T_x^b\strictembed T_x^a$, in which case $x\in B_0$ and $x\notin A_0$.

\begin{theorem} \label{separ}
  $\Pi_1^1$-reduction implies $\Sigma^1_1$-separation, i.e., for
  any two disjoint $\Sigma^1_1$ sets $A$ and $B$ there exists a
  $\Delta^1_1$-set $C$ which separates them. i.e., $A\subseteq C$ and
  $C\intersect B=\emptyset$.
\end{theorem}
\proof
Note that $\comp{A}\;\;\union\comp{B}=X$.  If $A_0$ and $B_0$ are
$\Pi_1^1$ sets reducing  $\comp{A}$ and $\comp{B}$, then
$\comp{A_0}=B_0$, so they are both $\Delta^1_1$. If we
set $C=B_0$ , then $C=B_0=\comp{A_0}\subseteq \comp{A}$ so
$C\subseteq \comp{A}$ and therefore $A\subseteq C$.  On the
other hand $C=B_0\subseteq\comp{B}$ implies $C\intersect B=\emptyset$.
\qed

\sector{$\Delta^1_1$-codes}

Using $\Pi_1^1$-reduction and universal sets it is possible to
get codes for $\Delta^1_1$ subsets of $\omega$ and $\bairespace$.
\sdex{$\Delta^1_1$-codes}

Here is what we mean by $\Delta_1^1$ codes for subsets of $X$
where $X=\omega$ or $X=\bairespace$.

There exists a $\Pi^1_1$ sets $C\subseteq\omega\cross\bairespace$ and
$P\subseteq \omega\cross\bairespace\cross X$ and a
$\Sigma^1_1$ set $S\subseteq\omega\cross\bairespace\cross X$ such that
\begin{itemize}
  \item for any $(e,u)\in C$
   $$\{x\in X: (e,u,x)\in P\}=\{x\in X:(e,u,x)\in S\}$$
   \item for any $u\in\bairespace$ and $\Delta^1_1(u)$ set
   $D\subseteq X$ there exists a $(e,u)\in  C$ such that
   $$D=\{x\in X: (e,u,x)\in P\}=\{x\in X:(e,u,x)\in S\}.$$
\end{itemize}

From now on we will write

``$e$ is a $\Delta^1_1(u)$-code for a subset of $X$''

to mean $(e,u)\in C$ and  remember that it is a $\Pi^1_1$ predicate.

We also write ``$D$ is the $\Delta^1_1(u)$ set coded by $e$'' if
``$e$ is a $\Delta^1_1(u)$-code for a subset of $X$'' and
$$D=\{x\in X: (e,x)\in P\}=\{x\in X:(e,x)\in S\}.$$
Note that
$x\in D$ can be said in either a $\Sigma^1_1(u)$ way or $\Pi^1_1(u)$ way,
using either $S$ or $P$.

\begin{theorem} (Spector-Gandy \site{spector}\site{gandy}) \label{deltacode}
  $\Pi^1_1$-reduction and universal sets implies $\Delta_1^1$ codes
  exist.
\end{theorem} \sdex{Spector-Gandy Theorem}
\proof
Let $U\subseteq \omega\cross\bairespace\cross X$ be a $\Pi^1_1$ set
which is universal for all $\Pi^1_1(u)$ sets, i.e., for every
$u\in\bairespace$ and $A\in\Pi^1_1(u)$ with $A\subseteq X$ there
exists $e\in\omega$ such that $A=\{x\in X: (e,u,x)\in U\}$.
For example, to get such a $U$ proceed as follows.  Let
$\{e\}^u$ be the partial function you get by using the $e^{th}$
Turing machine with oracle $u$.  Then define
$(e,u,x)\in U$ iff
$\{e\}^u$ is the characteristic function
of a tree $T\subseteq \Union_{n<\omega}(\omega^n\cross\omega^n)$ and
$T_x =\{s: (s,x\res |s|)\in T\}$ is well-founded.

Now get a doubly universal pair.
Let $e\mapsto (e_0,e_1)$ be the usual recursive unpairing function
from $\omega$ to $\omega\cross\omega$ and define
$$U^0=\{(e,u,x):(e_0,u,x)\in U\}$$
and
$$U^1=\{(e,u,x):(e_1,u,x)\in U\}.$$
The pair of sets $U^0$ and $U^1$ are $\Pi^1_1$ and doubly universal,
i.e., for any $u\in \bairespace$ and $A$ and $B$ which are
$\Pi^1_1(u)$ subsets
of $X$ there exists $e\in\omega$ such that
$$A=\{x: (e,u,x)\in U^0\}$$
and
$$B=\{x: (e,u,x)\in U^1\}.$$
Now apply reduction to obtain $P^0\subseteq U^0$ and $P^1\subseteq U^1$
which are $\Pi^1_1$ sets.  Note that the by the nature of taking
cross sections, $P^0_{e,u}$ and $P^1_{e,u}$  reduce
 $U^0_{e,u}$ and $U^1_{e,u}$.  Now we define
\begin{itemize}
\item ``$e$ is a $\Delta^1_1(u)$ code'' iff
$\forall x\in X (x\in P^0_{e,u}$ or $x\in P^1_{e,u})$, and
\item $P=P^0$ and $S=\comp{P^1}$.
\end{itemize}

Note that $e$ is a $\Delta^1_1(u)$ code is a $\Pi^1_1$ statement in $(e,u)$.
Also if $e$ is a $\Delta^1_1(u)$ code, then
$P_{(e,u)}=S_{e,u}$ and so its a  $\Delta^1_1(u)$ set.   Furthermore
if $D\subseteq X$ is a $\Delta_1^1(u)$ set, then since $U^0$ and
$U^1$ were a doubly universal pair, there exists $e$ such that
$U^0_{e,u}=D$ and $U^1_{e,u}=\comp{D}$.  For this $e$ it must
be that $U^0_{e,u}=P^0_{e,u}$ and $U^1_{e,u}=P^1_{e,u}$ since
the $P$'s reduce the $U$'s.  So this $e$ is a $\Delta_1^1(u)$ code
which codes the set $D$.\qed

\begin{corollary} \label{gandy}
  $\{(x,u)\in P(\omega)\cross \bairespace: x\in\Delta^1_1(u)\}$ is $\Pi^1_1$.
\end{corollary}
\proof
$x\in\Delta^1_1(u)$ iff $\exists e\in\omega$ such that
\begin{enumerate}
\item $e$ is a $\Delta^1_1(u)$ code,
\item $\forall n$ if $n\in x$, then $n$ is in the $\Delta^1_1(u)$-set
coded by $e$, and
\item $\forall n$ if $n$ is the $\Delta^1_1(u)$-set
coded by $e$, then $n\in x$.
\end{enumerate}
Note that clause (1) is $\Pi^1_1$.  Clause (2) is $\Pi^1_1$ if
we use that  $(e,u,n)\in P$ is equivalent to ``$n$ is in
the $\Delta^1_1(u)$-set
coded by $e$''.  While clause (3) is is $\Pi^1_1$ if
we use that  $(e,u,n)\in S$ is equivalent to
``$n$ is in the $\Delta^1_1(u)$-set
coded by $e$''. \qed

We say that $y\in\bairespace$ is $\Delta^1_1(u)$ iff
its graph $\{(n,m):y(n)=m\}$ is $\Delta^1_1(u)$.  Since being the
graph a function is a $\Pi^0_2$ property it is easy to see how
to obtain $\Delta^1_1(u)$ codes for functions $y\in\bairespace$.

\begin{corollary}\label{gandy2}
  Suppose $\theta(x,y,z)$ is a $\Pi^1_1$ formula, then
  $$\psi(y,z)=\exists x\in\Delta^1_1(y)\; \theta(x,y,z)$$
  is a $\Pi^1_1$ formula.
\end{corollary}
\proof
$\psi(y,z)$ iff

\noindent $\exists e\in\omega$ such that
\begin{enumerate}
  \item $e$ is a $\Delta^1_1(y)$ code, and
  \item $\forall x$ if $x$ is the set coded by $(e,y)$, then $\theta(x,y,z)$.
\end{enumerate}
This will be $\Pi^1_1$ just in case the clause
``$x$ is the set coded by $(e,y)$'' is $\Sigma^1_1$.  But this is
$\Delta^1_1$ provided that $e$ is a $\Delta^1_1(y)$ code, e.g.,
for $x\subseteq\omega$ we just say:
$\forall n\in\omega$
\begin{enumerate}
  \item if $n\in x$ then $(e,y,n)\in S$ and
  \item if $(e,y,n)\in P$, then $n\in x$.
\end{enumerate}
Both of these clauses are $\Sigma^1_1$ since $S$ is $\Sigma^1_1$ and
$P$ is $\Pi^1_1$.  A similar argument works for $x\in\bairespace$.\qed

The method of this corollary also works for the quantifier
$$\exists D\subseteq\bairespace\mbox{ such that }
D\in \Delta(y)\;\;\theta(D,y,z).$$
It is equivalent to say $\exists e\in\omega$ such that
 $e$ is a $\Delta^1_1(y)$ code for a subset of $\bairespace$ and
 $\theta(\ldots,y,z)$ where occurrences of the
 ``$q\in D$'' in the formula $\theta$ have been
 replaced by either $(e,y,q)\in P$ or
$(e,y,q)\in S$, whichever is necessary to makes $\theta$ come out $\Pi^1_1$.

\begin{corollary}
Suppose $f:\bairespace\arrow\bairespace$ is Borel, $B\subseteq\bairespace$
is Borel, and $f$ is one-to-one on $B$.  Then the image of $B$ under
$f$,  $f(B)$, is Borel.
\end{corollary}

\proof
By relativizing the following argument to an arbitrary parameter
we may assume that the graph of $f$ and the set $B$ are
$\Delta^1_1$.   Define
$$R=\{(x,y): f(x)=y\rmand x\in B\}.$$
Then for any $y$ the set
$$\{x: R(x,y)\}$$
is a $\Delta^1_1(y)$ singleton (or empty).   Consequently,
its unique element is $\Delta^1_1$ in $y$.
It follows that
$$y\in f(B) \rmiff \exists x\; R(x,y)\rmiff
\exists x\in\Delta^1_1(y)\; R(x,y)$$
and so $f(B)$ is both $\Sigma^1_1$ and $\Pi^1_1$.
\qed

Many applications of the Gandy-Spector Theorem exist.
For example, it is shown (assuming V=L in all three cases)
that
\begin{enumerate}
  \item there exists an uncountable
   $\Pi^1_1$ set which is concentrated on the rationals
   (Erdos, Kunen, and Mauldin \site{erdkm}),
  \item there
    exists a $\Pi^1_1$ Hamel basis
   (Miller \site{infcomb}), and
  \item there exists
   a topologically rigid $\Pi^1_1$ set
    (Van Engelen,
    Miller, and Steel \site{ems}).
\end{enumerate}

\addtocontents{toc}{\protect\newpage}
\newpage\part{Gandy Forcing}
\section{$\Pi^1_1$ equivalence relations}

\begin{theorem} \label{piequiv}
  (Silver \site{silver}) Suppose $(X,E)$ is  a $\bpi11$ equivalence
  relation, i.e. $X$ is a Borel set and $E\subseteq X^2$ is
  a $\bpi11$ equivalence relation on $X$.
  Then either $E$ has \sdex{$\Pi^1_1$ equivalence relations}
  countably many equivalence classes or there exists a
  perfect set of pairwise inequivalent elements.
\end{theorem}

Before giving the proof consider the following example.
Let \sdex{$WO$} $WO$ be the set of all characteristic
functions of well-orderings
of $\omega$.  This is a $\Pi^1_1$ subset of $2^{\omega\cross\omega}$.
Now define $x\isom y$ iff there exists an isomorphism taking
$x$ to $y$ or $x,y\notin WO$.  Note that
$(2^{\omega\cross\omega},\isom)$ is a $\Sigma_1^1$ equivalence relation
with exactly $\omega_1$ equivalence classes.  Furthermore,
if we restrict $\isom$ to $WO$, then $(WO,\isom)$ is
a $\Pi^1_1$ equivalence relation (since well-orderings
are isomorphic iff neither is isomorphic to an initial segment
of the other).   Consequently, Silver's theorem is the best
possible.

The proof we are going to give is due to Harrington \site{harrsilv}, see also
Kechris and Martin \site{kecmar}, Mansfield and Weitkamp
\site{manweit} and Louveau \site{louvsilv}.  A model theoretic proof is given in
Harrington and Shelah \site{harrshel}.

We can assume that $X$ is $\Delta^1_1$ and $E$ is $\Pi^1_1$, since
the proof readily relativizes to an arbitrary parameter.
Also, without loss, we may assume that $X=\bairespace$ since
we just make the complement of $X$ into one more equivalence class.

Let $\poset$ be the partial order of nonempty $\Sigma^1_1$ subsets
of $\bairespace$ ordered by inclusion.  This is known as
\dex{Gandy forcing}.

\begin{lemma}\label{real}
  If $G$ is $\poset$-generic over $V$, then there exists
  $a\in\bairespace$ such that $G=\{p\in \poset: a\in p\}$ and
  $\{a\}=\Intersect G$.
\end{lemma}

\proof
For every $n$ an easy density argument shows that there exists
a unique $s\in \omega^n$ such that $[s]\in G$ where
$[s]=\{x\in\bairespace: s\subseteq x\}$.  Define
$a\in \bairespace$ by $[a\res n]\in G$ for each $n$.
Clearly, $\Intersect G\subseteq \{a\}$.

Now suppose $B\in G$,  we need to show $a\in B$.
Let $B=p[T]$.

\claim There exists $x\in \bairespace$ such that $p[T^{x\res n,a\res n}]\in G$
for every $n\in\omega$.

\proof
This is by induction on $n$.  Suppose
$p[T^{x\res n,a\res n}]\in G$.  Then
$$p[T^{x\res n,a\res n+1}]\in G$$
 since
$$p[T^{x\res n,a\res n+1}]=[a\res n+1]\intersect p[T^{x\res n,a\res n}]$$
and both of these are in $G$.
But note that
$$p[T^{x\res n,a\res n+1}]=
\Union_{m\in\omega}p[T^{x\res n\concat m,a\res n+1}]$$
and so by a density argument there exists $m=x(n)$ such that
$$p[T^{x\res n\concat m,a\res n+1}]\in G.$$
This proves the Claim. \qed
 By the Claim we have that $(x,a)\in [T]$ (since elements of $\poset$
are nonempty) and so $a\in p[T]=B$. Consequently,
$\Intersect G=\{a\}$.  Now suppose that $a\in p\in\poset$ and
$p\notin G$.  Then since
$$\{q\in\poset: q\leq p \rmor q\intersect p=\emptyset \}$$
is dense there must be $q\in G$ with $q\intersect p=\emptyset$.
But this is impossible, because $a\in q\intersect p$, but
$q\intersect p=\emptyset$ is a $\Pi^1_1$ sentence and hence absolute.
\qed

We say that $a\in\bairespace$ is $\poset$-generic over $V$ iff
$G=\{p\in\poset: a\in p\}$ is $\poset$-generic over $V$.

\begin{lemma}\label{pair}
  If $a$ is $\poset$-generic over $V$ and $a=\langle a_0,a_1\rangle$
  (where $\langle,\rangle$ is the standard pairing function), then
  $a_0$ and $a_1$ are both $\poset$-generic over $V$.
\end{lemma}
\proof
The proof is symmetric so we just do it for $a_0$.  Note that we are
not claiming that they are product generic only that each is
separately generic.
Suppose $D\subseteq \poset$ is dense open.  Let
$$E=\{p\in\poset: \{x_0:x\in p\}\in D\}.$$
To see that $E$ is dense let $q\in\poset$ be arbitrary. Define
$$q_0=\{x_0:x\in q\}.$$
Since $q_0$ is a nonempty $\Sigma^1_1$ set and $D$ is dense, there
exists $r_0\leq q_0$ with $r_0\in D$.   Let
$$r=\{x\in q: x_0\in r_0\}.$$
Then $r\in E$ and $r\leq q$.

Since $E$ is dense we have that there exists $p\in E$ with
$a\in p$ and consequently,
$$a_0\in p_0=\{x_0:x\in p\}\in D.$$
\qed

\begin{lemma} \label{spill}
  Suppose $B\subseteq \bairespace$ is $\Pi^1_1$ and for every $x,y\in B$
  we have that $xEy$.  Then there exists a $\Delta^1_1$ set $D$ with
  $B\subseteq D\subseteq \bairespace$ and such that for every
  $x,y\in D$ we have that $xEy$.
\end{lemma}
\proof
Let $A=\{x\in \bairespace:\forall y \;\; y\in B \arrow x E y\}$.  Then
$A$ is a $\Pi^1_1$ set which contains the $\Sigma^1_1$ set $B$, consequently
by the Separation Theorem \ref{separation} or \ref{separ} there
exists a $\Delta^1_1$ set $D$ with
$B\subseteq D \subseteq A$. Since all elements of $B$ are equivalent,
so are all elements of $A$ and hence  $D$ is as required.
\qed

Now we come to the heart of Harrington's proof.  Let
$B$ be the union of all $\Delta^1_1$
subsets of $\bairespace$ which meet only one equivalence class of $E$, i.e.
$$B=\Union\{D\subseteq \bairespace: D\in\Delta^1_1 \rmand \forall x,y\in D
\;\;xEy\}.$$
Since $E$ is $\Pi^1_1$
we know that by using $\Delta_1^1$ codes that this union is $\Pi^1_1$, i.e.,

$z\in B$ iff $\exists e\in\omega$ such that
\begin{enumerate}
  \item $e$ is a $\Delta^1_1$ code for
   a subset of $\bairespace$,
  \item $\forall x,y$ in the set coded by $e$ we have $x E y$, and
  \item $z$ is in the set coded by $e$.
\end{enumerate}
Note that item (1) is $\Pi^1_1$ and (2) and (3) are both $\Delta^1_1$
(see Theorem \ref{deltacode}).

\bigskip

If $B=\bairespace$, then since there are only countably many
$\Delta^1_1$ sets, there would only be countably many $E$ equivalence
classes and we are done.  So assume $A=\comp{B}$ is a nonempty
$\Sigma^1_1$ set and in this case we will prove that there is a perfect
set of $E$-inequivalent reals.

\begin{lemma} \label{crux}
  Suppose $c\in\bairespace\intersect V$. Then
  $$A\forces_{\poset}\check{c}\notE \name{a}$$
  where  $\name{a}$ is a name for the generic real (Lemma \ref{real}).
\end{lemma}

\proof
Suppose not, and let $C\subseteq A$ be a nonempty $\Sigma^1_1$ set
such that $C\forces cEa$.  We know that there must exists
$c_0,c_1\in C$ with $c_0\notE c_1$.   Otherwise there would exists
a $\Delta^1_1$ superset of $C$ which meets only one equivalence class
(Lemma \ref{spill}).  But we these are all disjoint from $A$.
Let
$$Q=\{c:c_0\in C, c_1\in C, \rmand c_0\notE c_1\}.$$
Note that $Q$ is a nonempty $\Sigma^1_1$ set.  Let $a\in Q$ be
$\poset$-generic over $V$.  Then by Lemma \ref{pair} we have that
both $a_0$ and $a_1$ are $\poset$-generic over $V$ and
$a_0\in C$, $a_1\in C$, and $a_0\notE a_1$.  But
$a_i\in C$ and $C\forces a_i E c$ means that
$$a_0 E c, a_1 E c, \rmand a_0\notE a_1.$$  This contradicts
the fact that $E$ is an equivalence relation.

Note that
``$E$ is an equivalence relation'' is a $\Pi^1_1$ statement hence
it is absolute.  Note also that we don't need to assume that
there are $a$ which are $\poset$-generic over $V$.  To see this
replace
$V$ by a countable transitive model $M$ of ZFC$^*$ (a sufficiently large
fragment of ZFC) and use absoluteness. \qed

Note that the lemma implies that if $(a_0,a_1)$ is
$\poset\cross\poset$-generic over $V$ and $a_1\in A$, then
$a_0\notE a_1$.  This is because
$a_1$ is $\poset$-generic over $V[a_0]$ and so $a_0$ can
be regarded as an element of the ground model.

\begin{lemma} \label{perfect}
  Suppose $M$ is a countable transitive model of ZFC$^*$ and
  $\poset$ is a partially ordered set in $M$.  Then there
  exists a ``perfect'' set of $\poset$-filters $\{G_x: x\in\cantorset\}$
  such that for every $x\not=y$ we have that $(G_x,G_y)$ is
  $\poset\cross\poset$-generic over $M$.
\end{lemma}
\proof
Let $D_n$ for $n<\omega$ list all
dense open subsets of $\poset\cross\poset$ which
are in $M$. Construct $\langle p_s: s\in 2^{<\omega}\rangle$ by induction
on the length of $s$ so that
\begin{enumerate}
  \item $s\subseteq t$ implies $p_t\leq p_s$ and
  \item if $|s|=|t|=n+1$ and $s$ and $t$ are distinct, then
   $(p_s,p_t)\in D_n$.
\end{enumerate}
Now define for any $x\in\cantorset$
$$G_x=\{p\in\poset: \exists n \;\; p_{x\res n}\leq p\}.$$
\qed

Finally to prove Theorem \ref{piequiv} let $M$ be a countable transitive
set isomorphic to an elementary substructure of $(V_\kappa,\in)$ for
some sufficiently large $\kappa$.
Let $\{G_x:x\in\cantorset\}$ be
given by Lemma \ref{perfect} with
$A\in G_x$ for all $x$ and let
$$P=\{a_x: x\in\cantorset\}$$
be the corresponding
generic reals. By Lemma \ref{crux} we know that for every
$x\not= y\in\cantorset$ we have that $a_x\notE a_y$.
Note also that $P$ is perfect because the map $x\mapsto a_x$ is
continuous.   This is because for any $n\in\omega$ there exists
$m<\omega$ such that every $p_s$ with $s\in 2^m$ decides $a\res n$.
\qed

\begin{corollary}
Every  $\Sigma^1_1$ set which contains a real
which is not $\Delta^1_1$ contains a perfect subset.
\end{corollary}
\proof
Let $A\subseteq\bairespace$ be a $\Sigma^1_1$ set.
Define $x E y $ iff $x,y\notin A$ or x=y.
Then is $E$ is a $\Pi^1_1$ equivalence relation.
A $\Delta^1_1$ singleton is a $\Delta^1_1$ real,
hence Harrington's set $B$ in the above proof must be nonempty.
Any perfect set of $E$-inequivalent elements
can contain at most one element of $\comp{A}$.
\qed

\begin{corollary}
Every uncountable $\bsig11$ set
contains a perfect subset.
\end{corollary}

Perhaps this is not such a farfetched way of proving this
result, since one of the usual proofs looks like a combination of
Lemma \ref{real} and \ref{perfect}.

\sector{Borel metric spaces and lines in the plane}

We give two applications of Harrington's method of proof of Silver's
$\Pi^1_1$-equivalence relation theorem.    First let us begin
by isolating a principal which we call overflow.  It is an easy consequence
of the Separation Theorem.

\begin{lemma} \label{overflow}
  (Overflow) Suppose $\theta(x_1,x_2,\ldots,x_n)$ is a $\Pi^1_1$ formula
  and $A$ is a  $\Sigma^1_1$ set such that
  $$\forall x_1,\ldots, x_n\in A\;\;
     \theta(x_1,\ldots,x_n).$$
  Then there exists a $\Delta^1_1$ set $D\supseteq A$ such that
$$\forall x_1,\ldots, x_n\in D\;\;
     \theta(x_1,\ldots,x_n).$$
\end{lemma}
\proof
For $n=1$ this is just the Separation Theorem \ref{separation}.

For $n=2$ define
$$B=\{x:\forall y (y\in A\arrow \theta(x,y))\}.$$
Then $B$ is $\Pi^1_1$ set which contains $A$.  Hence by separation
there exists a $\Delta^1_1$ set $E$ with $A\subseteq E\subseteq B$.
Now define
$$C=\{x:\forall y (y\in E\arrow \theta(x,y)\}.$$
Then $C$ is a $\Pi^1_1$ set which also contains $A$.  By applying
separation again we get a $\Delta^1_1$ set $F$ with
$A\subseteq F\subseteq C$.  Letting $D=E\intersect F$ does the
job.  The proof for $n>2$ is similar.
\qed

We say that $(B,\delta)$ is a \dex{Borel metric space} iff $B$ is Borel,
$\delta$ is a metric on $B$, and
and for every $\epsilon\in\rationals$ the set
$$\{(x,y)\in B^2:\delta(x,y)\leq\epsilon\}$$
is Borel.

\begin{theorem}\label{bormet}
  (Harrington \site{harshemar})  If $(B,\delta)$ is a Borel metric space,
  then either $(B,\delta)$ is separable (i.e., contains a countable
  dense set) or for some $\epsilon>0$ there exists a perfect set
  $P\subseteq B$ such that $\delta(x,y)>\epsilon$ for every
  distinct $x,y\in P$.
\end{theorem}
\proof

By relativizing the proof to an arbitrary parameter we may assume that
$B$ and the sets $\{(x,y)\in B^2:\delta(x,y)\leq\epsilon\}$ are
$\Delta^1_1$.

\begin{lemma} \label{overmetric}
For any $\epsilon\in\rationals^+$ if $A\subseteq B$ is $\Sigma^1_1$ and
the diameter of $A$ is less than $\epsilon$, then there exists a
$\Delta^1_1$ set $D$ with diameter less than $\epsilon$ and
$A\subseteq D\subseteq B$.
\end{lemma}
\proof
This follows from Lemma \ref{overflow}, since
$$\theta(x,y)\;\; \rmiff\;\; \delta(x,y)<\epsilon \rmand x,y\in B$$
is a $\Pi^1_1$ formula.\qed

For any $\epsilon\in\rationals^+$ look at
$$Q_\epsilon=\Union\{D\in\Delta^1_1: D\subseteq B \rmand
 \diam(D)<\epsilon\}.$$
Note that $Q$ is a $\Pi_1^1$ set.  If for every $\epsilon\in\rationals^+$
$Q_\epsilon=B$, then since there are only countably many $\Delta^1_1$ sets,
 $(B,\delta)$ is separable and we are done. On the other hand
suppose for some $\epsilon\in\rationals^+$ we have that
$$P_\epsilon=B\setminus Q_\epsilon\not=\emptyset.$$

\begin{lemma}
For every $c\in V\intersect B$
$$P_\epsilon\forces \delta(\name{a},\check{c})>\epsilon/3$$
where $\forces$ is Gandy forcing and $\name{a}$ is a name for
the generic real (see Lemma \ref{real}).
\end{lemma}
\proof
Suppose not.  Then there exists $P\leq P_\epsilon$ such that
$$P\forces \delta(a,c)\leq\epsilon/3.$$
Since $P$ is disjoint from $Q_\epsilon$ by Lemma \ref{overmetric} we
know that the diameter of $P$ is $\geq\epsilon$.  Let
$$R=\{(a_0,a_1): a_0,a_1\in P\rmand \delta(a_0,a_1)>(2/3)\epsilon\}.$$
Then $R$ is in $\poset$ and by Lemma \ref{pair}, if $a$ is $\poset$-generic
over $V$ with $a\in R$, then $a_0$ and $a_1$ are each separately
$\poset$-generic over $V$. But $a_0\in R$ and $a_1\in R$ means
that $\delta(a_0,c)\leq \epsilon/3$ and $\delta(a_1,c)\leq \epsilon/3$.
But by absoluteness $\delta(a_0,a_1)>(2/3)\epsilon$.  This contradicts
the fact that $\delta$ must remain a metric by absoluteness.
\qed

Using this lemma and Lemma \ref{perfect} is now easy to get a perfect
set $P\subseteq B$ such that $\delta(x,y)>\epsilon/3$ for each
distinct $x,y\in P$. This proves Theorem \ref{bormet}.\qed

\begin{theorem}(van Engelen, Kunen, Miller \site{ekm}) \label{collinear}
For any $\bsig11$ set $A$ in the plane, either $A$ can
be covered by countably many lines or there exists a perfect set
$P\subseteq A$ such that no three points of P are collinear.
\end{theorem}\sdex{collinear points}
\proof

This existence of this proof was pointed out to me by
Dougherty, Jackson, and Kechris. The proof in \site{ekm} is
more elementary.

By relativizing the proof
we may as well assume that $A$ is $\Sigma^1_1$.

\begin{lemma} \label{line}
  Suppose $B$ is a $\Sigma^1_1$ set in the plane all of whose
  points are collinear. Then there exists a $\Delta^1_1$ set
  $D$ with $B\subseteq D$ and all points of $D$ are collinear.
\end{lemma}
\proof
This follows from Lemma \ref{overflow} since
$$\theta(x,y,z)\;\; \rmiff \;\; x,y,\rmand z \mbox{ are collinear}$$
is $\Pi^1_1$ (even $\Pi^0_1$).\qed

Define
$$\comp{P}=\Union\{D\subseteq \reals^2: D\in\Delta^1_1 \rmand
\mbox{ all points of $D$ are collinear}\}.$$
It is clear that $\comp{P}$ is $\Pi^1_1$ and therefore $P$ is
$\Sigma^1_1$. If $P\intersect A=\emptyset$, then $A$ can be
covered by countably many lines.

So assume that
$$Q=P\intersect A\not=\emptyset.$$
For any two distinct points
in the plane, $p$ and $q$, let $\lin(p,q)$ be the unique
line on which they lie.

\begin{lemma} \label{line2}
For any two distinct points
in the plane, $p$ and $q$, with $p,q\in V$
$$Q\forces \name{a}\notin \lin(\check{p},\check{q}).$$
\end{lemma}

\proof
Suppose for contradiction that there exists $R\leq Q$
 such that
$$R\forces \name{a}\in \lin(\check{p},\check{q}).$$
Since $R$ is disjoint from
$$\Union\{D\subseteq \reals^2: D\in\Delta^1_1 \rmand
\mbox{ all points of $D$ are collinear}\}$$
it follows from Lemma \ref{line} that not all triples of
points from $R$ are collinear.
Define the nonempty $\Sigma^1_1$ set
$$S=\{a:a_0,a_1,a_2 \in R\rmand a_0,a_1,a_2 \mbox{ are not collinear}\}$$
where $a=(a_0,a_1,a_2)$ via some standard tripling function.
Then $S\in\poset$ and by the obvious generalization of Lemma \ref{pair}
each of the $a_i$ is $\poset$-generic if $a$ is. But this is a
contradiction since all $a_i\in \lin(p,q)$ which makes them collinear.
\qed

The following Lemma is an easy generalization of Lemma \ref{perfect}
so we leave the proof to the reader.

\begin{lemma}\label{line3}
  Suppose $M$ is a countable transitive model of ZFC$^*$ and
  $\poset$ is a partially ordered set in $M$.  Then there
  exists a ``perfect'' set of $\poset$-filters $\{G_x: x\in\cantorset\}$
  such that for every $x,y,z$ distinct, we have that $(G_x,G_y,G_z)$ is
  $\poset\cross\poset\cross\poset$-generic over $M$.
\end{lemma}

Using Lemma \ref{line2} and \ref{line3} it is easy to get (just as in
the proof of Theorem \ref{piequiv}) a perfect set of triply generic
points in the plane, hence no three of which are collinear.  This proves
Theorem \ref{collinear}. \qed

Obvious generalizations of Theorem \ref{collinear} are:
\begin{enumerate}
  \item Any $\bsig11$ subset of $\reals^n$ which cannot be covered
  by countably many lines contains a perfect set all of whose points
  are collinear.

  \item Any $\bsig11$ subset of $\reals^2$ which cannot be covered
  by countably many circles contains a perfect set which does not
  contain four points  on the same circle.

  \item Any $\bsig11$ subset of $\reals^2$ which cannot be covered
  by countably many parabolas contains a perfect set which does not
  contain four points  on the same parabola.

  \item For any $n$ any $\bsig11$ subset of $\reals^2$ which
   cannot be covered by countably many polynomials of degree $<n$
   contains a perfect set which does not
   contain $n+1$ points  on the same polynomial of degree $<n$.

  \item Higher dimensional version of the above involving spheres or
   other surfaces.
\end{enumerate}

A very general statement of this type is due to Solecki \site{Solecki}.
Given any Polish space $X$, family of closed sets $Q$ in $X$, and
analytic $A\subseteq X$; either $A$ can be covered by countably many elements
of $Q$ or there exists a $G_\delta$ set $B\subseteq A$ such that
$B$ cannot be covered by countably many elements of $Q$.  Solecki
deduces Theorem \ref{collinear} from this.

Another result of this type is known as the Borel-Dilworth Theorem.
\sdex{Borel-Dilworth Theorem}
It is due to Harrington \site{harshemar}.  It says that if
$\poset$ is a Borel partially ordered set, then either $\poset$
is the union of countably many chains or there exist a perfect
set $P$ of pairwise incomparable elements.
One of the early Lemmas from \site{harshemar} is the following:

\begin{lemma}
Suppose $A$ is a $\Sigma^1_1$ chain in a $\Delta^1_1$ poset $\poset$.
Then there exists a $\Delta$ superset $D\supseteq A$ which is a
chain.
\end{lemma}
\proof
Suppose $\poset=(P,\leq)$ where $P$ and $\leq$ are $\Delta^1_1$.
Then
$$\theta(x,y)\;\;\rmiff\;\; x,y\in P\rmand (x\leq y\rmor y\leq x)$$
is $\Pi^1_1$ and so the result follows by Lemma \ref{overflow}.
\qed

For more on Borel linear orders, see Louveau \site{louvlin}.
Louveau \site{louvbor} is a survey paper on Borel equivalence
relations, linear orders, and partial orders.

\begin{question}
  Suppose we are given a Borel partition of the two
  element subsets of $\bairespace$.  Then is it true that
  either $\bairespace$ is the union of countably many
  $0-$homogenous sets or there exists a perfect $1-$homogeneous
  set?
\end{question}

Q.Feng \site{feng} has shown that the answer to this question is true
for open partitions.  An affirmative answer to this question
would give a proof of the Borel-Dilworth Theorem.

\sector{$\Sigma^1_1$ equivalence relations}

\begin{theorem} \sdex{$\Sigma^1_1$ equivalence relations}
  (Burgess \site{burgess}) Suppose $E$ is a $\bsig11$
  equivalence relation.  Then either $E$ has $\leq\omega_1$
  equivalence classes or there exists a perfect set of
  pairwise $E$-inequivalent reals.
\end{theorem}

\proof

We will need to prove the boundedness theorem for this result.
Define
$$WF=\{T\subseteq\finseq: T \mbox{is a well-founded tree}\}.$$
\sdex{$WF$}
For $\alpha<\omega_1$ define $WF_{<\alpha}$ to the subset of
$WF$ of all well-founded trees of rank $<\alpha$. \sdex{$WF_{<\alpha}$}
$WF$ is a complete $\Pi^1_1$ set, i.e., for every $B\subseteq\bairespace$
which is $\bpi11$ there exists a continuous map $f$ such that
$f^{-1}(WF)=B$  (see Theorem \ref{normal}).  Consequently, $WF$ is
not Borel.  On the other hand each of the $WF_{<\alpha}$ are Borel.

\begin{lemma}\label{treeborel}
  For each $\alpha<\omega_1$ the set $WF_{<\alpha}$ is Borel.
\end{lemma}
\proof
Define for $s\in\finseq$ and $\alpha<\omega_1$
$$WF^s_{<\alpha}=\{T\subseteq\finseq: T \mbox{ is a tree, }
s\in T,\;\; r_T(s)<\alpha\}.$$
The fact that $WF^s_{<\alpha}$ is Borel is
proved by induction on $\alpha$. The set of trees is $\Pi^0_1$.
For $\lambda$ a limit
$$WF^s_{<\lambda}=\Union_{\alpha<\lambda}WF^s_{<\alpha}.$$
For a successor $\alpha+1$
$$T\in WF^s_{<\alpha+1} \rmiff s\in T \rmand \forall n\;\;
( s^\concat n\in T\arrow T\in WF^{s\concat n}_{<\alpha}).$$
\qed

Another way to prove this is take a tree $T$ of rank $\alpha$ and
note that $WF_{<\alpha}=\{\hat{T}: \hat{T}\strictembed T\}$ and
this set is $\bdel11$ and hence Borel by Theorem \ref{sep}.

\begin{lemma}\label{bounded} \sdex{Boundedness Theorem}
  (Boundedness) If $A\subseteq WF$ is $\bsig11$, then
  there exists $\alpha<\omega_1$ such that $A\subseteq WF_\alpha$.
\end{lemma}
\proof
\noindent  Suppose no such $\alpha$ exists.  Then
\begin{center}
$T\in WF$ iff there exists $\hat{T}\in A$ such that $T\embed\hat{T}$.
\end{center}
But this would give a $\bsig11$ definition of $WF$, contradiction.
\qed

There is also a lightface version of the boundedness theorem, i.e.,
if $A$ is a $\Sigma^1_1$ subset of $WF$, then there exists a recursive
ordinal $\alpha<\omega_1^{CK}$ such that $A\subseteq WF_{<\alpha}$.
Otherwise,
\begin{center}
$\{e\in\omega: e $ is the code of a recursive well-founded tree $\}$
\end{center}
would be $\Sigma^1_1$.

\bigskip

Now suppose that $E$ is a $\bsig11$
  equivalence relation.
By the Normal Form Theorem \ref{normal} we know there exists
a continuous mapping $(x,y)\mapsto T_{xy}$ such that
$T_{xy}$ is always a tree and
\begin{center}
$xEy$ iff $T_{xy}\notin WF$.
\end{center}
Define
$$xE_{\alpha}y \rmiff T_{xy}\notin WF_{<\alpha}.$$
By Lemma \ref{treeborel} we know that the binary relation
$E_{\alpha}$ is Borel. Note that $E_\alpha$ refines
$E_\beta$ for $\alpha>\beta$.  Clearly,
$$E=\Intersect_{\alpha<\omega_1}E_\alpha$$
and for any limit ordinal $\lambda$
$$E_\lambda=\Intersect_{\alpha<\lambda}E_\alpha.$$

While there is no reason to expect that any of the
$E_\alpha$ are equivalence relations, we use the boundedness
theorem to show that many are.

\begin{lemma}
  For unboundedly many $\alpha<\omega_1$ the binary relation
  $E_\alpha$ is an equivalence relation.
\end{lemma}

\proof
Note that every $E_\alpha$ must be reflexive, since $E$ is reflexive and
$E=\Intersect_{\alpha<\omega_1}E_\alpha$.

The following claim
will allow us to handle symmetry.

\claim For every $\alpha<\omega_1$ there exists $\beta<\omega_1$ such
that for every $x,y$
\begin{center}
if $xE_\alpha y$ and $y\notE_\alpha x$, then
$x\notE_\beta y$.
\end{center}
\proof
Let
$$A=\{T_{xy}: xE_\alpha y \rmand y\notE_\alpha x\}.$$
Then $A$ is a Borel set.  Since $y\notE_\alpha x$ implies
$y\notE x$ and hence $x\notE y$, it follows that
$A\subseteq WF$.  By the Boundedness Theorem \ref{bounded}
there exists $\beta<\omega_1$ such that $A\subseteq WF_{<\beta}$.
\qed

The next claim is to take care of transitivity.

\claim For every $\alpha<\omega_1$ there exists $\beta<\omega_1$ such
that for every $x,y,z$
\begin{center}
if $xE_\alpha y$ and $yE_\alpha z$,
and  $x\notE_\alpha z$, then  either $x\notE_\beta y$ or
$y\notE_\beta z$.
\end{center}
\proof
Let
$$B=\{T_{xy}\oplus T_{yz}: xE_\alpha y,\;\; yE_\alpha z,
\rmand x\notE_\alpha z\}.$$
The operation $\oplus$ on a pair of trees
$T_0$ and $T_1$ is defined by
$$T_0\oplus T_1=\{(s,t): s\in T_0, \;\; t\in T_1, \rmand |s|=|t|\}.$$
Note that the rank of $T_0\oplus T_1$ is the minimum of the
rank of $T_0$ and the rank of $T_1$. (Define the rank function on
$T_0\oplus T_1$ by taking the minimum of the rank functions on the
two trees.)

The set $B$ is Borel because the relation $E_\alpha$ is.
Note also that since  $x\notE_\alpha z$ implies $x\notE z$ and
$E$ is an equivalence relation, then either
$x\notE y$ or $y\notE z$.  It follows that either
$T_{xy}\in WF$ or $T_{yz}\in WF$ and so in either case
$T_{xy}\oplus T_{yz}\in WF$ and so $B\subseteq WF$.  Again, by
the Boundedness Theorem there is a $\beta<\omega_1$ such
that $B\subseteq WF_{<\beta}$ and this proves the Claim.
\qed

Now we use the Claims
to prove the Lemma. Using  the usual Lowenheim-Skolem argument
we can find arbitrarily large countable ordinals $\lambda$ such
that for every $\alpha<\lambda$ there is a $\beta<\lambda$ which
satisfies both Claims for $\alpha$.  But this means that
$E_\lambda$ is an equivalence relation.  For suppose
$xE_\lambda y$ and $y\notE_\lambda x$.  Then since
$E_\lambda=\Intersect_{\alpha<\lambda}E_\alpha$ there must be
$\alpha<\lambda$ such that $xE_\alpha y$ and $y\notE_\alpha x$.
But by the Claim there exist $\beta<\lambda$ such that
$x\notE_\beta y$ and hence $x\notE_\lambda y$, a contradiction.
A similar argument using the second Claim works for transitivity.
\qed

Let $G$ be any generic filter over $V$ with the property that
it collapses $\omega_1$ but not $\omega_2$.  For example, Levy forcing
with finite partial functions from $\omega$ to $\omega_1$.
Then $\omega_1^{V[G]}=\omega_2^V$.  By absoluteness, $E$ is still
an equivalence relation and for any $\alpha$ if $E_\alpha$ was
an equivalence relation in $V$, then it still is one in $V[G]$.
Since
$$E_{\omega_1^V}=\Intersect_{\alpha<\omega_1^V}E_\alpha$$
and the intersection of equivalence relations is an equivalence
relation, it follows that the Borel relation $E_{\omega_1^V}$ is
an equivalence relation. So now suppose that
$E$ had more than $\omega_2$ equivalence classes in $V$. Let
$Q$ be a set of size $\omega_2$ in $V$ of pairwise E-inequivalent
reals.  Then $Q$ has cardinality  $\omega_1$ in $V[G]$ and
 for every $x\not= y\in Q$ there exists $\alpha<\omega_1^V$ with
$x\notE_\alpha y$.  Hence it must be that the elements of $Q$ are in different
$E_{\omega_1^V}$ equivalence classes.  Consequently, by Silver's Theorem
\ref{piequiv} there exists a perfect set $P$ of
 $E_{\omega_1^V}$-inequivalent reals.  Since in $V[G]$ the equivalence
 relation $E$ refines $E_{\omega_1^V}$, it must be that the elements
 of $P$ are pairwise $E$-inequivalent also.  The following is a
 $\bsig12$ statement:
$$V[G]\models \exists P \mbox{ perfect }\forall x\forall y\;\;
(x,y\in P \rmand x\not=y) \arrow x\notE y .$$
Hence, by Shoenfield Absoluteness \ref{shoenabs},
 $V$ must think that there is a perfect set of $E$-inequivalent
reals.

A way to avoid taking a generic extension of the universe is
to suppose Burgess's Theorem is false.   Then let
$M$ be the transitive collapse of an elementary substructure of
some sufficiently large $V_\kappa$ (at least large enough
to know about absoluteness and Silver's Theorem).
Let $M[G]$ be obtained as in the above proof by Levy collapsing
$\omega_1^M$.  Then we can conclude as above that $M$ thinks
$E$ has a perfect set of inequivalent elements, which contradicts
the assumption that $M$ thought Burgess's Theorem was false.
\qed

By Harrington's Theorem \ref{harring} it is consistent to
have $\bpi12$ sets of arbitrary cardinality, e.g it
is possible to have $c=\omega_{23}$ and
there exists a $\bpi12$ set $B$ with $|B|=\omega_{17}$.
Hence, if we define
$$xEy  \rmiff x=y \rmor x,y\notin B$$
then we get $\bsig12$ equivalence relation with exactly
\sdex{$?bsig12$ equivalence relation}
$\omega_{17}$ equivalence classes, but since the continuum is
$\omega_{23}$ there is no perfect set of $E$-inequivalent reals.

See Burgess \site{burgess2} \site{burgess3} for more results on
analytic equivalence relations.  For further results
concerning projective equivalence relations
see Harrington and Sami \site{harsami}, Sami \site{sami},
Stern \site{stern} \site{stern2}, Kechris \site{kechris2},
Harrington and Shelah
\site{harrshel}, Shelah \site{shelah}, and Harrington, Marker, and
Shelah \site{harshemar}.

\sector{Louveau's Theorem}

Let us define codes for Borel sets in our usual way of thinking of
them as trees with basic clopen sets attached to the terminal nodes.
\sdex{Louveau's Theorem}

Definitions
\begin{enumerate}
\item Define $(T,q)$ is an $\alpha$-code iff $T\subseteq \finseq$
\sdex{$\alpha$-code}
is a tree of rank $\leq\alpha$ and $q:\tzero\arrow\base$ is
a map from the terminal nodes, $\tzero$, of $T$ (i.e. rank zero nodes) to
a nice base, $\base$,   for the clopen sets of $\bairespace$, say
all sets of the form $[s]$ for $s\in \finseq$ plus the empty set.

\item Define $S^s(T,q)$ and $P^s(T,q)$ for $s\in T$ by induction on the
rank of $s$ as follows.  For $s\in\tzero$ define
$$P^s(T,q)=q(s) \rmand  S^s(T,q)=\comp{q(s)}.$$
For $s\in\tone$ define
$$P^s(T,q)=\Union\{ S(T,q)^{s\concat m}: {s\concat m}\in T\} \rmand
S^s(T,q)=\comp{P^s(T,q)}.$$

\item Define
$$P(T,q)=P^{\emptyseq}(T,q)\rmand S(T,q)=S^{\emptyseq}(T,q)$$
the $\Pi^0_\alpha$ set and the $\Sigma^0_\alpha$ set
coded by $(T,q)$, respectively.   ($S$ is short for
Sigma and $P$ is short for Pi.) \sdex{$P(T,q)$} \sdex{$S(T,q)$}

\item Define $C\subseteq\bairespace$ is $\Pi^0_\alpha(\hyp)$ iff
\sdex{$\Pi^0_\alpha(\hyp)$}
it has an $\alpha$-code which is hyperarithmetic.

\item $\omega_1^{CK}$ is the first nonrecursive ordinal.
\sdex{$\omega_1^{CK}$}

\end{enumerate}

\begin{theorem}\label{louv}
  (Louveau \site{louveau}) If $A,B\subseteq\bairespace$ are $\Sigma^1_1$
  sets, $\alpha<\omega_1^{CK}$, and $A$ and $B$ can be separated
  by $\bpi0{\alpha}$ set, then $A$ and $B$ can be separated by
  a $\Pi^0_\alpha(\hyp)$-set.
\end{theorem}

\begin{corollary}
  $\Delta^1_1\intersect {\bpi0\alpha}= \Pi^0_\alpha(\hyp)$
\end{corollary}

\begin{corollary}\label{section}
  (Section Problem) If $B\subseteq\bairespace\cross\bairespace$ is
  Borel and $\alpha<\omega_1$ is such that
  $B_x\in\bsig0\alpha$ for every
  $x\in\bairespace$, then
  $$B\in \bsig0\alpha(\{D\cross C: D\in\borel({\bairespace})
  \rmand C \mbox{is clopen}\}).$$
\end{corollary} \sdex{Section Problem}

Note that the converse is trivial.

This result was proved by Dellecherie for $\alpha=1$ who conjectured
it in general.  Saint-Raymond proved it for $\alpha=2$ and
Louveau and Saint-Raymond independently proved it for $\alpha=3$ and
then Louveau proved it in general.  In their paper \site{loust}
Louveau and Saint-Raymond give a different proof of it.
We will need the following lemma.

\begin{lemma} \label{defdelta}
  For $\alpha<\omega_1^{CK}$ the following sets are $\Delta_1^1$:

  $\{y:\;\; y$ is a $\beta$-code for some $\beta<\alpha\}$,

  $\{(x,y):\;\;y$  is a $\beta$-code for some $\beta<\alpha$ and $x\in P(T,q)\}$
,
  and

  $\{(x,y):\;\;y$  is a $\beta$-code for some $\beta<\alpha$ and $x\in S(T,q)\}.
$
\end{lemma}
\proof
For the first set it is enough to see that $WF_{<\alpha}$ the set of
trees of rank $<\alpha$ is $\Delta^1_1$.
Let $\hat{T}$ be a recursive tree of rank $\alpha$. Then
$T\in WF_{<\alpha}$ iff $T\strictembed \hat{T}$ shows that
$WF_{<\alpha}$ is $\Sigma^1_1$. But since $\hat{T}$ is well-founded
$T\strictembed \hat{T}$ iff $\neg(\hat{T}\embed T)$ and so it is
$\Pi^1_1$.  For the second set just use an  argument similar
to Theorem~\ref{hypcode}.  The third set is just the complement
of the second one.\qed

Now we prove Corollary \ref{section} by induction on $\alpha$.
By relativizing
the proof to a parameter we may assume $\alpha<\omega_1^{CK}$ and
that $B$ is $\Delta^1_1$.  By taking complements we may assume that
the result holds for $\Pi^0_\beta$ for all $\beta<\alpha$.
Define
$$R(x,(T,q)) \rmiff (T,q)\in \Delta^1_1(x),\;\;
(T,q) \mbox{ is an $\alpha$-code},
\rmand P(T,q)=B_x.$$
where $P(T,q)$ is the $\Pi^0_\alpha$ set coded by $(T,q)$.  Note that by the
relativized version of Louveau's Theorem for every $x$ there exists
a $(T,q)$ such that $R(x,(T,q))$.   By $\Pi^1_1$-uniformization
(Theorem~\ref{uniform})
there exist a $\Pi^1_1$ set $\hat{R}\subseteq R$ such that for every
$x$ there exists a unique $(T,q)$ such that $\hat{R}(x,(T,q))$.
 Fix $\beta<\alpha$ and
$n<\omega$ and define

$B_{\beta,n}(x,z)$ iff there exists $(T,q)\in \Delta^1_1(x)$ such that
\begin{enumerate}
\item $\hat{R}(x,(T,q))$,
\item $\rank_T(\langle n\rangle)=\beta$ and
\item $z\in P^{\langle n\rangle}(T,q)$.
\end{enumerate}
Since quantification over $\Delta^1_1(x)$ preserves $\Pi^1_1$
(Theorem \ref{gandy2}),
 $\hat{R}$ is $\Pi^1_1$, and the rest is $\Delta^1_1$ by Lemma \ref{defdelta},
 we see that $B_{\beta,n}$ is $\Pi^1_1$.
But note that
$\neg B_{\beta,n}(x,z)$ iff there exists $(T,q)\in \Delta^1_1(x)$ such that
\begin{enumerate}
\item $\hat{R}(x,(T,q))$,
\item $\rank_T(\langle n\rangle)\not=\beta$, or
\item $z\in S^{\langle n\rangle}(T,q)$.
\end{enumerate}
and consequently, $\comp{B_{\beta,n}}$ is $\Pi^1_1$ and therefore
$B_{\beta,n}$ is $\Delta_1^1$.  Note that every cross section
of $B_{\beta,n}$ is
a $\bpi0\beta$ set and so by induction (in case $\alpha>1$)
$$B_{\beta,n}\in \bpi0\alpha(\{D\cross C: D\in\borel(\bairespace)
  \rmand C \mbox{is clopen}\}).$$
But then $$B=\Union_{n<\omega,\beta<\alpha}B_{\beta,n}$$ and so
  $$B\in \bsig0\alpha(\{D\cross C: D\in\borel(\bairespace)
  \rmand C \mbox{is clopen}\}).$$

Now to do the case for $\alpha=1$, define for every
$n\in\omega$ and $s\in\finseq$

$B_{s,n}(x,z)$ iff there exists $(T,q)\in \Delta^1_1(x)$ such that
\begin{enumerate}
\item $\hat{R}(x,(T,q))$,
\item $\rank_T(\langle n\rangle)=0$,
\item $q(\langle n\rangle)=s$, and
\item $z\in [s]$.
\end{enumerate}
As in the other case $B_{s,n}$ is $\Delta^1_1$.
Let $z_0\in [s]$ be arbitrary, then define the Borel set
$C_{s,n}=\{x:(x,z_0)\in B_{s,n}\}$.
Then $B_{s,n}=C_{s,n}\cross [s]$ where
But now
$$B=\Union_{n<\omega, s\in\finseq} B_{s,n}$$
and so
$$B\in \bsig01(\{D\cross C: D\in\borel(\bairespace)
  \rmand C \mbox{ clopen}\}).$$
\qed

Note that for every $\alpha<\omega_1$ there exists a
$\bpi11$ set $U$ which is universal for all $\bdel0\alpha$
sets, i.e., every cross section of $U$ is $\bdel0\alpha$ and
every $\bdel0\alpha$ set occurs as a cross section of $U$.
To see this, let $V$ be a $\Pi^0_\alpha$ set which is universal for
$\bpi0\alpha$ sets.  Now put
$$(x,y)\in U \rmiff y\in V_{x_0} \rmand \forall z (z\in V_{x_0}\rmiff
z\notin V_{x_1}) $$
where $x=(x_0,x_1)$ is some standard pairing function.  Note also
that the complement of $U$ is also universal for all
$\bdel0\alpha$ sets, so there is a $\bsig11$ which
is universal for all $\bdel0\alpha$ sets.  Louveau's Theorem
implies that there can be no Borel set universal for all
$\bdel0\alpha$ sets.

\begin{corollary} \label{univdel}
There can be no Borel set universal for all
$\bdel0\alpha$ sets.
\end{corollary}\sdex{$?bdel0\alpha$-universal set}

In order to prove this corollary we will need the following lemmas.
A space is Polish iff it is a separable complete metric space.

\begin{lemma}\label{sierchar}
  If $X$ is a $0$-dimensional
  Polish space, then there exists a closed set $Y\subseteq\bairespace$
  such that $X$ and $Y$ are homeomorphic.
\end{lemma}
\proof
Build a tree $\langle C_s:s\in T\rangle$ of nonempty clopen sets
indexed by a tree $T\subseteq\finseq$ such that
\begin{enumerate}
  \item $C_{\emptyseq}=X$,
  \item the diameter of $C_s$ is less that $1/|s|$ for $s\not=\emptyseq$, and
  \item for each $s\in T$ the clopen set $C_s$ is the disjoint union
  of the clopen sets $$\{C_{s\concat n}:s\concat n\in T\}.$$
\end{enumerate}
If $Y=[T]$ (the infinite branch of $T$), then $X$ and $Y$ are homeomorphic.
\qed

I am not sure who proved this first.
I think the argument for
the next lemma comes
from a theorem about Hausdorff that lifts  the difference hierarchy
on the $\bdel02$-sets to the $\bdel0\alpha$-sets.  This
presentation is taken from
Kechris \site{kechris} {\it mutatis mutandis}.\footnote{
Latin for plagiarized.}

\begin{lemma} \label{haushomeo}
 For any sequence
  $\langle B_n:n\in\omega\rangle$
  of Borel subsets of $\bairespace$ there exists $0$-dimensional Polish
  topology,  $\tau$, which contains the standard topology and each
  $B_n$ is a clopen set in $\tau$.
\end{lemma}
\proof

This will follow easily from the next two claims.

\claim Suppose $(X,\tau)$ is a 0-dimensional Polish space and
$F\subseteq X$ is
closed,  then there exists a 0-dimensional Polish topology
$\sigma\supseteq \tau$ such that $F$ is clopen in $(X,\sigma)$.
(In fact, $\tau\union\{F\}$ is a subbase for $\sigma$.)
\proof
Let $X_0$ be $F$ with the subspace topology given by $\tau$
and $X_1$ be $\comp{F}$ with the subspace topology.  Since $F$ is
closed $X_0$ the complete metric on $X$ is complete on $X_0$.
Since $\comp{F}$ is open there is another metric which is complete
on $X_1$.  This is a special case of a Alexandroff's Theorem
which says that a $G_\delta$ set in a completely metrizable space is
completely metrizable in the subspace topology.  In this case the
complete metric $\hat{d}$ on $\comp{F}$ would be defined by
$$\hat{d}(x,y)=d(x,y)+\left|{1\over d(x,F)}- {1\over d(y,F)}\right|$$
where $d$ is a complete metric on $X$ and $d(x,F)$ is the distance
from $x$ to the closed set $F$.

Let
$$(X,\sigma)=X_0\oplus X_1$$
be the discrete topological sum, i.e.,
$U$ is open iff $U=U_0\union U_1$ where $U_0\subseteq X_0$ is open in $X_0$
and $U_1\subseteq X_1$ is open in $X_1$.
\qed

\claim If $(X,\tau)$ is a Hausdorff space and $(X,\tau_n)$ for
$n\in\omega$ are 0-dimensional Polish topologies extending $\tau$, then
there exists a 0-dimensional Polish topology $(X,\sigma)$ such that
$\tau_n\subseteq \sigma$ for every $n$.  (In fact
$\Union_{n<\omega}\tau_n$ is a subbase for $\sigma$.)
\proof
Consider the 0-dimensional Polish space
$$\prod_{n\in\omega}(X,\tau_n).$$
Let $f:X\arrow \prod_{n\in\omega}(X,\tau_n)$ be the embedding
which takes each $x\in X$ the the constant sequence $x$
(i.e., $f(x)=\langle x_n:n\in\omega\rangle$ where $x_n=x$ for every $n$).
Let $D\subseteq\prod_{n\in\omega}(X,\tau_n)$ be the range of $f$, the
set of constant sequences. Note that $f:(X,\tau)\arrow (D,\tau)$ is
a homeomorphism. Let $\sigma$ be the topology on $X$ defined by
\begin{center}
$U\in \sigma$ iff there exists $V$ open in $\prod_{n\in\omega}(X,\tau_n)$ with
$U=f^{-1}(V)$.
\end{center}
Since each $\tau_n$ extends $\tau$ we get that $D$ is a closed subset
of $\prod_{n\in\omega}(X,\tau_n)$.  Consequently, $D$ with the subspace
topology inherited from $\prod_{n\in\omega}(X,\tau_n)$ is Polish.
It follows that $\sigma$ is a Polish topology on $X$.
To see that $\tau_n\subseteq\sigma$ for every $n$ let
$U\in\tau_N$ and define
$$V=\prod_{n<N}X\cross U\cross\prod_{n>N}X.$$
Then $f^{-1}(V)=U$ and so $U\in\sigma$.\qed

We prove Lemma \ref{haushomeo} by induction on the rank of the
Borel sets.  Note that by the second Claim it
is enough to prove it for one Borel set at a time.  So suppose
$B$ is a $\bsig0\alpha$ subset of $(X,\tau)$. Let
$B=\Union_{n\in\omega}B_n$ where each $B_n$ is $\bpi0\beta$ for
some $\beta<\alpha$.  By induction on $\alpha$ there exists
a 0-dimensional Polish topology $\tau_n$ extending $\tau$ in
which each $B_n$ is clopen.   Applying the second Claim gives
us a 0-dimensional topology $\sigma$ extending $\tau$ in which each
$B_n$ is clopen and therefore $B$ is open.  Apply the first Claim
to get a 0-dimensional Polish topology in which $B$ is clopen.
\qed

\proof
(of Corollary \ref{univdel}).  The idea of this proof is to reduce it to
the case of a $\bdel0\alpha$ set universal for $\bdel0\alpha$-
sets, which is easily seen to be impossible by the standard diagonal
argument.

Suppose $B$ is a Borel set which is universal
for all $\bdel0\alpha$ sets.  Then by the Corollary~\ref{section}
$$B\in \Delta^0_\alpha(\{D\cross C: D\in\borel(\bairespace)
  \rmand C \mbox{ is clopen}\}).$$
By Lemma \ref{haushomeo} there exists a 0-dimensional Polish topology $\tau$
such that if $X=(\bairespace,\tau)$, then
$B$ is $\bdel0\alpha(X\cross\bairespace)$.  Now by Lemma \ref{sierchar}
there exists a closed set $Y\subseteq\bairespace$ and a homeomorphism
$h:X\arrow Y$. Consider
$$C=\{(x,y)\in X\cross X: (x,h(y))\in B\}.$$
The set $C$ is $\bdel0\alpha$ in $X\cross X$ because it is
the continuous preimage of the set $B$ under the map $(x,y)\mapsto (x,h(y))$.
The set $C$ is also universal for $\bdel0\alpha$ subsets of
$X$ because the set $Y$ is closed.  To see this for $\alpha>1$ if
$H\in\bdel0\alpha(Y)$, then $H\in\bdel0\alpha(\bairespace)$,
consequently there exists $x\in X$ with $B_x=H$.  For $\alpha=1$ just
use  that disjoint closed subsets of $\bairespace$ can be separated by
clopen sets.

Finally, the set $C$ gives a contradiction by the usual
diagonal argument:
$$D=\{(x,x): x\notin C\}$$
would be  $\bdel0\alpha$ in $X$ but cannot be a cross section of $C$.
\qed

\begin{question}
  (Mauldin) Does there exists a $\Pi^1_1$ set which is universal
  for all $\bpi11$ sets which are not Borel?
\end{question}

We could also ask for the complexity of a set which is universal
for $\bsig0\alpha\setminus \bdel0\alpha$ sets.

\sector{Proof of Louveau's Theorem}

Finally, we arrive at our last section.
The following summarizes how I feel now.

\begin{quote}
  You are walking down the street minding your own business and someone
  stops you and asks directions.  Where's xxx hall?  You don't
  know and you say you don't know.  Then they point at the next street and
  say: Is that xxx street?   Well by this time you feel kind of
  stupid so you say,  yea yea that's xxx street,   even though you haven't
  got the slightest idea whether it is or not.  After all, who wants
  to admit they don't know where they are going or where they are.
\end{quote}

\bigskip

For $\alpha<\omega_1^{CK}$ define $D\subseteq\bairespace$
is $\Sigma^0_\alpha(\semihyp)$  \sdex{$\Sigma^0_\alpha(\semihyp)$}
iff there exists $S$ a $\Pi_1^1$ set of hyperarithmetic reals
such that every
element of $S$ is a $\beta$-code for some $\beta<\alpha$ and
$$D=\Union\{P(T,q): (T,q)\in S\}.$$
A set is $\Pi^0_\alpha(\semihyp)$ iff it is the complement of
\sdex{$\Pi^0_\alpha(\semihyp)$}
a $\Sigma^0_\alpha(\semihyp)$ set.  The
$\Pi^0_0(\semihyp)$ sets are just the usual clopen basis
($[s]$ for $s\in\finseq$ together with
the empty set) and $\Sigma^0_0(\semihyp)$ sets are their complements.

\begin{lemma}
 $\Sigma^0_\alpha(\semihyp)$ sets are $\Pi^1_1$ and
 consequently $\Pi^0_\alpha(\semihyp)$ sets are $\Sigma^1_1$.
\end{lemma}
\proof
$x\in \Union\{P(T,q): (T,q)\in S\}$ iff there exists $(T,q)\in \Delta^1_1$
such that
$(T,q)\in S$ and $x\in P(T,q)$.
Quantification over
$\Delta_1^1$ preserves $\Pi^1_1$ ( see Corollary \ref{gandy2} )
and Lemma \ref{defdelta} implies that ``$x\in P(T,q)$'' is $\Delta^1_1$.
\qed

We will need the following reflection principle in order to prove
the Main Lemma \ref{louvmain}.

A predicate $\Phi\subseteq P(\omega)$ is called $\Pi_1^1$ on
$\Pi_1^1$ iff for any $\Pi^1_1$ set $N\subseteq\omega\cross\omega$
the set $\{e:\Phi(N_e)\}$ is $\Pi_1^1$
(where $N_e=\{n:(e,n)\in N\}$).  \sdex{$\Pi_1^1$ on $\Pi_1^1$}

\begin{lemma}\label{reflect} \sdex{$\Pi^1_1$-Reflection}
  (Harrington \site{harshemar} Kechris \site{kecref})
  $\Pi^1_1$-Reflection.   Suppose
  $\Phi(X)$  is $\Pi_1^1$ on
  $\Pi_1^1$ and Q is a $\Pi^1_1$ set.

  If $\Phi(Q)$, then
  there exists a $\Delta^1_1$ set $D\subseteq Q$ such that $\Phi(D)$.
\end{lemma}
\proof
By the normal form theorem \ref{normal} there is a recursive mapping
$e\arrow T_e$ such that $e\in Q$ iff $T_e$ is well-founded.
Define for $e\in\omega$
$$N_e^0=\{\hat{e}: T_{\hat{e}}\embed T_e\}$$
$$N_e^1=\{\hat{e}: \neg(T_e\strictembed T_{\hat{e}})\}$$
then  $N^0$ is $\Sigma^1_1$ and $N^1$ is $\Pi^1_1$.  For
$e\in Q$ we have $N_e^0=N_e^1=D_e\subseteq Q$ is  $\Delta^1_1$;
and for $e\notin Q$ we have that $N_e^1=Q$.  If we assume for
contradiction that $\neg\Phi(N_e^1)$ for all $e\in Q$, then
$$e\notin Q \rmiff \phi(N_e^1).$$
But this would mean that $Q$ is $\Delta^1_1$ and this proves
the Lemma.\qed

Note that a $\Pi^1_1$ predicate need not be $\Pi_1^1$ on
$\Pi_1^1$ since the predicate
$$\Phi(X)=\mbox{``}0\notin X\mbox{''}$$
is $\Delta^0_0$ but not $\Pi_1^1$ on
$\Pi_1^1$.  Some examples of $\Pi^1_1$ on $\Pi^1_1$ predicates
$\Phi(X)$ are
$$\Phi(X)\rmiff \forall x\notin X\;\; \theta(x)$$
or
$$\Phi(X)\rmiff \forall x,y\notin X\;\; \theta(x,y)$$
where $\theta$ is a $\Pi^1_1$ sentence.

\begin{lemma}\label{louvmain}
  Suppose $A$ is $\Sigma^1_1$ and $A\subseteq B\in\Sigma^0_\alpha(\semihyp)$,
  then there exists $C\in \Sigma^0_\alpha(\hyp)$ with
  $A\subseteq C\subseteq B$.
\end{lemma}
\proof
Let $B=\Union\{P(T,q): (T,q)\in S\}$ where $S$ is a $\Pi^1_1$ set of
hyperarithmetic $<\alpha$-codes.  Let $\hat{S}\subseteq\omega$
be the $\Pi^1_1$ set of $\Delta_1^1$-codes for elements of $S$, i.e.
\begin{center}
$e\in\hat{S}$ iff $e$ is a $\Delta^1_1$-code for  $(T_e,q_e)$ and
$(T_e,q_e)\in S$.
\end{center}
Now define the predicate $\Phi(X)$ for $X\subseteq\omega$ as follows:
\begin{center}
$\Phi(X)$ iff $X\subseteq \hat{S}$ and $A\subseteq\Union_{e\in X}P(T_e,q_e)$.
\end{center}

The predicate $\Phi(X)$ is $\Pi^1_1$ on $\Pi_1^1$ and
$\Phi(\hat{S})$.  Therefore by reflection (Lemma \ref{reflect})
there exists a $\Delta^1_1$ set $D\subseteq \hat{S}$ such that
$\Phi(D)$.   Define
$(T,q)$ by
$$T=\{e\concat s: e\in D \rmand s\in T_e\}\;\;\;\;\;\;\;
q(e\concat s)=q_e(s) \mbox{ for } e\in D \rmand s\in \tzero_e.$$
Since $D$ is $\Delta^1_1$ it is easy to check that $(T,q)$ is
$\Delta^1_1$ and hence hyperarithmetic.   Since $\Phi(D)$ holds
it follows that $C=S(T,q)$ the $\Sigma^0_\alpha(\hyp)$ set coded
by $(T,q)$ has the property that $A\subseteq C$ and since $D\subseteq\hat{S}$
it follows that $C\subseteq B$.
\qed

Define for $\alpha<\omega_1^{CK}$ the $\alpha$-topology by taking for
basic open sets the family                       \sdex{$\alpha$-topology}
$$\Union\{\Pi^0_\beta(\semihyp):\beta<\alpha\}.$$
As usual, $\cl_\alpha(A)$ denotes the closure of the set $A$ in the
$\alpha$-topology. \sdex{$?cl_\alpha(A)$}

The 1-topology is just the standard topology on $\bairespace$.  The
$\alpha$-topology has its basis certain special $\Sigma^1_1$ sets
so it is intermediate between the standard topology and the Gandy
topology corresponding to Gandy forcing.

\begin{lemma}\label{semihyp}
  If $A$ is $\Sigma^1_1$, then $\cl_\alpha(A)$ is $\Pi^0_\alpha(\semihyp)$.
\end{lemma}
\proof
Since the $\Sigma_{\beta}^0(\semihyp)$ sets for $\beta<\alpha$ form
a basis for the $\alpha$-closed sets,
$$\cl_\alpha(A)=\Intersect\{X\supseteq A:\exists\beta<\alpha\;\;
X\in\Sigma_{\beta}^0(\semihyp)\}.$$
By Lemma \ref{louvmain} this same intersection can be written:
$$\cl_\alpha(A)=\Intersect\{X\supseteq A:\exists\beta<\alpha\;\;
X\in\Sigma_{\beta}^0(\hyp)\}.$$
But now define $(T,q)\in Q$ iff $(T,q)\in \Delta^1_1$, $(T,q)$
is a $\beta$-code for some $\beta<\alpha$, and
$A\subseteq S(T,q)$.  Note that $Q$ is a $\Pi^1_1$ set and consequently,
$\cl_\alpha(A)$ is a $\Pi^0_\alpha(\semihyp)$ set, as desired.
\qed

Note that it follows from the Lemmas that for $A$ a  $\Sigma^1_1$
set, $\cl_\alpha(A)$ is a $\Sigma^1_1$ set which is a basic open
set in the $\beta$-topology for any $\beta>\alpha$.

Let $\poset$ be Gandy forcing, i.e., the partial order of all
nonempty $\Sigma^1_1$ subsets of $\bairespace$ and let
$\name{a}$ be a name for the real obtained by forcing with
$\poset$, so that by  Lemma \ref{real}, for any $G$ which is
$\poset$-generic, we have that $p\in G$ iff $a^G\in p$.

\begin{lemma} \label{louvforce}
  For any $\alpha<\omega_1^{CK}$, $p\in\poset$, and $C\in{\bpi0\alpha}$
  (coded in $V$) if
  $$p\forces \name{a}\in C,$$
  then
  $$cl_\alpha(p)\forces \name{a}\in C.$$
\end{lemma}
\proof
This is proved by induction on $\alpha$.

For $\alpha=1$ recall that the $\alpha$-topology is the standard topology
and $C$ is a standard closed set. If $p\forces \name{a}\in C$, then it
better be that $p\subseteq C$, else there exists $s\in\finseq$ with
$q=p\intersect [s]$ nonempty and $[s]\intersect C=\emptyset$. But
then $q\leq p$ and $q\forces\name{a}\notin C$.  Hence $p\subseteq C$
and since $C$ is closed, $\cl(p)\subseteq C$.  Since
$\cl(p)\forces a\in\cl(p)$, it follows that
$\cl(p)\forces a\in C$.

For $\alpha>1$ let
$$C=\Intersect_{n<\omega} \comp{C_n}$$
where each $C_n$ is
${\bpi0{\beta}}$ for some $\beta<\alpha$.  Suppose for contradiction
that
$$cl_\alpha(p)\not\forces \name{a}\in C$$
 Then for some $n<\omega$ and
$r\leq cl_\alpha(p)$ it must be that
$$r\forces \name{a}\in C_n.$$
Suppose that $C_n$ is ${\bpi0{\beta}}$ for some $\beta<\alpha$.  Then
by induction
$$\cl_\beta(r)\forces \name{a}\in C_n.$$
But $\cl_\beta(r)$ is a $\Pi^0_\beta(\semihyp)$ set by Lemma \ref{semihyp}
and hence a basic open set in the $\alpha$-topology.  Note that since
they force contradictory information
($\cl_\beta(r)\forces \name{a}\notin C$ and $p\forces \name{a}\in C$)
it must be that $\cl_\beta(r)\intersect p=\emptyset$, (otherwise
the two conditions would be compatible in $\poset$).
But since $\cl_\beta(r)$ is $\alpha$-open this means that
$$\cl_\beta(r)\intersect\cl_\alpha(p)=\emptyset$$
which contradicts the
fact that $r\leq \cl_\alpha(p)$.
\qed

Now we are ready to prove Louveau's  Theorem \ref{louv}.
Suppose $A$ and $B$ are $\Sigma^1_1$ sets and $C$ is
a $\bpi0\alpha$ set with $A\subseteq C$ and
$C\intersect B=\emptyset$.  Since $A\subseteq C$ it follows
that
$$A\forces \name{a}\in C.$$
By Lemma \ref{louvforce} it follows that
$$\cl_\alpha(A)\forces \name{a}\in C.$$
Now it must be that $\cl_\alpha(A)\intersect B=\emptyset$, otherwise
letting $p=\cl_\alpha(A)\intersect B$ would be a condition of
$\poset$ such that
$$p\forces \name{a}\in C$$
and
$$p\forces\name{a}\in B$$
which would imply that $B\intersect C\not=\emptyset$ in the generic
extension.  But by absoluteness $B$ and $C$ must remain disjoint.
So $\cl_\alpha(A)$ is a $\Pi_\alpha(\semihyp)$-set (Lemma \ref{semihyp})
which is disjoint from the set $B$ and thus by applying Lemma \ref{louvmain}
to its complement there exists a $\Pi^0_\alpha(\hyp)$-set $C$
with $\cl_\alpha(A)\subseteq C$ and $C\intersect B=\emptyset$.
\qed

The argument presented here is partially from Harrington \site{harlouv},
but contains even more simplification brought about by using
forcing and absoluteness.   Louveau's Theorem is also
proved in Sacks \site{sack} and Mansfield and Weitkamp \site{manweit}.

%%%%%%%% references and index %%%%%%%%%%%
\newpage
\addtocontents{toc}{ }
\addtocontents{toc}{\protect\bigskip}
\addtocontents{toc}{References \hfill \thepage}

\twocolumn
\pagestyle{myheadings}
\markboth{Index}{Index}
\markright{Index}
\addtocontents{toc}{ }
\addtocontents{toc}{Index \hfill  \thepage}
           \indexentry{$2^{\omega }$}{6}
            \indexentry{accumulation points}{62}
          \indexentry{$\alpha $-code}{131}
          \indexentry{$\alpha $-forcing}{26}
          \indexentry{$\alpha $-topology}{140}
            \indexentry{Aronszajn tree}{54}
          \indexentry{$[A]_I$}{36}
           \indexentry{$A_\alpha $}{55}
           \indexentry{$A_{<\alpha }$}{55}
            \indexentry{Baire space}{6}
            \indexentry{Borel metric space}{120}
            \indexentry{Borel-Dilworth Theorem}{124}
     \indexentry{${\prm Borel}(F)$}{14}
     \indexentry{${\prm Borel}(X)$}{9}
     \indexentry{${\prm Borel}(X)/{\prm meager}(X)$}{50}
            \indexentry{Boundedness Theorem}{126}
     \indexentry{$\sim {B}\discretionary {-}{}{}$}{9}
            \indexentry{Cantor space}{6}
            \indexentry{cBa}{31}
            \indexentry{characteristic function}{17}
          \indexentry{$\cl_\alpha (A)$}{140}
            \indexentry{code for a hyperarithmetic set}{106}
            \indexentry{collinear points}{122}
     \indexentry{${\prm cov}(I)$}{75}
     \indexentry{${\prm cov}({\prm meager}(2^{\omega }))$}{75}
         \indexentry{$\bdel0\alpha $-universal set}{134}
         \indexentry{$\bdel0\alpha $}{9}
         \indexentry{$\bdel0\alpha (F)$}{14}
          \indexentry{$\Delta ^1_1$-codes}{111}
          \indexentry{$\Delta ^1_2$ well-ordering}{84}
          \indexentry{$\Delta _0$-formulas}{86}
            \indexentry{direct sum}{47}
           \indexentry{$F/I$}{36}
            \indexentry{field of sets}{14}
     \indexentry{${\prm FIN}(\kappa ,2)$}{51}
            \indexentry{Fusion}{54}
           \indexentry{$F_\sigma $}{10}
            \indexentry{Gandy forcing}{115}
           \indexentry{$G_\delta $}{10}
          \indexentry{$\hat{Q}_\alpha $}{62}
            \indexentry{HC}{86}
            \indexentry{hereditarily countable sets}{86}
            \indexentry{hereditary order}{61}
            \indexentry{hyperarithmetic sets}{106}
           \indexentry{$I$-Luzin set}{38}
          \indexentry{$\kappa $-Borel}{104}
          \indexentry{$\kappa $-Souslin}{104}
          \indexentry{$\kappa $-Souslin}{88}
            \indexentry{Kleene Separation Theorem}{106}
            \indexentry{Louveau's Theorem}{131}
           \indexentry{$L_{\infty }(P_\alpha : \alpha <\kappa )$}{29}
            \indexentry{MA$_\kappa $(ctbl)}{75}
            \indexentry{Mansfield-Solovay Theorem}{90}
            \indexentry{Martin's Axiom}{20}
            \indexentry{Martin-Solovay Theorem}{41}
     \indexentry{${\prm MA}_\kappa $}{41}
     \indexentry{${\prm meager}(X)$}{51}
            \indexentry{Mostowski's Absoluteness}{82}
          \indexentry{$\mu $}{67}
           \indexentry{$m_\poset$}{78}
            \indexentry{nice $\alpha $-tree}{25}
            \indexentry{Normal form for $\Sigma ^1_1$}{82}
         \indexentry{$[\omega ]^{\omega }$}{7}
          \indexentry{$\omega ^{<\omega } $}{6}
          \indexentry{$\omega ^{\omega }$}{6}
          \indexentry{$\omega _1^{CK}$}{131}
         \indexentry{$\bool^+$}{33}
          \indexentry{$\ord(\bool)$}{32}
     \indexentry{${\prm ord}(X)$}{12}
           \indexentry{$P(T,q)$}{131}
            \indexentry{perfect set forcing}{74}
            \indexentry{perfect set}{13}
            \indexentry{perfect tree}{54}
          \indexentry{$\Pi ^0_\alpha ({\prm hyp})$}{131}
          \indexentry{$\Pi ^0_\alpha ({\prm semihyp})$}{138}
          \indexentry{$\Pi ^1_1$ equivalence relations}{115}
          \indexentry{$\Pi ^1_1$ singleton}{93}
          \indexentry{$\Pi ^1_1$-Reduction}{109}
          \indexentry{$\Pi ^1_1$-Reflection}{138}
          \indexentry{$\Pi _1^0$}{80}
          \indexentry{$\Pi _1^1$ on $\Pi _1^1$}{138}
          \indexentry{$\Pi _\beta $-sentence}{29}
         \indexentry{$\bpi0\alpha $}{9}
         \indexentry{$\bpi0\alpha (F)$}{14}
         \indexentry{$\bpi11$}{80}
         \indexentry{$\bpi12$}{100}
          \indexentry{$\poset * \name{\posetq}$}{43}
          \indexentry{$\poset(T)$}{56}
          \indexentry{$\posetq_\alpha $}{61}
          \indexentry{$\poset_\alpha$}{45}
            \indexentry{prewellorderings}{93}
            \indexentry{prewellordering}{109}
    \indexentry{${\goth p}$}{21}
            \indexentry{Q-sets}{21}
           \indexentry{$Q_\alpha $}{51}
            \indexentry{rank function}{25}
     \indexentry{${\prm rank}(p)$}{45}
           \indexentry{$S(T,q)$}{131}
            \indexentry{scale property}{92}
            \indexentry{second countable}{6}
            \indexentry{Section Problem}{132}
            \indexentry{separable}{6}
            \indexentry{separative}{31}
            \indexentry{Shoenfield Absoluteness}{89}
            \indexentry{Sierpinski set}{67}
         \indexentry{$\bsig0\alpha $}{9}
         \indexentry{$\bsig0\alpha (F)$}{14}
         \indexentry{$\bsig11$}{80}
         \indexentry{$\bsig12$ equivalence relation}{130}
          \indexentry{$\sigma $-field}{36}
          \indexentry{$\sigma $-ideal}{36}
          \indexentry{$\sigma $-ring}{11}
          \indexentry{$\Sigma ^0_\alpha ({\prm semihyp})$}{138}
          \indexentry{$\Sigma ^1_1$ equivalence relations}{126}
          \indexentry{$\Sigma ^1_1$}{80}
          \indexentry{$\Sigma ^1_1(x)$}{80}
          \indexentry{$\Sigma ^1_2$}{84}
          \indexentry{$\Sigma ^1_{2}=\Sigma ^{HC}_1$}{86}
          \indexentry{$\Sigma _1$-formula}{86}
            \indexentry{Silver forcing}{20}
            \indexentry{Souslin-Luzin Separation}{104}
            \indexentry{Spector-Gandy Theorem}{111}
            \indexentry{splitting node}{54}
     \indexentry{${\prm stone}(B)$}{36}
          \indexentry{$\sum _{\alpha <\omega _1}\poset_\alpha $}{47}
         \indexentry{$(\sum _{\alpha <\omega _1}\poset_\alpha ) *\name{\posetq}$}{47}
            \indexentry{super Luzin set}{51}
            \indexentry{super-$I$-Luzin}{38}
            \indexentry{switcheroo}{28}
           \indexentry{$s\concat n$}{6}
           \indexentry{$s\subset t$}{94}
          \indexentry{$[s]$}{6}
          \indexentry{$|s|$}{6}
            \indexentry{tree embedding}{94}
            \indexentry{tree}{25}
            \indexentry{two step iteration}{42}
           \indexentry{$T\embed\hat{T}$}{94}
           \indexentry{$T\leq _n T^\prime $}{54}
           \indexentry{$T\strictembed\mathaccent "705E {T}$}{94}
          \indexentry{$[T]$}{25}
           \indexentry{$T^0$}{26}
           \indexentry{$T^{>0}$}{26}
           \indexentry{$T_\alpha $}{51}
            \indexentry{uniformization property}{92}
            \indexentry{universal for $\bsig11$ sets}{81}
            \indexentry{universal set}{11}
            \indexentry{universal}{11}
            \indexentry{V=L}{84}
            \indexentry{well-founded}{25}
           \indexentry{$WF$}{126}
           \indexentry{$WF_{<\alpha }$}{126}
           \indexentry{$WO$}{115}
           \indexentry{$x<_cy$}{84}
            \indexentry{ZFC$^*$}{83}

\onecolumn
\pagestyle{myheadings}
\markboth{Elephant Sandwiches}{Elephant Sandwiches}
\markright{Elephant Sandwiches}
\addtocontents{toc}{ }
\addtocontents{toc}{Elephant Sandwiches \hfill  \thepage}

\newpage

A man walks by a restaurant.  Splashed all over are signs saying
``Order any sandwich'',  ``Just ask us, we have it'', and ``All
kinds of sandwiches''.

Intrigued, he walks in and says to the
proprietor, ``I would like an elephant sandwich.''

The proprietor
responds ``Sorry, but you can't have an elephant sandwich.''

``What do you mean?'' says the man, ``All your signs say to order
any sandwich and here the first thing I ask for you don't
have.''

Says the proprietor ``Oh we have elephant, its just that
here it is 5pm already and I just don't want to start another
elephant.''

\addtocontents{toc}{ }

\end{document}